\documentclass[11pt]{amsart}
\usepackage{txfonts}
\usepackage{amscd,amssymb,mathabx,mathtools}
\usepackage{graphicx,calrsfs}
\usepackage[colorlinks,plainpages,backref=page,urlcolor=blue,hyperindex,breaklinks]{hyperref}
\usepackage{tikz}
\usetikzlibrary{cd}
\usepackage{mdwlist}
\usepackage[shortlabels]{enumitem}
\usepackage{relsize}
\usepackage{xcolor}

\usepackage[capitalise]{cleveref}
\crefformat{section}{\S#2#1#3} 
\crefformat{subsection}{\S#2#1#3}
\crefformat{subsubsection}{\S#2#1#3}

\numberwithin{equation}{section}

\newtheorem{theorem}[equation]{Theorem}
\newtheorem{corollary}[equation]{Corollary}
\newtheorem{lemma}[equation]{Lemma}
\newtheorem{proposition}[equation]{Proposition}
\newtheorem{proposition-definition}[equation]{Proposition-Definition}

\theoremstyle{definition}
\newtheorem{definition}[equation]{Definition}
\newtheorem{notation}[equation]{Notation}
\newtheorem{convention}[equation]{Convention}
\newtheorem{example}[equation]{Example}
\newtheorem*{example*}{Example} %unnumbered version
\newtheorem{remark}[equation]{Remark}

\newtheorem*{remark*}{Remark}

\makeatletter
\newtheorem*{rep@theorem}{\rep@title}
\newcommand{\newreptheorem}[2]{%
\newenvironment{rep#1}[1]{%
 \def\rep@title{#2 \ref{##1}}%
 \begin{rep@theorem}}%
 {\end{rep@theorem}}}
\makeatother

\newreptheorem{theorem}{Theorem}
\newreptheorem{corollary}{Corollary}

\newcommand{\N}{\mathbb{N}}
\newcommand{\Z}{\mathbb{Z}}
\newcommand{\Q}{\mathbb{Q}}

\newcommand{\C}{\mathbb{C}}
\newcommand{\F}{\mathbb{F}}

\newcommand{\D}{\mathbb{D}}

\newcommand{\LL}{\mathfrak {L}}

\DeclareMathOperator{\res}{res}
\DeclareMathOperator{\gr}{gr}

\DeclareMathOperator{\id}{id}
\DeclareMathOperator{\ab}{{ab}}
\DeclareMathOperator{\abf}{{abf}}

\DeclareMathOperator{\Der}{Der}

\DeclareMathOperator{\End}{{End}}

\DeclareMathOperator{\proj}{pr}

\DeclareMathOperator{\ad}{ad}
\DeclareMathOperator{\Aut}{Aut}
\DeclareMathOperator{\Out}{Out}

\DeclareMathOperator{\pt}{pt}

\DeclareMathOperator{\Tors}{Tors}
\DeclareMathOperator{\BS}{BS}
\DeclareMathOperator{\Conf}{Conf}
\DeclareMathOperator{\Homeo}{Homeo}
\DeclareMathOperator{\MCG}{MCG}
\DeclareMathOperator{\rat}{\scalebox{0.6}{$\Q$}}
\DeclareMathOperator{\zz}{\scalebox{0.6}{$\Z$}}

\DeclareMathOperator{\ii}{i}

\newcommand{\surj}{\twoheadrightarrow}
\newcommand{\inj}{\hookrightarrow}
\newcommand{\isom}{\xrightarrow{
   \,\smash{\raisebox{-0.3ex}{\ensuremath{\scriptstyle\simeq}}}\,}}

\newcommand{\abs}[1]{\left| #1 \right|}

\newcommand{\ssqrt}{\!\sqrt}

\newcommand{\Sym}{\mathfrak S}

\newcommand{\qt}[1][q]{\mathsmaller(#1\mathsmaller)} %q-torsion
\newcommand{\pr}[1][p]{\mathsmaller[#1\mathsmaller]} %p-restricted

\makeatletter %define llangle and rrangle
\newsavebox{\@brx}
\newcommand{\llangle}[1][]{\savebox{\@brx}{\(\m@th{#1\langle}\)}%
  \mathopen{\copy\@brx\mkern2mu\kern-0.9\wd\@brx\usebox{\@brx}}}
\newcommand{\rrangle}[1][]{\savebox{\@brx}{\(\m@th{#1\rangle}\)}%
  \mathclose{\copy\@brx\mkern2mu\kern-0.9\wd\@brx\usebox{\@brx}}}
\makeatother

 \SetEnumitemKey{widestlabel}
  {labelwidth = \widthof{\ref{enum-\EnumitemId}},
   after      = \label{enum-\EnumitemId}}
   
\definecolor{lime}{HTML}{A6CE39}
\DeclareRobustCommand{\orcidicon}{
	\begin{tikzpicture}
	\draw[lime, fill=lime] (0,0) 
	circle [radius=0.16] 
	node[white] {{\fontfamily{qag}\selectfont \tiny ID}};
	\draw[white, fill=white] (-0.0625,0.095) 
	circle [radius=0.007];
	\end{tikzpicture}
	\hspace{-2mm}
}
\foreach \x in {A, ..., Z}{\expandafter\xdef\csname orcid\x\endcsname{\noexpand\href{%
https://orcid.org/\csname orcidauthor\x\endcsname}{\noexpand\orcidicon}}}

\makeatletter
\@namedef{subjclassname@2020}{%
  \textup{2020} Mathematics Subject Classification}
\makeatother

\begin{document}
%\date{June 25, 2021}
%\date{March 13, 2026}

\title[Lower central series and split extensions]%
{Lower central series and split extensions} 

\author[Jacques~Darn\'e]{Jacques~Darn\'e$^1$%\orcidB{}
}
\address{LAMFA, Université de Picardie Jules Verne,
33 rue Saint-Leu,
80039 Amiens Cedex 1, France}
\email{\href{mailto:jacques.darne@normalesup.org}{jacques.darne@normalesup.org}}
\urladdr{\href{https://sites.google.com/view/jacques-darne/}%
{sites.google.com/view/jacques-darne/}}
\thanks{$^1$Partially supported by the ANR Projects ChroK ANR-16-CE40-0003 
and AlMaRe ANR-19-CE40-0001-01}

\author[Alexander~I.~Suciu]{Alexander~I.~Suciu$^2$\orcidA{}}
\address{Department of Mathematics,
Northeastern University,
Boston, MA 02115, USA}
\email{\href{mailto:a.suciu@northeastern.edu}{a.suciu@northeastern.edu}}
\urladdr{\href{https://suciu.sites.northeastern.edu}%
{suciu.sites.northeastern.edu}}
\thanks{$^2$Partially supported by the Simons Foundation Collaboration 
Grants for Mathematicians \#354156 and \#693825 and by the project 
``Singularities and Applications" - CF 132/31.07.2023 funded by the 
European Union - NextGenerationEU - through Romania’s National Recovery 
and Resilience Plan.}



\begin{abstract}
Following Lazard, we study the $N$-series of a group $G$ and their 
associated graded Lie algebras. The main examples we consider 
are the lower central series (LCS), Stallings' rational and mod-$q$ 
versions, and Zassenhaus' mod-$p$ version of the LCS. We treat them 
as part of a general construction of the $\mathcal P$-LCS, for a 
property $\mathcal P$ of filtrations. We describe these $N$-series 
and the associated Lie algebras in the case when $G$ splits as a 
semi-direct product, in terms of the relevant data for the factors 
and the monodromy action. This allows us to generalize the well-known 
theorem of Falk--Randell regarding the LCS of split extensions to other 
versions of the LCS. In particular, we generalize the mod-$q$ version 
of Bellingeri--Gervais to any integer $q$, and we prove analogous results 
for the rational LCS and Zassenhaus' mod-$p$ LCS.

We then use the same tools to study residual properties of semi-direct 
products, and how they interact with residual properties of the factors. 
We also give a new proof of a classical theorem of Gruenberg. Finally, 
we apply our results to surface braid groups, which naturally split 
as semi-direct products, allowing us to recover and generalize  
known results about the residual nilpotency of groups of pure 
braids on surfaces.
\end{abstract}

\maketitle

\section*{Introduction}

A classical theorem of Falk and Randell~\cite{FR88} describes the lower 
central series of a split extension of groups in terms of those of the 
factors, under a homological condition on the monodromy action. This paper 
develops a general framework that extends this theorem to several 
variants of the lower central series, with applications to residual 
properties of groups and to the topology of fiber spaces. 
Specifically, we analyze the behavior of the LCS and its 
variants---the rational (torsion-free), $q$-torsion, and $p$-restricted 
LCS---under split extensions, focusing on their decompositions as 
semi-direct products. By introducing a general framework based on 
$\mathcal P$-envelopes, we extend classical results on the LCS to 
these variants, revealing new insights into their structure under 
split group extensions.

\subsection*{$N$-series} 

The study of the lower central series and the associated graded Lie algebra, 
pioneered by Hall \cite{Hall34} and Magnus \cite{Magnus35, Magnus40}, were 
generalized by Lazard \cite{Lazard} through the concept of an $N$-series, 
which serves as a powerful abstraction of the lower central series. 
An $N$-series for a group $G$ is a descending filtration by 
subgroups $G_* = \{G_n\}_{n\ge 1}$ which starts at $G_1=G$ and  
satisfies $[G_m, G_n] \subseteq G_{m+n}$, for all $m,n\ge 1$. 
Clearly, this is a central series (and thus a normal series), 
and so the quotient groups $G_n/G_{n+1}$ are abelian. 
The direct sum of these quotients, $\gr(G_*)  
= \bigoplus_{n\ge 1} G_n/G_{n+1}$, acquires the 
structure of a graded Lie algebra (over~$\Z$), 
whose Lie bracket is induced from the group commutator. 

\subsection*{Lower central series as minimal filtrations} 
The prototypical example of an $N$-series is the lower central 
series $\gamma_*(G)=\{\gamma_n(G)\}_{n\ge 1}$, that starts at 
$\gamma_1 (G)=G$ and is defined recursively by  
$\gamma_{n+1}(G) =[G,\gamma_n (G)]$. By construction, 
this is the minimal central filtration on the group $G$, 
with respect to inclusion. Variants of the lower central series 
can all be defined by their minimality with respect to some properties:
\begin{itemize}[itemsep = 2pt]
\item 
The \emph{rational LCS} $\gamma_*^\Q(G)$ is minimal among central filtrations 
$G_*$ with $G_1 = G$ and torsion-free quotients $G_i/G_{i+1}$. These are 
called \emph{central torsion-free filtrations}. 
\item 
The \emph{$q$-torsion LCS} $\gamma_*^{\qt}(G)$ (for an integer $q \geq 2$) 
is minimal among central filtrations $G_*$ such that $G = G_1$ and the 
quotients $G_i/G_{i+1}$ have all elements of order dividing $q$, that is, 
$g^q \in G_{i+1}$ if $g \in G_i$. These are called 
\emph{central $q$-torsion filtrations}. 
\item 
The \emph{$p$-restricted LCS} $\gamma_*^{\pr}(G)$ (for some prime $p \geq 2$) 
is minimal among $N_p$-series $G_*$ on $G = G_1$, where an $N$-series $G_*$ 
is an $N_p$-series if $g^p \in G_{pi}$ whenever $g \in G_i$.
\end{itemize}
These three examples appear as natural counterparts of the lower central 
series throughout the literature. On the one hand, $\gamma_*^{\pr}(G)$ 
(resp.~$\gamma_*^\Q(G)$) naturally appears in relation with 
filtrations on associative algebras over $\mathbb F_p$ (resp.~over $\Q$).  
Precisely, this series is known to coincide with the \emph{dimension series} $G \cap (1 + I^*)$, where 
$I$ is the augmentation ideal of the group ring $\mathbb F_pG$ (resp.~$\Q G$). This fact 
was already known to Zassenhaus~\cite{Zassenhaus}, and it explains the prominent roles 
played by these series in the works of Lazard~\cite{Lazard}, Quillen~\cite{Quillen}, 
Passi~\cite{Passi} and many others. 
On the other hand, Stallings \cite{Stallings} linked the lower central 
series to low-dimensional group homology, with $\gamma_*^\Q(G)$ and 
$\gamma_*^{\qt}(G)$ arising naturally for rational and $\Z/q$ coefficients, 
respectively (where $q$ need in fact not be a prime number). These 
series also appear as prominent characters in subsequent 
investigations by Cochran and Harvey \cite{CH-gt08, CH-jlms08}. 
Additionally, for a finitely generated group $G$ and a prime $p$, 
both $\gamma_*^{\pr}(G)$ and $\gamma_*^{\qt[p]}(G)$ detect whether $G$ 
is residually a finite $p$-group, with the latter series 
being more tractable in practice, as seen 
in the works of Paris~\cite{Paris} and Bellingeri and Gervais~\cite{Bellingeri-Gervais}.

\subsection*{Decomposition of minimal filtrations} 
Let $H \rtimes K$ denote a semi-direct product of groups. It is known from~\cite{Darne2} 
(see also Guaschi and Pereiro~\cite{Guaschi-Pereiro}) that the lower central series of a semi-direct product 
decomposes as a semi-direct product: 
\[
\gamma_*(H \rtimes K) = \gamma_*^K(H) \rtimes \gamma_*(K),
\]
where $\gamma_*^K(H)$ is the filtration on $H$ defined inductively by 
$\gamma_1^K(H) = H$ and $\gamma_{n+1}^K(H) = [\gamma_n^K(H), H \rtimes K]$. 
Moreover, we know that $\gamma_*^K(H) = \gamma_*(H)$ if and only if the semi-direct 
product is an \emph{almost-direct} one, in which case this decomposition 
recovers the classical result of Falk and Randell from \cite{FR}. Namely, 
$\gamma_*(H \rtimes K) = \gamma_*(H) \rtimes \gamma_*(K)$ if and only if 
$K$ acts trivially on $H^{\ab}$, and then we get a corresponding decomposition 
of the associated graded.

We want to understand how these results (recalled in~\cref{subsec:lcs-semi}) 
can be adapted to variants of the lower central series. In order to do this, 
it is helpful to introduce $\mathcal P$-envelopes and the $\mathcal P$-lower 
central series, for some properties $\mathcal P$ of filtrations. Namely, given 
a filtration $G_*$, we denote by $G_*^{\mathcal P}$ the minimal filtration 
containing $G_*$ that satisfies the property $\mathcal P$, called the 
\emph{$\mathcal P$-envelope} of $G_*$ (which may not exist in general). If $\mathcal P$ 
implies centrality, then $\gamma_*(G)^\mathcal P$, the $\mathcal P$-envelope 
of the lower central series, is the minimal filtration satisfying $\mathcal P$ 
on a group $G$, also denoted by $\gamma_*^{\mathcal P}(G)$ and called the 
\emph{$\mathcal P$-lower central series}. The variants of the lower central 
series introduced above are all of this form, for some property $\mathcal P$.

With this point of view, we see that the question of whether and how 
$\gamma_*^{\mathcal P}(H \rtimes K)$ decomposes as a semi-direct product is 
closely related to the more general question whether $(H_* \rtimes K_*)^\mathcal P = 
H_*^\mathcal P \rtimes K_*^\mathcal P$, for general filtrations $H_*$ and $K_*$ on 
$H$ and $K$. This leads us to introduce (in~\cref{subsec:P-filtrations-definitions}) 
a series of five hypotheses \ref{H0}--\ref{H4} on $\mathcal P$ under which the latter 
holds, at least for filtrations satisfying some weaker property $\mathcal Q$ (which 
we call our~\emph{context}). Moreover, in this general framework, we can generalize 
the proof of the above decomposition statements for the LCS to the $\mathcal P$-LCS, 
which allows us to understand when we have $\gamma_*^{\mathcal P}(H \rtimes K) =
\gamma_*^{\mathcal P}(H) \rtimes \gamma_*^{\mathcal P}(K)$. 
We observe that the latter equality can hold even if 
$\gamma_*(H \rtimes K) \neq \gamma_*(H) \rtimes \gamma_*(K)$, 
that is, even if $H \rtimes K$ is not an almost-direct product. 
In fact, such a decomposition holds whenever $H \rtimes K$ is merely 
a \emph{$\mathcal P$-almost direct product}, that is, if and only if 
$K$ acts trivially on $\gr_*^\mathcal P(H) = \gr(\gamma_*^\mathcal P(H))$.

\subsection*{Kaloujnine $N$-series} 
In the general framework described above, a prominent role is played by 
Kaloujnine's construction of an $N$-series on groups of automorphisms. 
Given a group $K$ acting by automorphisms on a group $H$ endowed with a filtration $H_*$ 
by normal subgroups, define $\mathcal A_j(K, H_*)$ as the subgroup of $K$ that 
preserves each $H_i$ and acts trivially on the quotients $H_i/H_{i+j}$. 
Then the filtration $\mathcal A_*(K, H_*) = \{\mathcal A_j(K, H_*)\}_{j \ge 1}$ forms 
an $N$-series on the automorphism group $\mathcal A_1(K, H_*)$, 
providing a versatile tool for analyzing group extensions.

\begin{example*}
The most famous example of such a filtration is the \emph{Johnson filtration}\/ 
on the mapping class group $K=\MCG(S)$ of a surface $S$, which corresponds to 
$H = \pi_1(S)$ and $H_* = \gamma_*(H)$ \cite{Johnson}. Similarly, 
the \emph{Andreadakis filtration} arises when $K = \Aut(F_n)$ acts 
on the free group $H = F_n$ with $H_* = \gamma_*(F_n)$ \cite{Andreadakis}.
\end{example*}

Let a group $K$ act on a group $H$ endowed with a filtration $H_*$ by normal 
subgroups.
We say that a filtration $K_*$ on a subgroup $K_1 \subseteq K$ \emph{acts on}\/ 
$H_*$ if $K_j \subseteq \mathcal A_j(K, H_*)$ for all $j \ge 1$. This means 
that each $K_j$ preserves the subgroups $H_i$ and acts trivially on $H_i/H_{i+j}$.
This is a generalization of a terminology already 
introduced in \cite{Darne1, Darne2} for $N$-series, for which the link with semi-direct 
products is most apparent: if $K_*$ and $H_*$ are $N$-series and $K_*$ acts on $H_*$, 
then the groups $H_i \rtimes K_i$ are the terms of an $N$-series $H_* \rtimes K_*$. 
This explains the relevance of Kaloujnine's construction to decompositions of lower 
central series. In particular, the filtration $\gamma_*^K(H)$ introduced above is 
minimal among the filtrations $H_*$ on $H = H_1$ which are acted upon by $\gamma_*(K)$, 
that is, such that $\gamma_*(K) \subseteq \mathcal A_*(K, H_*)$.

The concept of a \emph{$\mathcal P$-almost direct product} emerges naturally 
in this context. We say that a semi-direct product $H \rtimes K$ is a 
$\mathcal P$-almost direct product if $K$ acts trivially on $\gr_*^\mathcal P(H)$, 
that is, if $K = \mathcal A_1(K, \gamma_*^\mathcal P(H))$. 
If $\mathcal A_*(K, \gamma_*^\mathcal P(H))$ is a $\mathcal P$-filtration, 
this is equivalent to 
$\gamma_*^\mathcal P (K)$ acting on $\gamma_*^\mathcal P(H)$. 
This explains why we include in our hypotheses the assumption that 
$\mathcal A_*(K, \gamma_*^\mathcal P(H))$ is a $\mathcal P$-filtration. 
This assumption is part of a broader set of conditions (hypotheses \ref{H0}–\ref{H4}, 
detailed in \cref{subsec:P-filtrations-definitions}) ensuring that
$\mathcal A_*(K, H_*)$ is a $\mathcal P$-filtration whenever $H_*$ is, 
facilitating our analysis of semi-direct product decompositions.

\subsection*{The Johnson homomorphism} 
Johnson morphisms are powerful tools for studying the graded Lie rings 
$\gr(\mathcal A_*(K, H_*))$ associated with group actions on filtrations.
For a central filtration $H_*$, the {\em Johnson morphism}\/ is the Lie ring morphism 
\[\tau \colon \gr(\mathcal A_*(K, H_*)) \hookrightarrow  \End(\gr(H_*))\] defined by 
$\tau(\bar k)(\bar h) = \overline{(k \cdot h)k^{-1}}$, where $\End(\gr(H_*))$ 
is the graded ring of homogeneous endomorphisms of positive degree, with 
the Lie bracket induced by its associative structure.
When $H_*$ is an $N$-series, 
$\gr(H_*)$ is a Lie ring and $\tau$ takes values in $\Der(\gr(H_*))$. If a filtration 
$K_*$ acts on $H_*$, then the inclusion $K_* \subseteq \mathcal A_*(K, H_*)$ induces 
a map between the associated graded which can be composed with $\tau$ to give a map 
from $\gr(K_*)$ to $\End(\gr(H_*))$. The latter is a Lie morphism if $K_*$ is an $N$-series. 
Thus, our version of the Johnson morphism generalizes not only the classical case studied 
by Johnson \cite{Johnson} for the Torelli subgroup of the mapping class group, but also its 
standard generalization to the case where both $K_*$ and $H_*$ are $N$-series 
\cite{Habiro-Massuyeau, Papadima-Suciu-2012, Darne2}.
In the present study, the Johnson morphism is a key tool for verifying hypothesis \ref{H1}, 
which ensures that $\mathcal A_*(K, H_*)$ inherits properties from $H_*$. 
For example, if $\gr(H_*)$ is torsion-free, then $\End(\gr(H_*))$ is torsion-free, 
and thus so is $\gr(\mathcal A_*(K, H_*))$.

\subsection*{Main results} 
Our first main results are formal ones: we draw consequences from hypotheses 
\ref{H0}--\ref{H4}. Namely, we show that under \ref{H0}--\ref{H4}, taking 
$\mathcal P$-envelopes preserves central filtrations, $N$-series, and 
semi-direct products of filtrations (Lemmas~\ref{min_filtrations_are_N-series} 
and~\ref{psd_of_min_filtrations}), and we prove:
\begin{reptheorem}{dec_of_P-LCS}
Assume that $\mathcal P$-filtrations are central and that \ref{H0}--\ref{H4} hold. Then:
\begin{enumerate}[itemsep = 2pt]
    \item
    For every group $G$, the filtration $\gamma_*^{\mathcal P}(G)$ is an $N$-series.
\end{enumerate}
Moreover, if a group $K$ acts on a group $H$ by automorphisms, then:
\begin{enumerate}[itemsep = 2pt]
\setcounter{enumi}{1}
    \item
    $\gamma_*^{\mathcal P}(H \rtimes K) = 
    \gamma_*^K(H)^{\mathcal P} \rtimes \gamma_*^{\mathcal P}(K)$;
    \item
    If $H \rtimes K$ is a $\mathcal P$-almost direct product, then 
    $\gamma_*^K(H)^{\mathcal P} = \gamma_*^{\mathcal P}(H)$.
\end{enumerate}  
In particular, if $H \rtimes K$ is a $\mathcal P$-almost direct product, then 
we have a canonical isomorphism of graded Lie rings: 
$\gr_*^{\mathcal P}(H \rtimes K) \cong \gr_*^{\mathcal P}(H) 
\rtimes \gr_*^{\mathcal P}(K)$.
\end{reptheorem}
Then, we show successively that \ref{H0}--\ref{H4} hold:
\begin{itemize}[itemsep = 1.5pt]
\item 
For central torsion-free filtrations in the context of central filtrations, whence a 
decomposition theorem for the \emph{rational LCS} of semi-direct products. 
\item 
For central $q$-torsion filtrations in the context of central filtrations, 
whence a decomposition theorem for the \emph{$q$-torsion LCS} of semi-direct products. 
\item 
For $N_p$-series in the context of $N$-series, whence a decomposition 
theorem for the \emph{$p$-restricted LCS} of semi-direct products. 
\end{itemize}
Moreover, in each case, $\mathcal P$-almost direct product are well understood. 
Precisely, if $\mathcal P$ is ``being torsion-free central" (resp.~``being $q$-torsion 
central", resp.~``being an $N_p$-series"), then $H \rtimes K$ is a $\mathcal P$-almost 
direct product if and only if $K$ acts trivially on $H^{\ab} \otimes \Q$ 
(resp.~on $H^{\ab} \otimes \Z/q$, resp.~on $H^{\ab} \otimes \mathbb F_p$). 

In the case of $q$-torsion central filtrations, our decomposition theorem 
generalizes to a not necessarily prime integer~$q \geq 2$ Bellingeri and 
Gervais'~\cite[Thm.~3.3]{Bellingeri-Gervais}. In particular, we recover 
their notion of a $q$-almost direct product. A slight variation of this 
also recovers Bardakov, Bryukhanov, and 
Neshadim's~\cite[Thm.~4.1]{Bardakov-Bryukhanov-Neshchadim} 
(see Theorem~\ref{BBN_thm}).

\subsection*{Residual properties of groups} 
One of the many reasons for studying the lower central series and its rational and modular 
variants comes from the fact that these series control the corresponding residual properties 
of a group $G$. For instance, $G$ is residually nilpotent if and only if the intersection 
$\gamma_\infty(G)$ of the terms of $\gamma_*(G)$ is trivial. In order to deal with the various 
versions of the lower central series, we adapt our general framework to this setting. Namely, if 
$\mathcal P$ is a property of filtrations, we introduce the class $\mathcal R_\mathcal P$ of 
groups $G$ satisfying $\gamma_n^\mathcal P(G) = 1$ for some $n \geq 1$. Then we introduce 
hypotheses \ref{HR0}--\ref{HR2} on a property $\mathcal P$ of filtrations (two of which are 
implied by \ref{H0}--\ref{H4}) under which the intersection $\gamma_\infty^\mathcal P(G)$ of the 
terms of the series $\gamma_*^\mathcal P(G)$ is trivial if and only if $G$ is residually 
$\mathcal R_\mathcal P$. We then apply this to recover the following results:
\begin{itemize}[itemsep = 1.5pt]
\item 
A group $G$ is residually torsion-free nilpotent (RTFN) if and only if $\gamma_\infty^\Q(G) = 1$.
\item 
If $G$ is finitely generated, it is residually $p$ (that is, residually 
a finite $p$-group) if and only if $\gamma_\infty^{\qt[p]}(G) = 1$ 
(this recovers~\cite[Prop.~2.3]{Paris}).
\end{itemize}
The latter statement is part of a more general investigation of the residual 
property associated to $\gamma_*^{\qt}$ for any integer $q \geq 2$, and for 
not necessarily finitely generated groups. We call the groups satisfying 
$\gamma_\infty^{\qt} = 1$ \emph{residually nilpotent $q$-groups} (RTN$_q$ groups).

\begin{remark*}
\label{rem:cofinal}
For a prime $p$, the filtrations $\gamma^{\qt[p]}$ and $\gamma^{\pr}$ 
are \emph{cofinal}, so the corresponding residual properties are the same. 
Since working with the $p$-torsion LCS $\gamma^{\qt[p]}$ is easier, 
we do not deal with $\gamma^{\pr}$ when treating residual properties.
\end{remark*}

If we combine these new hypotheses with the previous ones, we can deduce from 
the decomposition theorem~\ref{dec_of_P-LCS} the following generalization of 
a result of Falk and Randell~\cite{FR88} to variants of the lower central series:
\begin{reptheorem}
{thm:res-P}
Under hypotheses \ref{H0}--\ref{H4} and \ref{HR0}--\ref{HR2}, 
if $H \rtimes K$ is a $\mathcal P$-almost direct product, then $H \rtimes K$ is 
residually $\mathcal R_{\mathcal P}$ if and only if $H$ and $K$ are.
\end{reptheorem}
The hypotheses are satisfied for torsion-free central filtrations and for  
$q$-torsion central filtrations, so we get:
\begin{itemize}[itemsep = 1.5pt]
\item 
If $K$ acts trivially on $H^{\ab} \otimes \Q$, then $H \rtimes K$ is RTFN 
if and only if $H$ and $K$ are.
\item 
If $H$ and $K$ are finitely generated and $K$ acts trivially on $H \otimes \F_p$, 
then $H \rtimes K$ is residually $p$ if and only if $H$ and $K$ are.
\end{itemize}
The latter result recovers a theorem of Bellingeri and 
Gervais~\cite[Cor.~3.5]{Bellingeri-Gervais}, 
which we generalize to (not necessarily finitely generated) 
RTN$_q$ groups, for any $q \geq 2$.

Our study of residual properties of groups features two 
new results. First, we provide a novel proof of Gruenberg’s theorem, 
which states that finitely generated RTFN groups are residually $p$ 
for every prime $p$ \cite[Thm.~2.1]{Gruenberg}. Our proof, while more technical, 
is constructive, offering explicit insights into the filtrations involved. 
Second, we establish Theorem~\ref{Residual_prop_of_sdp}, which shows that 
if $\mathcal P$ is a property concerning torsion in the associated graded 
of central filtrations, then a semi-direct product $H \rtimes K$ is residually 
$\mathcal R_\mathcal P$ under weaker conditions than that of Theorem~\ref{thm:res-P}. This generalizes a result of Hall~\cite{Hall58}.

\subsection*{Pure braids on surfaces}

We apply our methods (in particular, Theorem~\ref{Residual_prop_of_sdp}) to 
recover and generalize existing results about residual properties of pure 
braid groups on surfaces, by Bellingeri, Gervais, Guaschi, and 
Bardakov~\cite{Bellingeri-Gervais-Guaschi, Bardakov-Bellingeri, Bellingeri-Gervais}. 
We get a general result for any surfaces (Theorem~\ref{res_prop_of_Pn(S)}), 
whose consequence for surfaces of finite type (that is, surfaces homotopy 
equivalent to compact ones) is as follows:
\begin{repcorollary}
{res_prop_of_Pn(S)_fg}
Let $S$ be a connected surface of finite type different from the sphere 
and the projective plane, and let $n \geq 2$. Then the group $P_n(S)$ 
of pure braids on $S$ is:
\begin{itemize}[itemsep = 2pt]
    \item RTFN (hence residually~$p$ for all prime~$p$) if $S$ is orientable.
    \item residually~$2$ but not residually~$p$ for any prime $p \neq 2$ 
    if $S$ is not orientable.
\end{itemize}
\end{repcorollary}
This result was already known, except for the negative part, which 
answers a question asked byof Bellingeri and Gervais~\cite{Bellingeri-Gervais}. 
However, our proof is new, and it does not need any explicit computations 
with presentations. 

\subsection*{Further directions} 
The framework developed in this paper---particularly the results on 
semi-direct product decompositions of the LCS and its variants---has 
direct applications beyond pure group theory, and has been used in an 
essential way in two companion papers. In \cite{Suciu-pisa2024}, these 
techniques are applied to study the Alexander invariant $B(G) = G'/G''$ 
and its rational and $p$-torsion variants, together with the cohomology 
jump loci of finitely generated groups; the semi-direct product 
decomposition results established here allow one to relate these 
invariants to the LCS and its variants for a broad class of groups 
arising in topology. In \cite{Suciu-roum2024}, the same techniques 
are applied to Milnor fibrations of hyperplane arrangements with 
trivial algebraic monodromy---that is, arrangements for which the 
monodromy action on $H_1(F;\Q)$ (or even on $H_1(F;\Z)$) is 
trivial. Under these hypotheses, the decomposition theorems of this 
paper yield concrete consequences for the LCS, the Chen ranks, 
and the residual properties of the fundamental group $\pi_1(F)$ 
of the Milnor fiber $F$ (see also Example~\ref{ex:milnor}).

\subsection*{Outline of the paper} 
The first section (\cref{section_reminders}) 
introduces the main tools used in the paper: semi-direct products of groups 
and Lie rings, commutator calculus, almost-direct products and their variants, 
filtrations and associated graded, Kaloujnine's construction, Johnson morphisms, 
and decomposition of the LCS of semi-direct products. In \cref{sect:P-filtrations}, 
we turn to our general framework: we introduce our hypotheses \ref{H0}--\ref{H4}, 
we examine their validity on the basic example of central filtrations and $N$-series, 
and we draw consequences from them, showing that they imply that 
$\mathcal P$-envelopes behave nicely (Lemmas~\ref{min_filtrations_are_N-series} 
and~\ref{psd_of_min_filtrations}) and obtaining our decomposition 
theorem for the $\mathcal P$-LCS (Theorem~\ref{dec_of_P-LCS}). 
The next four sections are devoted to proving our hypotheses for 
torsion-free central filtrations (\cref{sect:lcs-rational}), for 
$q$-torsion central filtrations (\cref{sec_q-torsion}) 
and for $N_p$-series (\cref{sec_p-restricted}). In all three cases, we need a precise 
description of $\mathcal P$-envelopes. In the torsion-free case, this description is given 
by a result of Massuyeau~\cite[Lem.~4.4]{Massuyeau}, generalized to central filtrations 
(not only $N$-series). We also show (in~\cref{sec:pi-free}) that the results obtained 
for torsion-free filtrations generalize easily to $\pi$-free filtrations, where $\pi$ 
is a set of primes. In the $q$-torsion and the $p$-restricted cases, the description 
of $\mathcal P$-envelopes uses Dark's theorem to generalize results of Passi~\cite{Passi}. 
This is the subject of~\cref{sec_Dark}. We turn to residual properties of groups 
in~\cref{sect:res-semi}, which contains our results about residual properties 
of semi-direct products and a new proof of Gruenberg's theorem. Finally, as an 
application, we study the case of pure braids on surfaces in~\cref{sec:braids}, 
where we recover and generalize the results of Bardakov, Bellingeri, Gervais, 
and Guaschi about residual properties of pure braid groups on surfaces.

\tableofcontents

\section{Basic notions}
\label{section_reminders}

We start with some basic definitions and results that will be needed 
throughout the remainder of this paper.

\subsection{Semi-direct products}
\label{sec_psd} 
Given an action of a group $K$ on a group $H$ by automorphisms, that 
is, a group morphism $a\colon K \to \Aut(H)$, one can define the 
\emph{semi-direct product}\/ $H \rtimes_a K$ (often denoted 
$H \rtimes K$ for short) by the usual law 
on the set $H \times K$. Precisely, if we denote by $k \cdot h$ 
the element $a(k)(h)$ (for $k \in K$ and $h \in H$), this law 
is given by the formula:
\[
(h,k) (h',k') = (h(k \cdot h'), kk').
\]
If we identify $H$ with the normal subgroup $H \times 1$ and $K$ with 
the subgroup $1 \times K$, we see that the action of $K$ on $H$ is 
realized by conjugation in $H \rtimes K$; namely, $k h k^{-1} = k \cdot h$. 
Moreover, $H \rtimes K$ is universal for this property: given a group $G$ 
and morphisms $\alpha \colon K \to G$ and $\beta \colon H \to G$ 
such that $\beta(k \cdot h) = \alpha(k) \beta(h) \alpha(k)^{-1}$ 
(for all $k \in K$ and $h \in H$), there exists a unique morphism 
$\gamma \colon H \rtimes K \to G$ satisfying 
$\gamma|_K = \alpha$ and $\gamma|_H = \beta$. The morphism $\gamma$ is an 
isomorphism if and only if $\alpha$ and $\beta$ are injective (meaning 
that $H$ and $K$ are subgroups of $G$), $H \cap K = \{1\}$ and $HK = G$. 
In this case, we say that $G$ \emph{decomposes as a semi-direct product 
of $H$ and $K$}. Such a decomposition is equivalent to the data of a 
split short exact sequence of groups (up to isomorphism):
\[
\begin{tikzcd}
H \ar[r, hook]
&G \ar[r, two heads]
&K. \ar[l, start anchor = north west, end anchor = 
north east, bend right, shift left = 1]
\end{tikzcd}
\]
The action of $K$ on $H$ induces an action on the abelianization 
$H^{\ab}=H/H'$. If the action of $K$ on $H^{\ab}$ is trivial, 
we also say that the semi-direct product $H \rtimes K$ is 
an \emph{almost direct}\/ one.

The same story holds for Lie algebras, where actions by automorphisms 
are replaced by actions by derivations. Precisely, the data of a Lie 
morphism $a \colon \mathfrak k \to \Der(\mathfrak h)$ (denoted by 
$k \mapsto (h \mapsto k \cdot h)$, and called an~\emph{action}\/ 
of $\mathfrak k$ on $\mathfrak h$ by derivations) is equivalent 
to the data of a split short exact sequence of Lie algebras 
(up to isomorphism),
\[
\begin{tikzcd}
\mathfrak h \ar[r, hook]
&\mathfrak g \ar[r, two heads]
&\mathfrak k. \ar[l, start anchor = north west, end anchor = 
north east, bend right, shift left = 1]
\end{tikzcd}
\]
The semi-direct product $\mathfrak g = \mathfrak h \rtimes \mathfrak k$ 
is obtained from the action by endowing $\mathfrak h \oplus \mathfrak k$ 
with the Lie bracket given by
\[ 
[(h,k), (h',k')] = ([h,h'] + k \cdot h' - k' \cdot h, [k,k']).
\]
Conversely, one can recover the action $a$ from the split extension by 
sending $k \in \mathfrak k$ to $[k, -]|_{\mathfrak h} \in \Der(\mathfrak h)$, 
where the bracket is computed in $\mathfrak g$.

Several other versions of this story will appear in the present paper 
(though less explicitly), notably $p$-restricted Lie algebras (whose 
actions are via \emph{$p$-restricted derivations}, 
see~\cref{par_p-restricted_graded}), or $N$-series, whose semi-direct 
products are constructed in Proposition-Definition~\ref{psd_of_N-series}.

\subsection{Commutators}
\label{par_commutators}
We collect here some notation and basic facts related to commutator calculus.

\begin{notation}
Let $G$ be a group. If $x, y \in G$, we denote by $[x,y]$ their commutator 
$xyx^{-1}y^{-1}$, and we use the usual exponential notations $x^y = y^{-1}xy$ 
and ${}^y\! x = yxy^{-1}$ for conjugation in $G$. If $A, B \subseteq G$ are 
subsets of $G$ and $q \geq 1$ is an integer, we denote by $[A,B]$ the 
\emph{subgroup}\/ generated by commutators $[a,b]$ with $a \in A$ 
and $b \in B$, and by $A^q$ the \emph{subgroup}\/ generated by the 
elements $a^q$ with $a \in A$.
\end{notation}

The following identities hold for all $x, y, z \in G$:
\begin{itemize}[itemsep = 1.5pt]
\item $[x,y]^{-1} = [y,x]$;
\item $[x, yz]= [x, y] \cdot {}^y\! [x, z] = [x,y] \cdot [x,z] \cdot [[z,x],y]$;
\item $[[x,y],{}^y\! z] \cdot [[y,z],{}^z\! x] \cdot [[z,x],{}^x\! y ]= 1$ 
(the \emph{Witt--Hall identity}).
\end{itemize}
Of course, one can devise a lot of such identities; in particular, 
we have, for $x,y,z,t \in G$:
\begin{equation}
\label{[xy,zt]}
[xy, zt] = {}^x\! [y,z] \cdot [x,z] \cdot [zx, [y,t]] \cdot 
[y,t] \cdot {}^y\! [x,t]. 
\end{equation}

The Witt--Hall identity easiimmediately implies the so-called 
``Three Subgroup Lemma" of P.~Hall~\cite{Hall34}:

\begin{lemma}[Three Subgroups Lemma]
\label{lem:3subgroups}
Let $A,B,C$ be three subgroups of a group $G$. If two of the three 
following subgroups are trivial, then so is the third:
\[
[A, [B, C]],\ \ [B, [C, A]],\ \ [C, [A, B]].
\]
Equivalently, each one of them is contained in the normal closure 
of the two others.
\end{lemma}

\subsection{Almost-direct products and their variants}
\label{sec_ppd} 
Let $H \inj G \surj K$ be an extension of 
groups. Recall that it induces an outer action of $K$ on $H$, that is, 
a morphism $K \to \Out(H)$, sending $k$ to conjugation by a lift of 
$k$ in $G$ (which is well-defined up to an inner automorphism of $H$). 
This induces a genuine action of $K$ on $H^{\ab}$ by automorphisms, 
since inner automorphisms become trivial in $\Aut(H^{\ab})$. Note that 
this construction is functorial with respect to morphisms of 
extensions. Namely, a morphism of extensions,
\[
\begin{tikzcd}
H_1  \ar[r, hook] \ar[d, "f"]
&G_1 \ar[r, two heads] \ar[d, "f"]
&K_1  \ar[d, "f"] \\
H_2  \ar[r, hook] 
&G_2 \ar[r, two heads] 
&K_2,
\end{tikzcd}
\]
induces a morphism $f^{\ab} \colon H_1^{\ab} \to H_2^{\ab}$ which is 
equivariant with respect to $f \colon K_1 \to K_2$. 

An extension $H \inj G \surj K$ splits precisely when there is a coherent 
choice of lifts of elements of $K$, so that the associated outer action comes 
from a genuine action $K \rightarrow \Aut(H)$. It is however worth remembering 
that whereas this action depends on the choice of splitting, the corresponding 
outer action does not, and in particular \emph{the action of $K$ on $H^{\ab}$ 
does not depend on a choice of splitting}. In particular, this is true of the 
induced action of $K$ on $H_1(H, R) = H^{\ab} \otimes_\Z R$, for any commutative 
unitary ring $R$.

\begin{proposition-definition}
\label{def_ppd}
Let $R$ be a commutative unitary ring, and let us write $H_1(-)$ 
for $H_1(-,R)$. For a split extension $H \inj G \surj K$, the 
following conditions are equivalent:
\begin{enumerate}[itemsep = 2pt, label=(A\arabic*), wide=1em, leftmargin=*]
\item \label{adp1}
The induced action of $K$ on $H_1(H)$ is trivial.
\item \label{adp2}
The induced morphism $H_1(H) \to H_1(G)$ is injective.
\item \label{adp3}
Some (hence every) choice of splitting induces an isomorphism 
\[H_1(G) \cong H_1(H) \oplus H_1(K).\]
\end{enumerate}
When these equivalent conditions are satisfied, we say that $G$ in an 
\emph{$R$-almost direct product} of $H$ by $K$. If $R = \Z/q$ for some 
integer $q$, following~\cite{Bellingeri-Gervais} we speak of a 
\emph{$q$-almost direct product}, and if $R = \Z$, following~\cite{FR} 
we speak simply of an \emph{almost direct product}.
\end{proposition-definition}

We need the following lemma, where for a $R$-module $M$ acted upon by $K$, 
we denote by $M_K=H_0(K,M)$ the module of coinvariants under the action 
of $K$, i.e., the largest quotient of $M$ on which $K$ acts trivially.

\begin{lemma}
\label{H_1(psd)}
Let a group $K$ act on a group $H$ by automorphisms. Then the induced maps 
from $H_1(H)$ and $H_1(K)$ to $H_1(H \rtimes K)$ induce an isomorphism: 
\[H_1(H \rtimes K) \cong H_1(H)_K \oplus H_1(K).\]
\end{lemma}

\begin{proof}
If $R = \Z$, so that $H_1(-) = (-)^{\ab}$, this is a consequence of the usual formula 
$[x, yz] = [x, y] (x [x,z] x^{-1})$, which readily implies that the commutator 
subgroup $[H \rtimes K, H \rtimes K]$ is normally generated by $[H,H]$, $[H,K]$, 
and $[K,K]$. In order to compute the quotient $(H \rtimes K)^{\ab}$ of 
$H \rtimes K$, we can take the quotient by these three sets of relations 
successively: $(H \rtimes K)/[H,H]$ is isomorphic to $H^{\ab} \rtimes K$, 
then killing $[H,K]$ gives $(H^{\ab})_K \times K$, and finally, 
$((H^{\ab})_K \times K)/[K,K] \cong (H^{\ab})_K \times K^{\ab}$. 
This shows that:
\[
(H \rtimes K)^{\ab} \cong (H^{\ab})_K \oplus K^{\ab}.
\]
For a general ring $R$, we almost get the right formula by tensoring this with 
$R$. The conclusion follows from the existence of the canonical isomorphism 
$(H^{\ab})_K\otimes R \cong (H^{\ab}\otimes R)_K$.
\end{proof}

\begin{proof}[Proof of Proposition~\ref{def_ppd}]
Notice that the triviality of the action of $K$ on $H_1(H)$ is equivalent 
to the surjection $H_1(H) \surj H_1(H)_K$ being an isomorphism. Thus 
Lemma~\ref{H_1(psd)} readily implies that \ref{adp1} is equivalent 
to \ref{adp2}. It also implies that \ref{adp1} is equivalent 
to \ref{adp3} for any choice of section (and, since \ref{adp1} does 
not depend on the section, if \ref{adp3} holds for some choice of 
section, it does for every choice). 
\end{proof}

\begin{remark}
\label{rem:nonsplit-adp}
For non-split extensions, we still have that \ref{adp2} is equivalent 
to the exactness of the short exact sequence $H_1(H) \inj H_1(G) \surj H_1(K)$. 
This uses the fact that $H_1(H) \to H_1(G) \to H_1(K) \to 0$ is always 
exact, since the functor $H_1(-) = (-)^{\ab} \otimes_\Z R$ is left 
adjoint to the forgetful functor from $R$-modules to groups, so it 
is right-exact (it preserves cokernels). 
Moreover, we still have that \ref{adp2} implies \ref{adp1}, since 
the map $H_1(H) \to H_1(G)$ factors through $H_1(H) \surj H_1(H)_K$. 
However, \ref{adp1} does not imply \ref{adp2} anymore. For instance, 
if $G$ is $2$-nilpotent and $H = \gamma_2(G)$, then \ref{adp1} holds 
but the map $H_1(H) \to H_1(G)$ is trivial.
\end{remark}

\subsection{Filtrations on groups}
\label{par_filtrations}

Let us now introduce the main characters of this paper.

\begin{definition}
\label{def:filtration-group}
A \emph{filtration} $G_*$ on a group $G$ is a sequence of nested subgroups 
$G = G_1 \supseteq G_2 \supseteq G_3  \supseteq \cdots$. A filtration is 
\emph{central}\/ if $[G, G_i] \subseteq G_{i+1}$ for all $i \geq 1$. 
A central filtration is:
\begin{itemize}[itemsep = 1.5pt]
\item  an \emph{$N$-series}\/ if 
$[G_i, G_j] \subseteq G_{i+j}$, for all $i,j \geq 1$;
\item  {\em torsion-free}\/ if $G_i/G_{i+1}$ does not have 
elements of finite order, for all $i \geq 1$;
\item {\em of $q$-torsion} (for some integer $q \geq 1$) if 
$G_i^q \subseteq G_{i+1}$, for all $i \geq 1$;
\item {\em $p$-restricted} (for some prime number $p$) if 
$G_i^p \subseteq G_{pi}$ for all $i \geq 1$.
\end{itemize}
\smallskip
A torsion-free $N$-series is also called an $N_0$-series; a $p$-restricted 
$N$-series is also called an $N_p$-series.
\end{definition}

If $G_*$ is a central filtration, then, in particular, $[G, G_i] \subseteq G_i$, 
which means that each $G_i$ is a normal subgroup of $G$. Moreover, $[G_i, G_i] 
\subseteq G_{i+1}$ means that the quotients $G_i/G_{i+1}$ are abelian groups. 

If $G_*$ and $H_*$ are two filtrations on the same group $G$, we write 
$G_* \subseteq H_*$ if $G_i \subseteq H_i$ for all $i$. Let the lower 
central series $\gamma_*(G)$ of a group $G$ be defined, as usual, 
by $\gamma_1(G) = G$ and $\gamma_{i+1}(G) = [G,\gamma_i(G)]$ when 
$i \geq 1$. Since we will use variations on its proof, we recall
a proof of the following classical result.

\begin{proposition}
\label{LCS_min}
For any group $G$, the lower central series 
$\gamma_*(G)$ is the smallest $N$-series on $G$.
\end{proposition}

\begin{proof}
Set $\gamma_i=\gamma_i(G)$. 
The minimality of $\gamma_*(G)$ is easy to prove by induction: if $G_*$ is any 
$N$-series of $G$, we have $\gamma_1 = G = G_1$, and if $\gamma_i \subseteq G_i$ 
for some $i \geq 1$, then $\gamma_{i+1} =  [G, \gamma_i] \subseteq [G, G_i] = G_{i+1}$. 

It remains to prove that $\gamma_*(G)$ is an $N$-series. We prove by induction 
on $i \geq 1$ that for all $j \geq 1$, we have $[\gamma_i, \gamma_j] \subseteq 
\gamma_{i+j}$. This holds for $i = 1$ by definition of $\gamma_*$. Suppose that 
it holds for $i \leq k$, for some $k \geq 1$. Then, by the Three Subgroups Lemma 
(Lemma~\ref{lem:3subgroups}), for all $j \geq 1$, the subgroup 
$[\gamma_{k+1}, \gamma_j] = [[G, \gamma_k], \gamma_j]$ is contained 
in the normal subgroup of $G$ generated by $[\gamma_k, [G, \gamma_j]]$ 
and $[G, [\gamma_k, \gamma_j]]$. Using the induction hypothesis and 
the definition of $\gamma_*$, we see that both of these subgroups 
are contained in $\gamma_{j+k+1}$, which is normal in $G$. As a consequence, 
$[\gamma_{k+1}, \gamma_j] \subseteq \gamma_{j+k+1}$ for all $j \geq 1$, as required.
\end{proof}

\subsection{Associated graded}
\label{subsec:gr}
If $G_*$ is a central filtration, we have seen that the quotients 
$G_i/G_{i+1}$ are abelian, so it makes sense to consider the graded 
abelian group $\gr(G_*)\coloneqq \bigoplus_{i\ge 1} G_i/G_{i+1}$, 
whose degree $k$ part is $\gr_k(G_*) = G_k/G_{k+1}$. We call $\gr(G_*)$ 
the \emph{associated graded}\/ of the filtration $G_*$. This construction 
defines a functor from the category of central filtrations and 
filtration-preserving group morphisms to the category of graded 
abelian groups and homogeneous morphisms. We will use the fact 
that this functor is \emph{exact}. Precisely, 

\begin{definition}
\label{def:filtration}
A \emph{short exact sequence of filtrations}, $H_* \inj G_* \surj K_*$, 
is a short exact sequence $H_1 \inj G_1 \surj K_1$ 
of groups, such that the morphisms are filtration-preserving, and such that 
they induce a short exact sequence of groups
$H_i \inj G_i \surj K_i$ for each $i \geq 1$. 
\end{definition}

\begin{lemma}[{\cite[Prop.~1.24]{Darne1}}]
\label{exactness_Lie}
The functor $\gr\colon G_* \mapsto \gr(G_*)$ sends a short exact sequence 
of central filtrations to a short exact sequence of graded abelian groups.
\end{lemma}

\begin{proof}
Given a short exact sequence $H_* \inj G_* \surj K_*$ of central 
filtration and an integer $i\ge 1$, consider the commutative diagram
\[
\begin{tikzcd}
H_{i+1} \ar[r, hook] \ar[d, hook]
&G_{i+1} \ar[r, two heads] \ar[d, hook]
&K_{i+1} \ar[d, hook] \phantom{\, .}\\
H_i \ar[r, hook] 
&G_i \ar[r, two heads]
&K_i \, .
\end{tikzcd}
\]
Using the Snake Lemma, we get an induced short exact sequence of (abelian) groups, 
$\gr_i(H_*) \inj \gr_i(G_*) \surj \gr_i(K_*)$, whence our conclusion.
\end{proof}

\begin{remark}
\label{rk_gr_not_ab}
The above lemma actually makes sense for filtrations satisfying only 
$G_{k+1} \triangleleft G_k$ (for all $k$), which is the minimal requirement 
for the quotient $\gr_k(G_*) = G_k/G_{k+1}$ to be well-defined. If one wants 
to put these groups in a single graded group, one must then 
replace the direct sum by the restricted direct product.
\end{remark}

Before going further, here is a helpful convention for dealing with elements 
of the associated graded:
\begin{convention}
\label{conv_deg_wrt_filtration}
Let $G$ be a group endowed with a filtration $G_*$. Let $g$ be an element of $G$. If there is an 
integer $d$ such that $g \in G_d \setminus G_{d+1}$, it is obviously unique; we then call $d$ 
the \emph{degree} of $g$ with respect to $G_*$. The notation $\overline g$ will denote the class 
of $g$ in some quotient $G_i/G_{i+1}$ ; if the integer $i$ is not specified, it will be assumed 
that $i = d$, which means that $\overline g$ denotes the only non-trivial class induced by $g$ 
in some $G_i/G_{i+1}$. If such a $d$ does not exist (that is, if $g \in \bigcap_{i\ge 1} G_i$), 
we say that $g$ has degree $\infty$ and we put $\overline g = 0$.
\end{convention}

We will also shorten notations when referring to lower central series:
\begin{convention}
\label{conv:gr}
When no filtration is specified on a group $G$, it is implied that $G$ 
is endowed with its lower central series $\gamma_*(G)$. In particular, 
we denote $\gr(\gamma_*(G))$ simply by $\gr(G)$.
\end{convention}

If $G_*$ is an $N$-series, then $\gr(G_*)$ has a nice structure \cite[Thm.~(2.1)]{Lazard}:
\begin{proposition-definition}
\label{def_Lie_struct}
If $G_*$ is an $N$-series, then the graded abelian group 
$\gr(G_*)\coloneqq \bigoplus_{i\ge 1} G_i/G_{i+1}$ becomes a graded Lie ring (that is, 
a graded Lie algebra over~$\Z$) when endowed with the Lie bracket $[-,-]$ induced 
by commutators in $G$. Precisely, using Convention~\ref{conv_deg_wrt_filtration}, 
this bracket is defined by:
\[
[\overline x, \overline y]\coloneqq \overline{[x,y]} \in \gr_{i+j}(G_*),\ 
\text{ for all $x \in G_i$  and  $y \in G_j$.}\ 
\]
\end{proposition-definition}

\begin{proof}[Sketch of proof]
It is an easy exercise to deduce from the identities of \cref{par_commutators} that the 
commutator map $[-,-] \colon G_i \times G_j \to G_{i+j}$ induces a well-defined map from 
$\gr_i(G_*) \times \gr_j(G_*)$ to $\gr_{i+j}(G_*)$, and that this defines a Lie bracket 
on $\gr(G_*)$. In particular, the Jacobi identity is deduced from the Witt--Hall identity.
\end{proof}

Of course, filtration-preserving group morphisms do preserve commutators, so they induce 
graded linear maps preserving the bracket on the associated graded. Hence the construction 
$\gr(-)$ described above defines a functor from the category of $N$-series and 
filtration-preserving group morphisms to the category of graded Lie rings. Moreover, it is 
still \emph{exact}\/ (Lemma~\ref{exactness_Lie}). Precisely, it turns exact sequence of 
$N$-series (exact sequences of filtrations which are $N$-series) into exact sequences of 
graded Lie rings (sequences of morphisms of graded Lie rings that are exact as sequences 
of $\Z$-linear maps).

Since products of commutators become sums of brackets inside the Lie ring, the following 
fundamental property follows easily from the definition of the lower central series:
\begin{proposition}
\label{engdeg1} 
The Lie ring $\gr(G)$ is \emph{generated in degree $1$}. Precisely, it is generated 
(as a Lie ring) by $\gr_1(G) = G^{\ab}$. As a consequence, if $G$ is finitely generated 
(as a group), then each $\gr_n(G)$ is also finitely generated (as an abelian group).
\end{proposition}

\begin{example}
\label{ex:free}
Let $F_n$ be the free group on $n$ generators, and let $\LL_n$ be the 
free Lie algebra of rank $n$. Work of P.~Hall \cite{Hall34}, 
W.~Magnus  \cite{Magnus35, Magnus40}, and E.~Witt from the 1930s 
(see \cite{MKS,Serre}) shows that the canonical map $\Z^n\isom \gr_1(F_n)$ 
induces an isomorphism of graded Lie algebras, $\LL_n\isom \gr(F_n)$. 
Consequently, the groups $\gr_i(F_n)$ are torsion-free of rank equal to 
$\tfrac{1}{i}\sum_{d\mid n} \mu(d) n^{i/d}$, 
where $\mu\colon \N\to \{0,\pm 1\}$ is the M\"{o}bius function.
\end{example}

\subsection{Kaloujnine $N$-series}
\label{par_Kaloujnine}

Kaloujnine \cite{Kaloujnine1, Kaloujnine2} introduced 
fundamental $N$-series on subgroups of automorphism 
groups of filtered groups. These series were later 
studied notably by Hall~\cite{Hall58}, by Andreadakis 
for automorphisms of free groups~\cite{Andreadakis}, 
and by Johnson for Mapping Class Groups~\cite{Johnson}. 

Let $H = H_1 \supseteq H_2 \supseteq H_3  \supseteq \cdots$ 
be a filtration of a group $H$ by \emph{normal}\/ subgroups 
(but not necessarily an $N$-series). Let a group $K$ act 
on $H$ by automorphisms; this action will be denoted by 
$(k,g) \mapsto k \cdot g$. We can consider the subgroup 
$\mathcal A_0(K, H_*)$ of $K$ consisting of those elements $k$ 
acting by a filtration-preserving automorphism
whose inverse 
is also filtration-preserving. Then we can consider its subgroup 
$\mathcal{A}_1(K, H_*)$, which consists of those elements of 
$\mathcal{A}_0(K, H_*)$ acting trivially on $\gr(H)$. More generally, 
if $j \geq 1$, define $\mathcal{A}_j(K, H_*)$ as the subgroup of 
elements of $\mathcal{A}_0(K, H_*)$ acting trivially on all the quotients 
$H_i/H_{i+j}$. It is clear that $\mathcal{A}_{j+1} \subseteq 
\mathcal A_j$ for every $j \geq 0$, and that $\mathcal A_j$ 
is a normal subgroup of $\mathcal{A}_0$, since it is the kernel:
\begin{equation}
\label{eq:AjKH}
\mathcal A_j(K, H_*) = \ker\Big(\!
\mathcal{A}_0(K, H_*) \longrightarrow 
\prod\limits_{i\ge 1} \Aut(H_i/H_{i+j}) 
\!\Big).
\end{equation}
We can reformulate this definition in a way more adapted to 
commutator calculus, provided we consider commutators in the semi-direct 
product $H \rtimes K$, where for $k \in K$ and $h \in H$, we have: 
\[
[k,h] = (k \cdot h) h^{-1}.
\]
This allows us to write our definition as:
\begin{equation}
\label{eq:AjKH-bis}
\mathcal{A}_j(K, H_*) = \big\{ k \in K \mid 
[k, H_i] \subseteq H_{i+j}\: \text{ for all $i \geq 1$}\big\}.
\end{equation}
The filtration $\mathcal A_*(\Aut(H), H_*)$, denoted simply by 
$\mathcal A_*(H_*)$, is called the \emph{Kaloujnine filtration 
associated to $H_*$}. Observe that if the morphism 
$a\colon K \to \Aut(H)$ represents the action of $K$ on $H$, 
then $\mathcal{A}_*(K,H_*) = a^{-1}(\mathcal{A}_*(H_*))$. 

When furthermore $H_* = \gamma_* (H)$, we denote this filtration 
by $\mathcal A_*(H)$, and we call it the \emph{Andreadakis 
filtration associated to $H$}. Notice that $\mathcal A_1(H)$ 
is the group $\mathrm{IA}_H$ of automorphisms of $H$ acting trivially 
on $H^{\ab}$. Indeed, elements of $\Aut(H)$ act on $\gr(H)$ 
by Lie ring automorphisms, so we deduce from 
Proposition~\ref{engdeg1} that acting trivially on all the 
quotients $\gamma_k H/\gamma_{k+1} H = \gr_k(H)$ is equivalent 
to acting trivially on $H^{\ab} = \gr_1(H)$.

\begin{proposition}
\label{A*_N-series}
Let $H_*$ be a filtration of a group $H = H_1$ by normal 
subgroups and let $K$ be a group acting by automorphisms 
of $H$. Setting $\mathcal{A}_*\coloneqq \mathcal{A}_*(K, H_*)$, 
we have
\[
\left[\mathcal{A}_j, \mathcal{A}_k \right] \subseteq 
\mathcal{A}_{j+k},\ \text{ for all $j, k \geq 0$. }
\]
In particular, $\mathcal{A}_*$ is an $N$-series on $\mathcal{A}_1$.
\end{proposition}

\begin{proof}
We work in the semi-direct product $H \rtimes \mathcal A_0 \subseteq H \rtimes K$. 
In there, the subgroups $H_i$, which are normal in $H$ and stable under the action 
of $\mathcal A_0$, are normal subgroups. Now, by the Three Subgroup Lemma 
(Lemma~\ref{lem:3subgroups}), $[[\mathcal A_j, \mathcal A_k], H_i]$ is contained 
in the normal subgroup generated by 
$[\mathcal A_j, [\mathcal A_k, H_i]] \subseteq [\mathcal A_j, H_{i+k}] 
\subseteq H_{i+j+k}$ and $[\mathcal A_k, [\mathcal A_j, H_i]] 
\subseteq [\mathcal A_k, H_{i+j}] \subseteq H_{i+j+k}$. 
As a consequence, $[[\mathcal A_j, \mathcal A_k], H_i] \subseteq H_{i+j+k}$, which 
means exactly that $[\mathcal A_j, \mathcal A_k] \subseteq \mathcal A_{j+k}$.
\end{proof}

\begin{corollary}[\textup{see also~\cite[(12)--(14)]{Hall58}}]
\label{LCS_acts_on_central}
Let $G_*$ be a central filtration on a group $G = G_1$. 
Then, for all $i,j \geq 1$:
\[
[\gamma_i(G), G_j] \subseteq G_{i+j}.
\]
\end{corollary}

\begin{proof}
Let $G$ act on itself by conjugation: $g \cdot h = ghg^{-1}$. 
Notice that commutators in $G \rtimes G$ are computed as in $G$, since 
$(g \cdot h)h^{-1} = ghg^{-1}h^{-1}$. As a consequence, the centrality 
of $G_*$ (which means that $[G, G_i] \subseteq G_{i+1}$ for all 
$i \geq 1$) is equivalent to $G = \mathcal A_1(G,G)$. 
Hence $\mathcal A_*(G,G)$ is a filtration on $G$, and by 
Proposition~\ref{A*_N-series}, it is central. By construction, 
$\gamma_*(G)$ is the minimal \emph{central} filtration 
of $G$. This implies that $\gamma_*(G) \subseteq \mathcal A_*(G,G)$, 
which is equivalent (by definition of $\mathcal A_*$) to the 
desired conclusion.
\end{proof}

\begin{remark}
We can \emph{recover}\/ the fact that $\gamma_*(G)$ is an $N$-series 
(Proposition~\ref{LCS_min}) as a corollary of this result, whose proof 
uses only the definition of $\gamma_*(G)$. In fact, this would be nothing 
more than a reformulation of the proof of Proposition~\ref{LCS_min} 
using the language of the Kaloujnine filtration.
\end{remark}

\begin{remark}
Another direct consequence of Corollary~\ref{LCS_acts_on_central} (up to a slight 
adaptation) is the classical formula $[\gamma_i(G), \mathcal Z^j(G)] \subseteq 
\mathcal Z^{j-i}(G)$ if $j \geq i$ (the case $i = j$ is \cite[Cor.~0.33]{Baumslag}). 
Here, $\mathcal Z^*$ is the usual sequence of iterated centers, defined by 
$\mathcal Z^0(G) = \{1\}$ and $\mathcal Z^{i+1}(G)$ is the subgroup of elements 
$g \in G$ that are central modulo $\mathcal Z^i(G)$, that is, such that 
$[g, G] \subseteq \mathcal Z^i(G)$. Then $\mathcal Z^{-*}$ is a central 
filtration, to which the corollary may be applied if we extend it to 
filtrations indexed by $\Z$, which does not present any difficulty.
\end{remark}

The following statement will allow us to construct $N$-series on semi-direct 
products from $N$-series on the factors:

\begin{proposition-definition}
\label{psd_of_N-series}
Let $K$ be a group acting on another group $H$ by automorphisms, and let $H_*$ and 
$K_*$ be $N$-series, respectively, on $H = H_1$ and $K = K_1$. Then the following 
assertions are equivalent:
\begin{enumerate}[itemsep = 2pt]
    \item \label{sdp1}
    The subgroups $H_i \rtimes K_i$ form an $N$-series on $H \rtimes K$;
    \item \label{sdp2}
    $[K_j, H_i] \subseteq H_{i+j}$, for all $i,j \geq 1$;
    \item \label{sdp3}
    $K_j \subseteq \mathcal A_j(K, H_*)$, for all $j \geq 1$.
\end{enumerate}
When any one of these conditions is satisfied, we say that $K_*$ \emph{acts on} $H_*$, 
and then the $N$-series $(H_i \rtimes K_i)_{i \geq 1}$, denoted by $H_* \rtimes K_*$, is 
called the \emph{semi-direct product}\/ of $H_*$ by $K_*$. 
\end{proposition-definition}

\begin{proof}
The equivalence between \eqref{sdp2} and \eqref{sdp3} is by definition of $\mathcal A_*$, 
and \eqref{sdp1} clearly implies \eqref{sdp2}. Let us now suppose that \eqref{sdp2} holds. 
Applying it for $j = 1$, we get $[K, H_i] \subset H_i$, which means that the subgroups $H_i$ 
are stable under conjugation by $K$, hence normal in $H \rtimes K$. In particular, 
$H_i \rtimes K_i$ is a well-defined subgroup of $H \rtimes K$, for each $i\ge 1$. 
Now, for $h \in H_i$, $h' \in H_j$, $k \in K_i$ and $k' \in K_j$, formula~\eqref{[xy,zt]} 
gives:
\[
[hk,h'k'] = {}^h [k, h'] \cdot [h,h'] \cdot [h'h,[k,k']] \cdot [k,k'] \cdot {}^{h'} [h,k'].
\]
Using the hypothesis and the normality of the subgroups $H_i$, one sees that 
$[hk,h'k'] \in H_{i+j} \rtimes K_{i+j}$, whence the conclusion.
\end{proof}

\begin{remark}
The construction of~Proposition-Definition~\ref{psd_of_N-series} is 
also the semi-direct product construction in the category of 
$N$-series from~\cite{Darne2}.
\end{remark}

When $K_* = \gamma_*(K)$, the last condition can be simplified:
\begin{corollary}
\label{action_of_LCS}
Let $K$ act on $H$ by automorphisms, and let $H_*$ be an $N$-series on $H = H_1$ 
and $K = K_1$. Then the following assertions are equivalent:
\begin{enumerate}[itemsep = 1.5pt]
    \item The subgroups $H_i \rtimes \gamma_i(K)$ form an $N$-series on $H \rtimes K$;
    \item $[\gamma_j(K), H_i] \subseteq H_{i+j}$, for all $i,j \geq 1$;
    \item $\mathcal A_1(K, H_*) = K$ (that is, $[K, H_i] \subseteq H_{i+1}$ for all $i \geq 1$).
\end{enumerate}
\end{corollary}

\begin{proof}
If the third condition is satisfied, then $\mathcal A_*(K, H_*)$ is a central 
filtration on the whole of $K$, so by minimality of the lower central series, 
$\gamma_*(K) \subseteq \mathcal A_*(K, H_*)$.
\end{proof}

\begin{remark}
\label{rem:simp}
The same kind of simplification can be used if $K_*$ is minimal for a certain 
property satisfied by $\mathcal A_*(K, H_*)$. 
\end{remark} 

Filtrations of the form $\mathcal A_*(K, H_*)$ can be studied using the corresponding 
Johnson homomorphism, which we introduce next.

\subsection{Associated graded and Johnson homomorphisms}
\label{par_Johnson_Lie}

Let $H_*$ and $K_*$ be $N$-series, and suppose that $K_*$ acts on $H_*$ 
(as in Proposition-Definition~\ref{psd_of_N-series}). Then, 
by construction, we have a split short exact sequence of $N$-series:
\[\begin{tikzcd}
H_* \ar[r, hook]
&H_* \rtimes K_* \ar[r, two heads]
&K_*, \ar[l, start anchor = north west, end anchor = north east, bend right, shift left = 1]
\end{tikzcd}\]
where the semi-direct product is the filtration from 
Proposition-Definition~\ref{psd_of_N-series}. Using Lemma~\ref{exactness_Lie}, 
we get an associated split exact sequence of graded Lie rings:
\[
\begin{tikzcd}
\gr(H_*) \ar[r, hook]
&\gr(H_* \rtimes K_*) \ar[r, two heads]
&\gr(K_*). \ar[l, start anchor = north west, end anchor = north east, bend right, shift left = 1] 
\end{tikzcd}
\]
Thus, $\gr(H_* \rtimes K_*)$ decomposes as a semi-direct product of graded Lie rings  (see~\cref{sec_psd}):
\[
\gr(H_* \rtimes K_*) \cong \gr(H_*) \rtimes \gr(K_*).
\]
Such a semi-direct product is determined by an \emph{action} of Lie ring by derivations, 
which is a restriction of the adjoint action of $\gr(H_* \rtimes K_*)$ on itself 
(see~\cref{sec_psd}). Precisely, an element $\bar k$ acts via $\bigl[\bar k, -\bigr]$, 
and since brackets in $\gr(H_* \rtimes K_*)$ come from commutators in $H \rtimes K$, 
the action is described explicitly by the formula:
\[
\overline k \cdot \overline h = \overline{(k \cdot h)h^{-1}},\ 
\text{ for all $k \in K$, $h \in H$} .
\]
This action can also be seen as a morphism from $\gr(K_*)$ to the graded Lie ring $\Der(\gr(H_*))$ of homogeneous derivations of $\gr(H_*)$, called the \emph{Johnson homomorphism}\/ associated to the action of  $K_*$ on $H_*$:
\[
\tau\colon \left\{\renewcommand{\arraystretch}{1.5}
\begin{array}{clc}
\gr(K_*) &\longrightarrow &\Der(\gr(H_*)) \\
\overline k &\longmapsto & \left(\overline h \mapsto \overline{(k \cdot h)h^{-1}}\right).
\end{array}
\right.
\]

\begin{proposition}
\label{Johnson_inj_Lie}
Let $H_*$ be an $N$-series on a group $H = H_1$, and let some group $K$ 
act on $H$ by automorphisms. Then the Johnson homomorphism associated to 
the action of $\mathcal A_*(K, H_*)$ on $H_*$ is injective.
\end{proposition}

\begin{proof}
Let $k \in \mathcal A_j$ such that $\overline k \in \gr_j(\mathcal A_*)$ is in 
$\ker(\tau)$. The components of the derivation $\tau(\bar k)$ are the maps 
$H_i/H_{i+1} \rightarrow H_{i+j}/H_{i+j+1}$ induced by $[k, -] \colon H \rightarrow H$ 
(computed in $H \rtimes K$). They are trivial if and only if $[k, H_i] \subseteq H_{i+j+1}$ 
for all $i \geq 1$, which means exactly that $k \in \mathcal A_{j+1}$ or, equivalently, 
that $\overline k = 0 \in \gr_j(\mathcal A_*)$. Thus, $\ker(\tau) = \{0\}$.
\end{proof}

\begin{corollary}
\label{A_*_and_torsion_Lie}
If $H_*$ is an $N_0$-series (resp.~an $N$-series of $q$-torsion), then so is 
$\mathcal A_*(K, H_*)$.
\end{corollary}

\begin{proof}
The associated graded $\gr(\mathcal A_*(K, H_*))$ injects into $\Der(\gr(H_*))$. 
If $\gr(H_*)$ is torsion-free (resp.~of $q$-torsion), then so is $\Der(\gr(H_*))$, 
and in turn, so is $\gr(\mathcal A_*(K, H_*))$.
\end{proof}

\subsection{Lower central series of semi-direct products}
\label{subsec:lcs-semi}

We now recall the description of the lower central series of semi-direct products 
of groups. This description is based on the following general construction:

\begin{proposition-definition}[{\cite[Def.~3.3]{Darne2}}]
\label{SCD_relative}
Let $G$ be a group, and let $H$ be a normal subgroup of $G$. The following formulas 
define an $N$-series $\gamma_*^G (H)$ on $H$:
\[
\gamma_1^G (H) \coloneqq H, \qquad 
\gamma_{k+1}^G(H) \coloneqq [G, \gamma_k^G (H)].
\]
\end{proposition-definition}

\begin{proof}
The proof of Proposition~\ref{LCS_min} works verbatim in this context.
\end{proof}

If a group $K$ acts on another group $H$ by automorphisms, then we denote 
$\gamma^{H \rtimes K}_*(H)$ simply by  $\gamma^K_*(H)$. In this context, 
the inductive definition becomes:
\begin{equation}
\label{eq:gammaKH}
\gamma_{k+1}^K(H) = [H \rtimes K, \gamma_k^K (H)] = 
[H, \gamma_k^K (H)] [K, \gamma_k^K (H)].
\end{equation}
Notice that by construction, $[K, \gamma^K_i(H)] \subseteq \gamma^K_{i+1}(H)$ 
for all $i \geq 1$. Moreover, an easy induction shows that $\gamma^K_*(H)$ is 
in fact the smallest $N$-series $H_*$ on $H$ satisfying $[K, H_i] \subseteq H_{i+1}$ 
for all $i \geq 1$. This condition is equivalent to $K = \mathcal A_1(K, H_*)$, 
which in turn is equivalent to $\gamma_*(K) \subseteq \mathcal A_*(K, H_*)$, 
by minimality of the lower central series. In other words:

\begin{lemma}
\label{minimality_of_LCS^K}
Let $K$ act on $H$ by group automorphisms. Then $\gamma^K_*(H)$ is the minimal 
filtration on $H$ which is acted upon by $\gamma_*(K)$ \emph{via} the action 
of $K$ on $H$.
\end{lemma}

\begin{remark}
\label{rem:gp-filt}
In particular, we find that $[\gamma_i(K), H] \subseteq \gamma^K_{i+1}(H)$. 
As a consequence, we have that 
$\gamma^K_{k+1}(H) = [H, \gamma_k^K (H)] [K, \gamma_k^K (H)][\gamma_k(K), H]$, 
which was the description of this filtration used by Guaschi and Pereiro 
in~\cite[Thm.~1.1]{Guaschi-Pereiro}.
\end{remark}

\begin{proposition}[{\cite[Prop.~3.4]{Darne2}; see also~\cite[Thm.~1.1]{Guaschi-Pereiro}}]
\label{lcs_of_sdp}
Let $K$ be a group acting on another group $H$ by automorphisms. 
Then, for all $i \geq 1$, 
\[
\gamma_i(H \rtimes K) = \gamma^K_i(H) \rtimes \gamma_i(K).
\]
\end{proposition}

\begin{proof}
Since $\gamma_*(K)$ acts on $\gamma^K_*(H)$, the left term defines an $N$-series,
which contains the minimal $N$-series $\gamma_*(H \rtimes K)$ on $H \rtimes K$. 
But $\gamma_*(K)$ and $\gamma^K_*(H)$ are contained in $\gamma_*(H \rtimes K)$ 
by construction, whence the equality.
\end{proof}

At the level of associated Lie rings, by applying the construction 
from~\cref{par_Johnson_Lie}, we get a semi-direct product decomposition:
\begin{corollary}
\label{cor:gr-semidirect}
Let $K$ be a group acting on another group $H$ by automorphisms. 
Then:
\[
\gr(H \rtimes K) \cong \gr(\gamma_*^K(H)) \rtimes \gr(K).
\]
\end{corollary}

\begin{example}
\label{ex:baumslag-solitar}
Let $n \neq 0$ be an integer, and $G = \BS(1,n)$ be the Baumslag--Solitar 
group with presentation $\langle t,a\mid tat^{-1}=a^n\rangle$ (for $n = -1$, 
this is the fundamental group of the Klein bottle). The group $G$ is metabelian, with 
$G^{\ab}=\Z$ generated by $t$ and $G'\cong \Z[1/n]$ normally generated by $a$. 
Thus, $G \cong \Z[1/n] \rtimes \Z$, where $t$ acts by multiplication by $n$ 
(which corresponds to $a \mapsto a^n$). Then if $x \in \Z[1/n]$, we have 
$[t^k,x] = (t \cdot x) - x = (n-1)x$. We easily deduce that 
$\gamma^{\Z}_k(\Z[1/n]) = (n-1)^{k-1} \Z[1/n]$  for all $k \geq 1$, 
and thus $\gamma_k(G) =(n-1)^{k-1} \Z[1/n]$ whenever $k \geq 2$. Hence, we 
conclude that $\gr(G) \cong A \rtimes \Z$, where $\Z$ is an abelian Lie 
ring concentrated in degree one, generated by $\tau = \overline t$, and 
the abelian Lie ring $A$ has one copy of $\Z/(n-1)$ in each degree 
$k \geq 0$, generated by the class $\alpha_k$ of $(n-1)^{k-1} \in \Z[1/n]$. 
The non-trivial Lie brackets (corresponding to the action of $\Z$ on $A$) 
come from $[\tau, \alpha_k] = \alpha_{k+1}$.
\end{example}

Under the right conditions, $\gamma^K_*(H)$ is in fact equal to 
$\gamma_*(H)$. This recovers the decomposition result of 
Falk and Randell~\cite{FR}, with a shorter proof, due to 
our taking into account the Lie ring structure on the 
associated graded:

\begin{proposition}[{\cite[Prop.~3.5]{Darne2}}]
\label{lcs_of_adp}
Let $K$ be a group acting on another 
group $H$ by automorphisms. Suppose that the action of $K$ on $H^{\ab}$ is trivial 
(that is, $H \rtimes K$ is an almost direct product). Then $\gamma_*^K(H) = \gamma_*(H)$, 
and we have a decomposition $\gr(H \rtimes K) \cong \gr(H) \rtimes \gr(K)$ 
into a semi-direct product of graded Lie rings.
\end{proposition}

\begin{proof}
The action of $K$ on $H$ induces an action of $K$ by automorphisms of Lie 
rings on $\gr(H)$, which is generated by $\gr_1(H)=H^{\ab}$. If $K$ 
acts trivially on $H^{\ab}$, then $K$ acts trivially on $\gr(H)$, that is, 
$K = \mathcal A_1(K, \gamma_*(H))$. This means that $\gamma_*(K)$ 
acts on $\gamma_*(H)$ (see Corollary~\ref{action_of_LCS}), so 
$\gamma_*^K(H) \subseteq \gamma_*(H)$ by Lemma~\ref{minimality_of_LCS^K}. 
The reverse inclusion holds by minimality of the lower central series, 
whence the conclusion.
\end{proof}

\begin{example}
\label{ex:pure braids}
Let $P_n$ be the Artin pure braid group on $n$ strands. Forgetting the last 
strand induces a split exact sequence, $F_{n-1} \inj P_n \surj P_{n-1}$, 
with monodromy given by the Artin embedding $P_{n-1}\inj \Aut (F_{n-1})$. 
Since pure braids act trivially on $F_{n-1}^{\ab}$, the extension is 
an almost direct product, and thus $\gr(P_n) \cong \gr(F_{n-1}) \rtimes \gr(P_{n-1})$, 
from which it follows that 
$\gr(P_n)\cong \LL_{n-1}\rtimes ( \cdots \rtimes (\LL_{2}\rtimes \LL_1) \cdots)$. 
More generally, fundamental groups of complements of supersolvable 
arrangements can be realized as iterated almost direct products of 
finitely generated free groups \cite{FR}, and thus their associated 
graded Lie algebras are iterated semi-direct products of finitely 
generated free Lie algebras. Another generalization is to pure braid 
groups on planar surfaces (see~\cref{par:planar}).
\end{example}
 
\begin{example}
\label{ex:raag-bb}
Let $G_{\Gamma}=\langle v\in V \mid \text{$[v,w] = 1$ if $\{v,w\} \in E$}\rangle$ 
be the right-angled Artin group associated to a finite simplicial graph 
$\Gamma=(V,E)$. The associated graded Lie algebra $\gr(G_{\Gamma})$ 
was determined in \cite{Duchamp-Krob, Papadima-Suciu-2006}. 
The corresponding Bestvina--Brady group, $N_{\Gamma}$,  
is the kernel of the homomorphism $G_{\Gamma}\surj \Z$ that sends each 
generator $v\in V$ to $1\in \Z$, see \cite{BB}. As shown in 
\cite{Papadima-Suciu-2007}, if $\Gamma$ is connected, then 
the split extension $G_{\Gamma}=N_{\Gamma}\rtimes \Z$ is an 
almost direct product, and so $\gr(G_{\Gamma})\cong \gr(N_{\Gamma}) \rtimes \Z$, 
with $\Z$ concentrated in degree $1$. In particular, 
$\gr_{\geq 2}(N_{\Gamma})\cong \gr_{\geq 2}(G_{\Gamma})$.
\end{example}

\begin{example}
\label{ex:artin-kernel}
More generally, one may assign an 
integer weight $m_v$ to each vertex $v$ of $\Gamma$, 
and define the corresponding {\em Artin kernel}\/ $N_{\Gamma,\chi}$ as 
the kernel of the homomorphism $\chi\colon G_{\Gamma} \to \Z$ 
that sends each generator $v$ to $m_v$. 
In \cite[Thm.~6.2]{Papadima-Suciu-2009}, a necessary and sufficient condition  
for the triviality of the action of $\Z$ on $H_1(N_{\Gamma,\chi};\Q)$ is given. 
As shown in \cite[Lem.~9.1]{Papadima-Suciu-2009}, when this condition 
is satisfied, the group $N_{\Gamma,\chi}$ is finitely generated, and, 
in fact, $\Z$ acts trivially on $N_{\chi}^{\ab}$. Thus, once again, 
$\gr(G_{\Gamma}) \cong \gr(N_{\Gamma, \chi}) \rtimes \Z$ with $\Z$ 
concentrated in degree $1$, and in particular, 
$\gr_{\geq 2}(N_{\Gamma, \chi})\cong \gr_{\geq 2}(G_{\Gamma})$.
\end{example}

\subsection{Actions of central filtrations}
\label{subsec:act-central}
Most of the filtrations that we want to deal with are $N$-series. For 
such a series, we can fully take advantage of the Lie ring structure 
on the associated graded, as we just did. However, motivated 
by the fact that $N$-series occur from any filtration by normal 
subgroups (Proposition~\ref{A*_N-series}) and inspired by 
\cite{Bass-Lubotzky}, we look at what happens with less 
structured filtrations. We are thus led to the following definition.

\begin{definition}
\label{def_action_filtrations}
An \emph{action}\/ of a filtration $K_*$ on a filtration $H_*$ is a group 
action of $K = K_1$ on $H = H_1$ satisfying the equivalent conditions:
\begin{enumerate}[itemsep=1.5pt]
    \item  $[K_j, H_i] \subseteq H_{i+j}$, for all $i, j \geq 1$.
    \item  $K_j \subseteq \mathcal A_j(K, H_*)$, for all $j \geq 1$.
\end{enumerate}
These conditions are equivalent by the definition of $\mathcal A_*$ 
(see~\cref{par_Kaloujnine}).
\end{definition}

\begin{remark}
\label{rk_action_central}
This may seem too strong a requirement. In particular, the conditions needed to get 
well-defined split short exact sequences of filtrations from the short exact sequence 
$H \inj H \rtimes K \surj K$ are weaker. Precisely, for 
$H_i \rtimes K_i$ to be well-defined subgroups, we only need $K_i$ to stabilize $H_i$ 
(that is, $[K_i, H_i] \subseteq H_i$). If we want them to be normal subgroups, assuming 
that the subgroups $K_i$ (resp.~$H_i$) are normal in $K$ (resp.~in $H$), we need 
$[K, H_i] \subseteq H_i$ and $[K_i, H] \subseteq H_i$. If we want $H_* \rtimes K_*$ 
to be central, assuming that $H_*$ and $K_*$ are, we need $[K, H_i] \subseteq H_{i+1}$ 
and $[K_i, H] \subseteq H_{i+1}$. However, our stronger requirement is a reasonable 
one if we want to get interesting structure at the level of associated graded. 
\end{remark}

\begin{example}
\label{ex_LCS_acts_on_central}
Corollary~\ref{LCS_acts_on_central} says that $\gamma_*(G)$ acts on $G_*$ through 
the conjugation action of $G$ on itself, provided $G_*$ is a central filtration on $G$.
\end{example}

We first investigate the behavior of the action of a single automorphism with 
respect to a given filtration.

\begin{lemma}
\label{general_Johnson(phi)}
Let $H_*$ be a filtration by normal subgroups, and let $\varphi \in \mathcal A_j(H_*)$. 
Then $[\varphi,-]$ (computed in $H \rtimes \Aut(H)$, that is: $[\varphi, h] = \varphi(h)h^{-1}$) 
induces well-defined maps $\tau_i(\varphi)\colon \gr_i(H_*) \to \gr_{i+j}(H_*)$ with the 
following properties.
\begin{itemize}[itemsep=1.5pt]
    \item These maps depend only on the class of $\varphi$ modulo $\mathcal A_{j+1}(H_*)$.
    \item They are all trivial if and only if $\varphi \in \mathcal A_{j+1}(H_*)$.
    \item If $H_*$ is central, then the maps $\tau_i(\varphi)$ are linear; together 
    they define a linear endomorphism $\tau(\varphi)$ of $\gr_*(H_*)$ of degree $j$. 
    \item If $H_*$ is in fact an $N$-series and $j \geq 1$, then $\tau(\varphi)$ is 
    a derivation of the Lie ring $\gr(H_*)$.
\end{itemize}
\end{lemma}

\begin{proof}
Since $\varphi \in \mathcal A_j$, the map $[\varphi, -]$ sends $H_i$ to $H_{i+j}$, 
for all $i \geq 1$. Now, if $h \in H_i$ and $u \in H_{i+1}$, then 
$[\varphi, hu] = [\varphi, h]({}^h [\varphi,u])$. We have $[\varphi,u] \in H_{i+j+1}$ 
and $H_{i+j+1} \triangleleft H$, so modulo $H_{i+j+1}$, we find that 
$[\varphi, hu] \equiv [\varphi,h]$. Hence $[\varphi, -]$ induces a well-defined map 
$\tau_i(\varphi)$ from $\gr_i(H_*)$ to $\gr_{i+j}(H_*)$, for each $i \geq 1$. Now, 
if $\nu \in \mathcal A_{j+1}$, then $[\varphi \nu, h] = {}^\varphi [\nu,h] [\varphi ,h]$, 
and $[\nu,h] \in H_{i+j+1}$, which is stable under the action of $\varphi$, so 
$[\varphi \nu, h] \equiv [\varphi ,h]$ modulo $H_{i+j+1}$. This means that the 
$\tau_i(\varphi)$ depend only on the class of $\varphi$ modulo $\mathcal A_{j+1}(H_*)$. 
Moreover, $\tau_i(\varphi) = 1$ for all $i \geq 1$ if and only if 
$[\varphi, H_i] \subseteq H_{i+j+1}$ for all $i \geq 1$, that is, 
if and only if $\varphi \in \mathcal A_{j+1}$.

Suppose now that $H_*$ is central. Then action of $H$ on the abelian group 
$\gr_{i+j}(H)$ induced by conjugation is trivial, so if $h, h' \in H_i$, then 
$[\varphi, hh'] = [\varphi, h] \cdot {}^h[\varphi,h'] \equiv  [\varphi, h] [\varphi,h']$ 
modulo $H_{i+j+1}$. Finally, if $H_*$ is an $N$-series, then we can apply the 
Witt--Hall identity from~\cref{par_commutators} to get, for $h \in H_i$ and 
$h' \in H_k$, $[[h, h'],{}^{h'} \varphi] \cdot [[h',\varphi], {}^\varphi h] 
\cdot [[\varphi,h],{}^h h' ]= 1$. The first factor is also 
${}^{h'} [[h, h']^{h'}, \varphi]$. Using the centrality of $H$, together 
with the triviality of the action of $\varphi$ on $\gr(H_*)$ (since 
$j \geq 1$), we see that in $\gr_{i+j+k}(H_*)$, this equality becomes 
$[[h, h'], \varphi]  + [[h',\varphi], h] + [[\varphi,h], h' ] \equiv  0$, 
which is the desired conclusion.
\end{proof}

\begin{remark}
\label{rk_H_psd_Z}
The last part of the proof can also be obtained as an application of the constructions 
of~\cref{par_Johnson_Lie} to $K = \Z$ acting on $H$ through the powers of $\varphi$, 
where $H$ is endowed with the filtration $H_*$ and $\Z$ is endowed with the filtration 
$K_*$ given by $K_k = \Z$ if $k \leq j$ and $K_k = 0$ if $k > j$. Then 
$\varphi \in \mathcal A_j(H_*)$ if and only if $K_*$ acts on $H_*$.
\end{remark}

\begin{remark}
\label{rem:not-central}
Although we will not use this, more can be said in the case when $H_*$ is not central, 
but only a filtration by normal subgroups. For instance, $\tau(\varphi)$ is no longer linear, 
but it still has to satisfy $\tau(\varphi)(hh') = 
\tau(\varphi)(h) \cdot {}^h \tau(\varphi)(h')$.
\end{remark}

\begin{proposition}\label{general_Johnson_morphism}
Let $K_*$ and $H_*$ be central filtrations, and let $K_*$ act on $H_*$ 
(see Definition~\ref{def_action_filtrations}). Then we have a well-defined 
graded linear map:
\[
\tau \colon\left\{\renewcommand{\arraystretch}{1.5}
\begin{array}{clc}
\gr(K_*) &\longrightarrow &\End(\gr(H_*)) \\
\overline k &\longmapsto & \left(\overline h \mapsto \overline{(k \cdot h)h^{-1}}\right),
\end{array}
\right.
\]
where $\End(\gr(H_*))$ is the graded abelian group of homogeneous endomorphisms 
of positive degree. Moreover, if $K_*$ is an $N$-series, then $\tau$ is a morphism 
of Lie rings (for the structure of Lie ring coming from the ring structure on 
$\End(\gr(H_*))$); in other words, $\gr(H_*)$ is then a graded $\gr(K_*)$-module.
\end{proposition}

\begin{remark}
The situation of~\cref{par_Johnson_Lie} is the case where $H_*$ is an $N$-series; 
then $\gr(H_*)$ is a Lie algebra, $\tau$ takes value in derivations, and we recover 
the previous construction.
\end{remark}

\begin{definition}
\label{def:Johnson-hom}
The morphism $\tau\colon \gr(K_*) \to \End(\gr(H_*))$ from Proposition~\ref{general_Johnson_morphism} 
is called the \emph{Johnson homomorphism}\/ associated with the action of $K_*$ on $H_*$. 
\end{definition}

\begin{proof}[Proof of Proposition~\ref{general_Johnson_morphism}]
Recall that $(k \cdot h)h^{-1}$ is the commutator $[k, h]$, computed in $H \rtimes K$. Each 
element $k \in K_j$ acts through an element of $\mathcal A_j(H_*)$, so we can apply 
Lemma~\ref{general_Johnson(phi)} to conclude that the map $h \mapsto (k \cdot h)h^{-1}$ 
(which depends only on the image of $k$ in $\Aut(H)$) induces a well-defined linear 
endomorphism $\tau(k)$ of degree $j$ of $\gr(H_*)$. The latter depends only on the 
class of the automorphism $k \cdot (-)$ modulo $\mathcal A_{j+1}(H_*)$, which depends 
only on the class of $k$ modulo $K_{j+1}$, so the map $\tau$ of the statement is well-defined. 
Moreover, if $k, k' \in K_j$, then $[kk', h] = {}^k [k',h] [k,h]$; the condition 
$[K, H_i] \subseteq H_{i+1}$ implies that the action of $K$ on $\gr(H_*)$ by 
conjugation is trivial, hence $[kk', h] \equiv [k',h] [k,h]$ 
modulo $H_{i+j+1}$, which implies that $\tau$ is linear. 

The only thing left to prove is the last statement. This uses the Witt--Hall identity 
from~\cref{par_filtrations}. Namely, if $k \in K_j$, $k' \in K_k$ and $h \in H_i$, then 
$[[k,k'],{}^{k'} h] \cdot {}^h [[k',h]^h, k] \cdot [[h,k],{}^k k' ]= 1$. 
Using the first part of the proof, and the fact that the actions by conjugation 
do not change the classes considered in $\gr(K_*)$ or $\gr(H_*)$ (because $K_*$ 
acts on $H_*$ in the first case, and because $H_*$ and $K_*$ are central in the 
other two), we deduce that $[[k,k'], h] = [k, [k',h]] - [k',[k,h]]$ in 
$\gr_{i+j+k}(H)$, which is the desired conclusion.
\end{proof}

\begin{example}
\label{gr(G)-modules}
If $G_*$ is any central filtration on a group $G$, then the canonical action 
of $\gamma_*(G)$ on $G_*$ (Example~\ref{ex_LCS_acts_on_central}) induces a 
graded $\gr(G)$-module structure on $\gr(G_*)$. Observe that if $G_*$ is an 
$N$-series, this module structure is the one coming from the morphism 
of Lie rings $\gr(G) \rightarrow \gr(G_*)$ induced by the inclusion 
$\gamma_*(G) \subseteq G_*$, together with the adjoint representation 
of $\gr(G_*)$. Observe further that this graded linear map still exists 
if $G_*$ is only central (we still have $\gamma_*(G) \subseteq G_*$), 
but it is no longer a Lie morphism, since there is no canonical Lie 
structure on its target. It is, however, a $\gr(G)$-equivariant 
map (with respect to the adjoint action of $\gr(G)$ on itself), since 
the module structure on $\gr(G_*)$ is obviously functorial in $G_*$.
\end{example}

With the same proofs as Proposition~\ref{Johnson_inj_Lie} and Corollary~\ref{A_*_and_torsion}, 
where derivations and Lie morphism are replaced by linear maps, we deduce from 
Proposition~\ref{general_Johnson_morphism} the following:

\begin{proposition}
\label{Johnson_inj}
Let $H_*$ be a central filtration on a group $H = H_1$, and let some group $K$ act on $H$ 
by automorphisms. Then the Johnson homomorphism associated to the action of 
$\mathcal A_*(K, H_*)$ on $H_*$ is injective.
\end{proposition}

\begin{corollary}
\label{A_*_and_torsion} 
If a central filtration $H_*$ is torsion-free (resp.~of $q$-torsion), then the filtration 
$\mathcal A_*(K, H_*)$ is an $N_0$-series (resp.~an $N$-series of $q$-torsion).
\end{corollary}

\begin{remark}
\label{rem:weak-J}
The centrality hypothesis on $H_*$ can be weakened, at the cost of losing part of the 
structure. For instance, if we take $H_*$ to be any filtration by normal subgroups 
such that $\gr(H_*)$ is abelian, then $\tau$ is still a well-defined graded linear 
map that becomes injective if $K_* = \mathcal A_*(K, H_*)$, but it does not take 
values in endomorphisms of $\gr(H_*)$ anymore: it takes values in the $H$-derivations 
of the (non-trivial) $H$-module $\gr(H_*)$, and it is no longer a Lie morphism 
when $K_*$ is an $N$-series.
\end{remark}

\section{Decompositions of \texorpdfstring{$\mathcal P$}{P}-filtrations}
\label{sect:P-filtrations}

\subsection{\texorpdfstring{$\mathcal P$}{P}-filtrations}
\label{subsec:P-filtrations-definitions}

Given a property $\mathcal P$ of filtrations (such as being central, 
torsion-free, of $q$-torsion, etc.), we call \emph{$\mathcal P$-filtrations}\/ 
the filtrations satisfying $\mathcal P$. Let us fix $\mathcal P$ and $\mathcal Q$ 
two properties of filtrations, such that $\mathcal P \Rightarrow \mathcal Q$, 
which means that $\mathcal P$-filtrations are in particular $\mathcal Q$-filtrations. 
We refer to $\mathcal Q$ as our \emph{context}. In our main applications, 
$\mathcal Q$ will be a centrality property (being central, or being an 
$N$-series), and $\mathcal P$ will be obtained from $\mathcal Q$ by adding some 
property concerning the behavior of powers with respect to the filtration 
(being torsion-free, or of $q$-torsion, or $p$-restricted). The goal of 
this section is to investigate the consequences of the following hypotheses:
\begin{enumerate}[itemsep = 1.5pt, label=(H\arabic*), wide=1em, leftmargin=*, start=0]
\item \label{H0}
If $f\colon K \to G$ is a group morphism and $G_*$ is a $\mathcal P$-filtration 
on $G$, then $f^{-1}(G_*)$ is a $\mathcal P$-filtration on $K$.
\item \label{H1}
If $G_*$ is a $\mathcal P$-filtration, then $\mathcal A_*(G_*)$ is too.
\item \label{H2}
If $K_*$ acts on $H_*$ and $K_*$ and $H_*$ satisfy $\mathcal P$, 
then so does $H_* \rtimes K_*$.
\item \label{H3}
For every $\mathcal Q$-filtration $G_*$, there is a smallest 
$\mathcal P$-filtration $G_*^{\mathcal P}$ containing $G_*$.
\item \label{H4}
For every $\mathcal Q$-filtration $G_*$, we have 
$\mathcal A_*(G_*) \subseteq \mathcal A_*(G_*^{\mathcal P})$.
\end{enumerate}
We call $G_*^{\mathcal P}$ the \emph{$\mathcal P$-envelope} of the filtration $G_*$.

\begin{remark}
\label{rk_restricting_context}
It is clear that these hypotheses remain true if we restrict our context. Namely, if 
$\mathcal P \Rightarrow \mathcal Q' \Rightarrow \mathcal Q$, then the hypotheses hold 
in the context $\mathcal Q'$ if they hold in the context $\mathcal Q$, and 
$G_*^{\mathcal P}$ is the same in both contexts. 
\end{remark}

\begin{remark}
The first three hypotheses depend only on $\mathcal P$. They are also clearly 
stable under conjunction: if $\mathcal P = \mathcal P_1 \wedge\mathcal P_2$, and if they hold for $\mathcal P_1$ and for $\mathcal P_2$, 
then they hold for $\mathcal P$. 
\end{remark}

\begin{remark}
Hypotheses \ref{H0} and \ref{H3} imply that the category of $\mathcal P$-filtrations 
is a reflexive subcategory of the category of $\mathcal Q$-filtrations (where both are 
seen as full subcategories of the category of group filtrations).  
\end{remark}

\begin{remark}
\label{rk_H4_for_action}
Recall that if $K$ acts on a group $G$ endowed with a filtration $G_*$ via a morphism 
$K \to \Aut(G)$, then $\mathcal A_*(K, G_*) = a^{-1}(\mathcal A_*(G_*))$. So if \ref{H0} 
holds, then \ref{H1} implies that $\mathcal A_*(K, G_*)$ is a $\mathcal P$-filtration 
as soon as $G_*$ satisfies $\mathcal P$. Moreover, if \ref{H4} holds, then $\mathcal A_*(K, G_*) \subseteq \mathcal A_*(K, G_*^{\mathcal P})$ for every action of $K$ on $G$. 
\end{remark}

\begin{remark}
\label{rk_normality_for_(H4)}
Hypothesis \ref{H2} (resp.~\ref{H4}) only makes sense ifwhen 
$\mathcal P$-filtrations (resp.~.~$\mathcal Q$-filtrations) are 
filtrations by normal subgroups, which is the hypothesis under 
which the construction of $\mathcal A_*$ works.
\end{remark}

\subsection{Checking the hypotheses}
\label{subsec:checking_hyp}
Let us discuss how such hypotheses can be checked. Notice first that some 
of them can come for free. In particular, when $\mathcal P$ is centrality together 
with some property of the associated graded abelian group, and the latter is stable 
by taking abelian subgroups, taking graded direct sums and taking graded endomorphisms, 
then \ref{H0}--\ref{H2} are satisfied. Indeed, remark first that $f^{-1}(-)$, 
$\mathcal A_*(-)$, and $(-) \rtimes (-)$ turn central filtrations into central 
filtrations. Then, if $f\colon K \to G$ is a group morphism, and $G_*$ is 
a central filtration on $G$, then $f$ induces a linear injection of 
$\gr(f^{-1}(G_*))$ into $\gr(G_*)$, proving \ref{H0}; also, 
$\gr(\mathcal A_*(G_*))$ injects linearly into $\End(\gr(G_*))$ 
(cf.~Proposition~\ref{Johnson_inj}), proving \ref{H1}. 
Finally, if $K_*$ acts on $H_*$, then 
$\gr(H_* \rtimes K_*) \cong \gr(H_*) \oplus \gr(K_*)$ as graded 
abelian groups, proving \ref{H2}. As an instance of this, we have:

\begin{lemma}
\label{easy_hyp_torsion}
Torsion-free (resp.~$q$-torsion) central filtrations satisfy \ref{H0}--\ref{H2}.
\end{lemma}

As for the existence of $\mathcal P$-envelopes \ref{H3}, it will necessitate 
constructions that will depend very much 
on the property considered; a construction by induction will often work. This 
will also be the case for \ref{H4}, but there is still a general pattern for 
this one: it will hold if the canonical map 
$\gr(G_*) \rightarrow \gr(G_*^{\mathcal P})$ is ``surjective enough''. For instance, 
$\gr(G_*^{\mathcal P})$ can have some canonical structure under which it is generated 
by the image of this morphism. Since $\mathcal A_*$ is made of group automorphisms 
acting ``more and more trivially'' on the associated graded, if maps induced by group 
automorphisms must preserve the structure, then acting ``more and more trivially'' on 
$\gr(G_*)$ will imply acting ``more and more trivially'' on $\gr(G_*^{\mathcal P})$, 
which is generated by $\gr(G_*)$. This rough outline could probably be made into a 
precise general statement, but we do not think it worth the effort (at least given 
our goals), and we only invite the reader to keep it in mind as a guide through 
our investigation of particular cases.

\subsection{Centrality properties}
\label{subsection:central-prop}

As a first example, let us take $\mathcal Q$ to be empty (so that 
$\mathcal Q$-filtrations are all filtrations), and $\mathcal P$ 
to be ``being central'' (resp. ``being an $N$-series''). 

\begin{lemma}
\label{lem_centrality_from_nothing}
The hypotheses \ref{H0}--\ref{H3} are satisfied in these two cases.
\end{lemma}

\begin{proof}
The hypothesis \ref{H0} is clear, \ref{H1} comes from Proposition~\ref{A*_N-series}, 
and \ref{H2} follows from Remark~\ref{rk_action_central} and 
Proposition-Definition~\ref{psd_of_N-series}. 
We prove \ref{H3} by using an inductive construction, in both cases. Namely, let 
$G_*$ be a filtration (on $G = G_1$). For central filtrations, let us define 
$K_*$ inductively by $K_1 = G$ and:
\begin{equation}
\label{induction_formula_for_min_central}
K_{n+1} \coloneqq \llangle G_{n+1} \rrangle{} [G, K_n], 
\end{equation}
where $\llangle - \rrangle{}$ denotes the normal closure in $G$. 
By induction, each $K_n$ is a normal subgroup of $G$. In particular, 
$[G, K_n] \subseteq K_n$ and $\llangle G_{n+1} \rrangle{} \subseteq K_n$, 
hence $K_{n+1} \subseteq K_n$. Moreover, the filtration $K_*$ is central by construction, 
and an easy induction shows that if $H_*$ is another central filtration 
containing $G_*$, then $H_* \subseteq K_*$. A similar construction works 
for $N$-series, replacing the formula~\eqref{induction_formula_for_min_central} by:
\begin{equation}
\label{induction_formula_for_min_N}
K_{n+1} \coloneqq \llangle G_{n+1} \rrangle{} 
\prod\limits_{i = 1}^n [K_i, K_{n+1-i}]. 
\end{equation}
Thus, \ref{H3} is satisfied in both cases. 
\end{proof}

As remarked earlier in Remark~\ref{rk_normality_for_(H4)}, the hypothesis \ref{H4} 
only makes sense if we include normality of the subgroups in our context, so let 
us replace our context $\mathcal Q$ by ``being a filtration by normal subgroups''. 
Recall that this does not change the fact that \ref{H0}--\ref{H3} are satisfied 
(see Remark~\ref{rk_restricting_context}). We then have:

\begin{proposition}
\label{prop_(H4)_for_N-series}
The hypothesis \ref{H4} is satisfied for $N$-series in the context 
of filtrations by normal subgroups.
\end{proposition}

\begin{proof}
Let us use notations from the proof of Lemma~\ref{lem_centrality_from_nothing}: 
let $K_* = G_*^{\mathcal P}$. We prove that 
$\mathcal A_j(G_*) \subseteq \mathcal A_j(K_*)$ for all $j \geq 0$ using the 
strategy described in~\cref{subsec:checking_hyp} above. 
 
Since $K_*$ is an $N$-series, $\gr(K_*)$ admits a Lie ring structure 
(Proposition-Definition~\ref{def_Lie_struct}). 
As a Lie ring, $\gr(K_*)$ is generated by the 
image of $\iota_*\colon \gr(G_*) \to \gr(K_*)$ (this makes sense even 
if $\gr(G_*)$ is not abelian; see Remark~\ref{rk_gr_not_ab}). 
This is by the construction of $K_*$ by induction 
(formula~\eqref{induction_formula_for_min_N}, where $\llangle G_{n+1} 
\rrangle{} = G_{n+1}$): a direct induction shows that the smallest 
sub-Lie ring $A$ of $\gr(K_*)$ containing the 
image of $\gr(G_*)$ must contain $\gr_k(K_*)$, for all $k \geq 1$. 
It also follows directly from this formula that an automorphism of $G$ 
preserving the subgroups $G_i$ also has to preserve the subgroups $K_i$, 
whence the conclusion for $j = 0$. 

If $\varphi \in \mathcal A_1(G_*)$, then in particular it preserves $G_*$, 
hence also $K_*$. Then it induces automorphisms $\varphi_*$ of both 
$\gr(G_*)$ and $\gr(K_*)$, with respect to which $\iota_*$ is equivariant; 
the hypothesis $\varphi \in \mathcal A_1(G_*)$ then says that the former 
one is the identity of $\gr(G_*)$. The second one is a Lie algebra 
automorphism, and if $A$ is its set of fixed points, $A$ is then a 
sub-Lie algebra containing $\iota_*(\gr(G_*))$, so $A = \gr_*(K_*)$, 
which means that $\varphi$ acts trivially on $\gr_*(K_*)$. 
In other words, $\varphi \in \mathcal A_1(K_*)$, and the 
case $j = 1$ is proved. 

Let us now assume that the conclusion holds for some fixed 
$j \geq 1$. Let $\varphi \in \mathcal A_{j+1}(G_*)$. Then 
$\varphi \in \mathcal A_j(G_*) \subseteq \mathcal A_j(K_*)$; 
Lemma~\ref{general_Johnson(phi)} shows that $[\varphi, -]$ induces a 
well-defined endomorphism of degree $j$ of both $\gr(G_*)$ and $\gr(K_*)$, 
and it is a Lie algebra derivation for the latter. It is easily seen that 
we have commutative squares:
\begin{equation}
\label{eq:gr-square}
\begin{tikzcd}
\gr_i(G_*) \ar[r, "\iota_*"] \ar[d, "{[\varphi, -]}"] 
&\gr_i(K_*)  \ar[d, "{[\varphi, -]}"] \\
\gr_{i+j}(G_*) \ar[r, "\iota_*"]
&\gr_{i+j}(K_*).
\end{tikzcd}
\end{equation}
The fact that $\varphi \in \mathcal A_{j+1}(G_*)$ means that 
$[\varphi, G_i] \subseteq [\varphi, G_{i+j+1}]$, which means 
exactly that the vertical map on the left is trivial. 
Hence $A \coloneqq \ker([\varphi, -])$ contains $\iota_*(\gr(G_*))$. 
Since $[\varphi, -]$ is a derivation, it is also a sub-Lie algebra 
of $\gr(K_*)$, so $A$ is equal to $\gr(K_*)$. In other words, 
for every $i \geq 1$, we have $[\varphi, K_i] \subseteq K_{i+j+1}$, 
and this means that $\varphi \in \mathcal A_{j+1}(K_*)$, 
finishing the induction and the proof.
\end{proof}

For central filtrations (in the context of filtrations by normal 
subgroups), \ref{H4} is not satisfied. In fact, we now describe an 
example of a filtration $G_*$ such that $\mathcal A_1(G_*) \not\subset 
\mathcal A_1(G_*^{\mathcal P})$. This $G_*$ will be a filtration on 
$G = F_{3,2} \coloneqq F_3/\gamma_3(F_3)$, the free $2$-step nilpotent group 
on three generators $x$, $y$, and $z$. Recall that this group is a central 
extension of $F_3^{\ab} = \gamma_1(F_3)/\gamma_2(F_3) \cong \Z^3$ (with basis 
$\overline x, \overline y, \overline z$) by the central subgroup 
$\mathcal{Z} = \gamma_2(F_3)/\gamma_3(F_3) \cong \Z^3$ 
(with basis $[x,y], [x,z], [y,z]$), which are the degree $1$ and $2$ parts 
of the free Lie ring on three elements. Our goal is the following:

\begin{proposition}
There is a filtration $G_*$ on $G = F_{3,2}$ such that if $K_*$ denotes the 
smallest central filtration containing $G_*$, then $\mathcal A_1(G_*) 
\not\subset \mathcal A_1(K_*)$.
\end{proposition}

\begin{proof}
Let $G_*$ be the filtration on $G_1 = G$ defined by $G_2 = \llangle y,z \rrangle{}$ 
(which is the kernel of the projection $G \surj \Z$ defined by 
$x \mapsto 1$ and $y, z \mapsto 0$), $G_3 =\gamma_2(G_2)$ and $G_k = 1$ if $k > 3$. 
Let $\varphi$ be the automorphism of $G = F_{3,2}$ fixing $y$ and $z$ and sending 
$x$ to $xy$. The proof will be completed in the next two lemmas, which show that 
$\varphi \in \mathcal A_1(G_*) \setminus \mathcal A_1(K_*)$.
\end{proof}

\begin{lemma}
\label{lem:phi-1}
We have $\varphi \in \mathcal A_1(G_*)$.
\end{lemma}

\begin{proof}
Notice first that $\varphi$ stabilizes the subgroup $G_2 = \llangle y,z \rrangle{}$, 
hence it also stabilizes its characteristic subgroup $G_3 = \gamma_2(G_2)$. 
Moreover, $\varphi$ 
acts trivially on $G_1/G_2$, generated by the class of $x$. Then notice that 
$G_2$ is generated by $y$, $z$, $[x,y]$, $[x,z]$, and $[y,z]$. Indeed, the normality 
of $G_2$ in $G$ can be written as $[G,G_2] \subseteq G_2$, showing that $G_2$ contains 
the last three elements, and hence $\mathcal Z$. Moreover, 
$\langle y, z, \mathcal Z \rangle = 
\pi^{-1}(\langle y,z \rangle)$ is normal in $G$, where $\pi$ is the canonical 
projection from $G$ onto $G^{\ab}$. This implies that $G_3 = [G_2, G_2]$ is 
the cyclic central subgroup generated by $[y,z]$. Now $\varphi$ fixes $y$, $z$, 
and $[y,z]$ (in particular, it acts trivially on $G_3/G_4 = G_3$), and for 
$t \in \{y,z\}$, we have:
\[
\varphi([x,t]) = [\varphi(x),\varphi(t)] = [xy,t] = [x,t][y,t],
\]
where the last equality uses the centrality of commutators in $G$. For 
$t = y$, this gives that $[x,y]$ is also fixed. For $t = z$, this says 
that $[x,z]$ is fixed modulo $G_3 = \langle [x,y] \rangle$, and finally 
we see that $\varphi$ acts trivially on $G_2/G_3$, whence the result.
\end{proof}

Now, let us determine the smallest central series $K_*$ containing $G_*$. 
We must have $K_1 = G$, and $K_2 = G_2$ (which contains $[G,K_1] = \mathcal Z$). 
Then $K_3 = G_3 [G, G_2] = \mathcal Z$. Indeed, $[G, G_2]$ is already 
equal to $\mathcal Z$, since it contains $[x,y]$, $[x,z]$, and $[y,z]$. 
Finally, $K_4 = G_4[G,K_3] = 1$.

\begin{lemma}
\label{lem:phi-2}
We have $\varphi \notin \mathcal A_1(K_*)$.
\end{lemma}

\begin{proof}
It is enough to show that $\varphi$ does not act trivially on 
$K_3/K_4 = \mathcal Z$. But this has already been seen in the proof 
of the previous lemma: $\varphi([x,z]) = [x,z][y,z] \neq [x,z]$.
\end{proof}

\begin{remark}
\label{rem:notpreserved}
The reason why the previously described strategy fails here is that even 
if $\gr(G_*^{\mathcal P})$ has some algebraic structure under which it is 
generated by the image of $\gr(G_*)$ (namely the $\gr(G)$-module structure 
from Example~\ref{gr(G)-modules}), the action of group automorphisms do not 
preserve this structure in general (the relevant formula being: $\varphi([g,x]) = 
[\varphi(g), \varphi(x)] \neq [g,\varphi(x)]$).
\end{remark}

\subsection{Consequences of the hypotheses} 
\label{subsec:conseq-hyp}

Let us now investigate the consequences of our hypotheses. The 
first one can seem to be surprising. It is one of the reasons 
why we do not want to assume from the start that our filtrations 
are $N$-series, or even central.

\begin{lemma}
\label{min_filtrations_are_N-series}
Under \ref{H0}--\ref{H4}, if a ${\mathcal Q}$-filtration $G_*$ is central 
(resp.~an $N$-series), then its $\mathcal P$-envelope $G_*^{\mathcal P}$ is too. 
\end{lemma}

\begin{proof}
A filtration $G_*$ is central if and only if $G \subseteq \mathcal A_1(G, G_*)$. 
But \ref{H4}, together with \ref{H0}, implies that $\mathcal A_1(G, G_*) 
\subseteq \mathcal A_1(G, G_*^{\mathcal P})$ (using Remark~\ref{rk_H4_for_action} 
applied to the conjugation action of $G$ on itself). Thus, if $G_*$ is central, 
then $G \subseteq \mathcal A_1(G, G_*^{\mathcal P})$, so $G_*^{\mathcal P}$ 
is central too. Now suppose that $G_*$ is an $N$-series, that is, for every 
$j \geq 1$, we have $G_j \subseteq \mathcal A_j(G, G_*)$ (where $G_1$ 
acts on itself by conjugation). Then \ref{H4} implies that $\mathcal A_*(G, G_*) 
\subseteq \mathcal A_*(G, G_*^{\mathcal P})$. The latter satisfies 
property $\mathcal P$ by \ref{H1} (using Remark~\ref{rk_H4_for_action}), so it 
contains $G_*^{\mathcal P}$. But 
$G_*^{\mathcal P} \subseteq \mathcal A_*(G, G_*^{\mathcal P})$ means 
exactly that $G_*^{\mathcal P}$ is an $N$-series, and we are done.
\end{proof}

\begin{lemma}
\label{psd_of_min_filtrations}
Assume \ref{H0}--\ref{H4}. Let $K_*$, $H_*$ be $\mathcal Q$-filtrations, 
and let $K_*$ act on $H_*$. Then  
$(H_* \rtimes K_*)^{\mathcal P} = H_*^{\mathcal P} \rtimes K_*^{\mathcal P}$.
\end{lemma}

\begin{proof}
That $K_*$ act on $H_*$ means that $K_* \subseteq \mathcal A_*(K, H_*)$. 
But \ref{H4} implies (through Remark~\ref{rk_H4_for_action}) that $\mathcal A_*(K, H_*) 
\subseteq \mathcal A_*(K, H_*^{\mathcal P})$. The latter filtration satisfies property 
$\mathcal P$ by \ref{H1}, so it contains $K_*^{\mathcal P}$. But 
$K_*^{\mathcal P} \subseteq \mathcal A_*(K, H_*^{\mathcal P})$ 
means that $K_*^{\mathcal P}$ acts on $H_*^{\mathcal P}$, so 
\ref{H2} implies that $H_*^{\mathcal P} \rtimes K_*^{\mathcal P}$ 
satisfies $\mathcal P$. Moreover, it certainly contains $H_* \rtimes K_*$, 
so it must contain $(H_* \rtimes K_*)^{\mathcal P}$. To prove the converse, 
we use \ref{H0} to deduce that the intersection of $(H_* \rtimes K_*)^{\mathcal P}$ 
with $K$ (resp.\ with~$H$) satisfies $\mathcal P$. Since it contains 
$K_*$ (resp.~$H_*$), it contains $K_*^{\mathcal P}$ (resp.~$H_*^{\mathcal P}$). 
As a consequence, $(H_* \rtimes K_*)^{\mathcal P}$ contains 
$H_*^{\mathcal P} \rtimes K_*^{\mathcal P}$, and the proof is complete.
\end{proof}

\subsection{The $\mathcal P$-lower central series} 
\label{subsec:p-lcs}
We now introduce $\mathcal P$-lower central series, and we prove 
our main decomposition theorem (Theorem~\ref{dec_of_P-LCS}), which 
concerns their behaviour with respect to semi-direct product decompositions. 
This will be applied in the subsequent sections to the rational LCS 
(\cref{sect:lcs-rational}), the $q$-torsion LCS (\cref{sec_q-torsion}) 
and the $p$-restricted one (\cref{sec_p-restricted}).

Let our context $\mathcal Q$ be such that for every group $G$, $\gamma_*(G)$ is a 
$\mathcal Q$-filtration, and suppose that $\mathcal Q$-filtrations are filtrations by 
normal subgroups. For instance, $\mathcal Q$ can be ``being central'', ``being an 
$N$-series'', or only ``being a filtration by normal subgroups". Let us assume 
further that $\mathcal P$-filtrations are central, and that $\mathcal P$ 
satisfies \ref{H0} and \ref{H3}. If $G_*$ is a 
$\mathcal P$-filtration on a group $G = G_1$, then it is central, so it contains 
$\gamma_*(G)$, which is a $\mathcal Q$-filtration. As a consequence, it 
also contains its $\mathcal P$-envelope $\gamma_*(G)^{\mathcal P}$, 
which we denote by $\gamma_*^{\mathcal P}(G)$. Thus, 
$\gamma_*^{\mathcal P}(G)$ is the minimal $\mathcal P$-filtration on $G$. We call 
it its \emph{$\mathcal P$-lower central series}. We also denote by $\gr_*^{\mathcal P}(H)$ 
its associated graded.

\begin{lemma}
Under \ref{H0} and \ref{H3}, the $\gamma_i^{\mathcal P}$ are functorial subgroups, 
i.e., $\varphi(\gamma^{\mathcal P}_i(G))\subseteq \gamma^{\mathcal P}_i(H)$ for every 
homomorphism $\varphi\colon G\to H$. In particular, for every group $G$, the subgroups 
$\gamma^{\mathcal P}_i(G)$ are fully invariant in $G$.  
\end{lemma}

\begin{proof}
Let $\varphi\colon G\to H$ be a group morphism. Then \ref{H0} ensures that 
$\varphi^{-1}(\gamma^{\mathcal P}_*(H))$ is a $\mathcal P$-filtration on $G$. 
Hence it contains $\gamma^{\mathcal P}_*(G)$, which is the minimal 
$\mathcal P$-filtration on $G$.
\end{proof}

In particular, if $K$ acts on $H$ by automorphisms, then these automorphisms and their 
inverses have to preserve $\gamma^{\mathcal P}_*(H)$, and we get an induced action of 
$K$ on $\gr_*^{\mathcal P}(H)$. 

The following theorem will be applied throughout the paper to decompose $\mathcal P$-lower 
central series of semi-direct products, for various $\mathcal P$, giving several analogues of 
Proposition~\ref{lcs_of_adp} (and in particular of the result of Falk and Randell~\cite{FR}) 
in these contexts. Remark, however, that Proposition~\ref{lcs_of_adp} is used in its proof.

\begin{theorem}
\label{dec_of_P-LCS}
Assume \ref{H0}--\ref{H4}. Then:
\begin{enumerate}[itemsep = 2pt]
    \item \label{i1}
    For every group $G$, the filtration $\gamma_*^{\mathcal P}(G)$ is an $N$-series.
\end{enumerate}
Moreover, if a group $K$ acts on a group $H$ by automorphisms, then:
\begin{enumerate}[itemsep = 2pt]
\setcounter{enumi}{1}
    \item \label{i2}
    $\gamma_*^{\mathcal P}(H \rtimes K) = 
    \gamma_*^K(H)^{\mathcal P} \rtimes \gamma_*^{\mathcal P}(K)$;
    \item \label{i3}
    If $K$ acts trivially on $\gr_*^{\mathcal P}(H)$, then 
    $\gamma_*^K(H)^{\mathcal P} = \gamma_*^{\mathcal P}(H)$.
\end{enumerate}  
In particular, if $K$ acts trivially on $\gr_*^{\mathcal P}(H)$, then 
we have a canonical isomorphism of graded Lie rings: 
$\gr_*^{\mathcal P}(H \rtimes K) \cong \gr_*^{\mathcal P}(H) 
\rtimes \gr_*^{\mathcal P}(K)$.
\end{theorem}

\begin{proof}
Since $\gamma_*(G)$ is an $N$-series (Proposition~\ref{LCS_min}), so is 
$\gamma_*^{\mathcal P}(G)$ by Lemma~\ref{min_filtrations_are_N-series}. 
Then we can apply Lemma~\ref{psd_of_min_filtrations} to the decomposition of 
Proposition~\ref{lcs_of_sdp} to get the second statement. 

The hypothesis 
of the third one means that $K = \mathcal A_1(K, \gamma_*^{\mathcal P}(H))$. 
But $\mathcal A_*(K, \gamma_*^{\mathcal P}(H))$ is the inverse image of 
$\mathcal A_*(\gamma_*^{\mathcal P}(H))$ by the morphism $K \rightarrow \Aut(H)$. 
As a consequence of \ref{H1} and \ref{H0}, it is a $\mathcal P$-filtration, so it 
contains $\gamma^{\mathcal P}_*(K)$. This means that $\gamma^{\mathcal P}_*(K)$ 
acts on $\gamma^{\mathcal P}_*(H)$. But then \ref{H2} implies that 
$\gamma_*^{\mathcal P}(H) \rtimes \gamma_*^{\mathcal P}(K)$ is a 
$\mathcal P$-filtration on $H \rtimes K$, which must contain 
$\gamma^{\mathcal P}_*(H \rtimes K)$. The reverse inclusion is obtained by 
applying \ref{H0} twice: the intersection of $\gamma^{\mathcal P}_*(H \rtimes K)$ 
with $H$ (resp.~with $K$) is a $\mathcal P$-filtration, so it has to contain 
$\gamma_*^{\mathcal P}(H)$ (resp.~$\gamma_*^{\mathcal P}(H)$). Hence, for all 
$i \geq 1$, the group $\gamma^{\mathcal P}_i(H \rtimes K)$ contains the product 
$\gamma_i^{\mathcal P}(H) \gamma_i^{\mathcal P}(K)$. Consequently,  
$\gamma^{\mathcal P}_*(H \rtimes K) = \gamma_*^{\mathcal P}(H) 
\rtimes \gamma_*^{\mathcal P}(K)$, which is equivalent to 
$\gamma_*^K(H)^{\mathcal P} = \gamma_*^{\mathcal P}(H)$, in view 
of the second statement. 

Finally, the statement about Lie algebras follows directly from 
statement \eqref{i3}, using that if an $N$-series $K_*$ 
acts on an $N$-series $H_*$, then $\gr_*(H_* \rtimes K_*) \cong 
\gr_*(H_*) \rtimes \gr_*(K_*)$ (see the beginning of~\cref{par_Johnson_Lie}).
\end{proof}

\begin{remark}
\label{weaker_hyp_for_LCS}
We do not need the full force of \ref{H4} to get these results about 
$\mathcal P$-lower central series. In fact, we only need to assume 
the weaker:
\begin{enumerate}[itemsep = 1.5pt, label=(H\arabic*') , wide=1em,  
leftmargin=*, start=4]
\item \label{H4-prime}
For every $\mathcal Q$-filtration $G_*$, we have 
$\mathcal A_1(G_*) \subseteq \mathcal A_1(G_*^{\mathcal P})$.
\end{enumerate}
Indeed, whenever we use \ref{H4} in the proof of Theorem~\ref{dec_of_P-LCS}, we 
use it to show that, for some action of $K$ on some group $H = H_1$ endowed with a 
$\mathcal Q$-filtration $H_*$, if $\gamma_*(K) \subseteq \mathcal A_*(K, H_*)$, 
then $\gamma_*^{\mathcal P}(K) \subseteq \mathcal A_*(K, H_*^{\mathcal P})$. 
But since we are dealing with minimal $\mathcal P$-filtrations, and since \ref{H0} 
and \ref{H1} imply that $\mathcal A_*(K, H_*^{\mathcal P})$ is a $\mathcal P$-filtration, 
the latter inclusion follows directly from $K \subseteq \mathcal A_1(K, H_*^{\mathcal P})$, 
which is implied by $K \subseteq \mathcal A_1(K, H_*)$ (hence also by 
$\gamma_*(K) \subseteq \mathcal A_*(K, H_*)$) under our weaker hypothesis. 
\end{remark}

The hypothesis of the third statement of Theorem~\ref{dec_of_P-LCS} can often 
be reduced to a weaker one. Precisely, in all our applications, it will be enough 
to check that $K$ acts trivially on $\gr_1^{\mathcal P}(H)$ to get the triviality 
of the action on $\gr_*^{\mathcal P}(H)$ and the decomposition of the 
$\mathcal P$-lower central series that ensues. This is due to Proposition~\ref{engdeg1}, 
mixed with the fact that $\gr_*(H) \rightarrow \gr_*^{\mathcal P}(H)$ will 
be ``surjective enough'' (in the sense alluded to at the end of \cref{subsec:checking_hyp}). 
Indeed, if $\gr_*^{\mathcal P}(H)$ is, in some sense, generated by the image of $\gr_*(H)$, 
since the latter is generated as a Lie ring by $\gr_1(H) = H^{\ab}$, then 
$\gr_*^{\mathcal P}(H)$ will be ``generated in degree one'' too (by the quotient 
$\gr_1^{\mathcal P}(H)$ of $H^{\ab}$).

\section{Torsion-free filtrations}
\label{sect:lcs-rational}

Our main goal in this section is to obtain a decomposition theorem for rational lower central 
series of semi-direct products (Theorem~\ref{dec_of_rat_LCS}). The rational lower 
central series $\gamma_*^{\Q}$ is the minimal torsion-free central filtration on 
a group. With the point of view of \cref{sect:P-filtrations}, it is $\gamma_*^{\mathcal P}$, 
when the property $\mathcal P$ is ``being torsion-free central'' and  $\mathcal Q$ is 
``being central''. We are thus going to deduce Theorem~\ref{dec_of_rat_LCS} from the 
general Theorem~\ref{dec_of_P-LCS} by showing that our hypotheses \ref{H0}--\ref{H4} are 
satisfied in this situation. 

We have already remarked in Lemma~\ref{easy_hyp_torsion} that the only hypotheses 
that we need to check are \ref{H3} and \ref{H4}. Thus, our first aim is to build a minimal 
torsion-free central filtration containing a given central filtration $G_*$. This 
will be achieved in Proposition~\ref{min_torsion-free_filtrations} by taking isolators, 
generalizing slightly a result of Massuyeau~\cite[Lem.~4.4]{Massuyeau}. Then, 
in~\cref{par_torsion-free_graded}, we investigate the associated graded to show 
that \ref{H4} is satisfied, using the strategy described at the end of~\cref{subsec:checking_hyp}.

Finally, we will show in~\cref{sec:pi-free} that all this (and, in particular, 
Theorem~\ref{dec_of_rat_LCS}) can be generalized to filtrations having no torsion 
only at a given family $\pi$ of primes (see in particular Theorem~\ref{dec_of_pi-free_LCS}).

\subsection{Isolators}
\label{subsec:isolator}
For a subset $S\subseteq G$, we let 
\begin{equation}
\label{eq:isolator}
\ssqrt{S}\coloneq \ssqrt[G]{S}=\{g\in G\mid 
\text{$g^m \in S$ for some $m\in \N$} \}
\end{equation}
be the {\em isolator}\/ (or, {\em root set}) of $S$ in $G$. Clearly, 
$S\subseteq \ssqrt{S}$ and $\ssqrt{\!\smash[b]{\ssqrt{S}}}= \ssqrt{S}$. 
Moreover, if $\varphi\colon G\to H$ is a homomorphism and 
$\varphi(S)\subseteq T$, then $\varphi(\ssqrt[G]{S})\subseteq \ssqrt[H]{T}$.  

The isolator of a subgroup of $G$ need not be a subgroup; for instance, 
$\ssqrt[G]{\{1\}}$ is equal to $\Tors(G)$, the set of torsion elements 
in $G$, which is not a subgroup in general. However, we recall the classical:

\begin{proposition}
\label{isolator_subgroup}
Let $G$ be a nilpotent group and $N\triangleleft G$ be a normal subgroup. 
Then $\ssqrt[G]{N}$ is a normal subgroup of $G$.
\end{proposition}

Since we clearly have $\ssqrt[G]{N}=q^{-1}(\Tors(G/N))$, where 
$q\colon G\surj G/N$ is the canonical projection, 
the proposition follows from the weaker statement:

\begin{lemma}
\label{torsion_subgroup}
Let $G$ be a nilpotent group. Then $\Tors(G)$ is a normal subgroup of $G$. 
\end{lemma}

\begin{proof}
The set of torsion elements is clearly stable under conjugation by elements of $G$; 
we need to show that it is a subgroup. Let $H = \langle \Tors(G) \rangle$ be the subgroup 
generated by torsion elements in $G$. Then $H$ is nilpotent since $G$ is, so 
$\gamma_{c+1}(H) = \{1\}$ for some $c \geq 0$. Moreover, $\gr(H)$ is a Lie ring 
generated by $H^{\ab}$, and the latter is generated by torsion elements, so all 
its elements have finite order. We deduce that $\gr(H)$ is a torsion abelian group. 
Now, if $h \in H$, then its order is finite in $\gr_1(H) = H^{\ab}$, so there exists 
$n_1 \geq 1$ such that  $h^{n_1} \in \gamma_2(H)$. Then $h^{n_1}$ gives an element in 
$\gr_2(H)$, which is torsion, so there exists $n_2 \geq 1$ such that  $h^{n_1n_2} \in 
\gamma_3(H)$. This process ends at $h^{n_1 \cdots n_c} \in \gamma_{c+1}(H) = \{1\}$, 
so $h$ is of finite order. Finally, $H \subseteq \Tors(G)$, which means exactly that 
$\Tors(G)$ is a subgroup.
\end{proof}

A subgroup $H\le G$ is said to be {\em isolated in $G$}\/ if $\ssqrt{H}=H$. If $N$ is 
a normal subgroup, then $\ssqrt[G]{N}/N\cong \Tors(G/N)$; in particular, $N$ is isolated 
if and only if $G/N$ is torsion-free. Furthermore, we have the following easy

\begin{lemma}
\label{G/G_i_torsion-free}
Let $G_*$ be a filtration of a group $G$ by normal subgroups. Then $G_*$ is torsion-free 
if and only if the subgroups $G_i$ are isolated in $G$, that is, if and only if the groups 
$G/G_i$ are torsion-free.
\end{lemma}

\begin{proof}
If $G_i$ is isolated in $G$, then it is isolated in $G_{i-1}$, so $G_{i-1}/G_i$ is 
torsion-free. Conversely, if $G_*$ is torsion-free, and $g \in G \setminus G_i$ for some $i$, 
then there exists $j < i$ such that $g \in G_j \setminus G_{j+1}$. Then for all $n \geq 1$, 
$g^n \notin G_{j+1}$, so $g^n \notin G_i$.
\end{proof}

\subsection{A construction by induction} 
\label{par_inductive_construction}
Let $G_*$ be a central filtration on a group $G = G_1$. We can easily build the minimal 
torsion-free central filtration $K_*$ containing $G_*$ by induction. Namely, we must begin 
with $K_1 = G$. Then, if $K_n$ is constructed, we must choose $K_{n+1}$ such that it contains 
$G_{n+1}$ and $[G, K_n]$. In other words, $K_{n+1}$ must contain the product 
$G_{n+1} [G, K_n]$, which is a normal subgroup of $G$, hence also of $K_n$. 
Then $K_n/K_{n+1}$ must be a torsion-free 
quotient of $K_n/(G_{n+1} [G, K_n])$, which is an abelian group (since 
$[K_n, K_n] \subseteq [G, K_n]$). In order for $K_{n+1}$ to be minimal, $K_n/K_{n+1}$ must 
be the largest torsion-free quotient of $K_n/(G_{n+1} [G, K_n])$, that is, the quotient by 
its torsion subgroup. In other words:
\begin{equation}
\label{induction_formula_for_min_torsion-free}
K_{n+1} = \ssqrt[K_n]{G_{n+1} [G, K_n]}.
\end{equation}
Observe that since $G/K_n$ is torsion-free by Lemma~\ref{G/G_i_torsion-free}, taking an  
isolator in $K_n$ is in fact the same as taking an isolator in $G$. By construction, 
$K_*$ is central, torsion-free, containing $G_*$, and it is minimal for these properties. 

\subsection{Minimal torsion-free central filtrations} 
\label{subsec:min-tff}
The previous construction is not the most useful one. We now present another construction, 
extending to central filtrations Massuyeau's construction of $N_0$-series from 
\cite[Lem.~4.4]{Massuyeau}. Namely, if $G_*$ is a central filtration on $G$, 
each quotient $G/G_i$ is nilpotent, and so by Proposition~\ref{torsion_subgroup} 
the set of torsion elements  forms a normal subgroup of $G/G_i$; hence, $\ssqrt{G_i}$ 
is a normal subgroup of $G$; together they form a filtration $\ssqrt{G_*}$ of $G$. 

\begin{proposition}
\label{min_torsion-free_filtrations}
Let $G_*$ be a central filtration on a group $G$. Then $\ssqrt{G_*}$ is the smallest 
torsion-free central filtration on $G$ containing $G_*$.
\end{proposition}

\begin{proof}
We have already seen that it is a descending filtration by normal subgroups. Is is also 
torsion-free, by construction. Moreover, it is minimal among torsion-free filtrations 
containing $G_*$. Indeed, if $K_*$ is such a filtration, then the subgroups $K_i$ are 
isolated (Lemma~\ref{G/G_i_torsion-free}), so $G_i \subseteq K_i$ implies 
$\ssqrt{G_i} \subseteq \ssqrt{K_i} = K_i$.  We are left with showing that 
$\ssqrt{G_*}$ is central, that is, that $[G, \ssqrt{G_i}] \subseteq \ssqrt{G_{i+1}}$ 
for all $i \geq 1$. This means that in the torsion-free nilpotent quotient 
$\Gamma \coloneqq G/\ssqrt{G_{i+1}}$, we want to show that the isolator of $G$ is central. 
We know that $G$ is central modulo $G_{i+1}$, hence its image in $\Gamma$ is too. So the 
proof is finished if we know that the center $\mathcal Z$ of $\Gamma$ contains the isolators 
of its subgroups, which is true if $\mathcal Z$ is isolated, that is, if $\Gamma/\mathcal Z$ 
is torsion-free. But it is a classical result by Mal'cev that if $\Gamma$ is a torsion-free 
nilpotent group, then $\Gamma/\mathcal{Z}$ is too.
\end{proof}

\begin{remark}
\label{rem:Malcev}
The proof of Mal'cev's result is not a very difficult one; the interested reader is 
referred for instance to~\cite[\S2.2,~Cor.~2.11]{Baumslag}, a more general version 
which we quote below as Lemma~\ref{G/Z_pi-free}.
\end{remark}

\subsection{Associated graded}
\label{par_torsion-free_graded}
Let us now investigate the relations between the associated grad\-ed Lie rings $\gr(G_*)$ 
and $\gr(\ssqrt{G_*})$. The inclusion $\iota$ of $G_*$ into $\ssqrt{G_*}$ induces 
a morphism $\iota_*\colon \gr(G_*) \to \gr(\ssqrt{G_*})$ of graded abelian groups. 
A quick look at the definition allows us to give a first description of the kernel 
and the cokernel at each level $i \geq 1$:
\[
\begin{tikzcd}[column sep = 1.3em]
0 \ar[r] & (\ssqrt{G_{i+1}} \cap G_i)/G_{i+1} \ar[r] 
&G_i/G_{i+1} \ar[r, "\iota_*"] &\ssqrt{G_i}/\ssqrt{G_{i+1}} \ar[r] 
&\ssqrt{G_i}/(G_i\ssqrt{G_{i+1}}) \ar[r] &0.
\end{tikzcd}
\]
We see that the kernel of $\iota_*$ is just the torsion subgroup of $G_i/G_{i+1}$, 
and that the cokernel is a quotient of the torsion subgroup $\ssqrt{G_i}/G_i$ of $G/G_i$. 
In particular, both are torsion abelian groups, and hence $\iota_*$ induces a 
canonical isomorphism:
\begin{equation}
\label{gr(G*)_rat}
\gr(G_*) \otimes \Q \cong \gr(\ssqrt{G_*}) \otimes \Q.
\end{equation}
This generalizes part of~\cite[Prop.~7.2]{Bass-Lubotzky}, which dealt with the case when 
$G_* = \gamma_*(G)$. A consequence of this is the fulfillment of our hypothesis \ref{H4}:

\begin{proposition}
\label{(H4)_for_torsion-free}
Let $G_*$ be a central filtration. Then $\mathcal A_j(G_*) \subseteq 
\mathcal A_j(\ssqrt{G_*})$ for all $j \geq 0$.
\end{proposition}

\begin{proof}
It is clear that an automorphism of $G$ preserving the subgroups $G_i$ must preserve 
their isolators  $\ssqrt{G_i}$, whence the conclusion for $j = 0$. Then, we remark that since 
$\gr(\ssqrt{G_*})$ is torsion-free, the canonical morphism from $\gr(\ssqrt{G_*})$ to 
$\gr(\ssqrt{G_*}) \otimes \Q$ is injective. 
Thus, if an automorphism of $G$ acts trivially on $\gr(G_*)$, then it acts trivially on 
$\gr(G_*) \otimes \Q \cong \gr(\ssqrt{G_*}) \otimes \Q$, hence also on $\gr(\ssqrt{G_*})$. 
This proves our result for $j = 1$. Now, let us assume that the conclusion holds for some 
fixed $j \geq 1$. Let $\varphi \in \mathcal A_{j+1}(G_*)$. Then 
$\varphi \in \mathcal A_j(G_*) \subseteq \mathcal A_j(\ssqrt{G_*})$; 
Lemma~\ref{general_Johnson(phi)} shows that $[\varphi, -]$ induces a well-defined 
endomorphism of degree $j$ of both $\gr(G_*)$ and $\gr(\ssqrt{G_*})$. This 
construction is easily seen to be functorial, so we get commutative diagrams:
\begin{equation}
\label{eq:graded-i-squares}
\begin{tikzcd}
\gr_i(G_*) \ar[r, "\iota_*"] \ar[d, "{[\varphi, -]}"] 
&\gr_i(\ssqrt{G_*})  \ar[d, "{[\varphi, -]}"] \\
\gr_{i+j}(G_*) \ar[r, "\iota_*"]
&\gr_{i+j}(\ssqrt{G_*}).
\end{tikzcd}
\end{equation}
The fact that $\varphi \in \mathcal A_{j+1}(G_*)$ means that 
$[\varphi, G_i] \subseteq [\varphi, G_{i+j+1}]$, which means exactly that the 
vertical map on the left is zero. Since the top and bottom arrows induce isomorphisms 
on rationalizations, the vertical map on the right must also be zero after tensoring 
by $\Q$. But this map is a homomorphism between torsion-free abelian groups, so it is 
trivial, which means that $[\varphi,\ssqrt{G_i}] \subseteq \ssqrt{G_{i+j+1}}$, 
hence $\varphi \in \mathcal A_{j+1}(\ssqrt{G_*})$. This finishes the induction.
\end{proof}

\begin{remark}
\label{rem:decomp-rat-lcs}
The reader interested only in the proof of Theorem~\ref{dec_of_rat_LCS} about 
decompositions of rational lower central series (and not in more general filtrations) 
only needs to know that hypothesis \ref{H4-prime} from  Remark~\ref{weaker_hyp_for_LCS} holds 
(and not the full \ref{H4}), which corresponds only to the case $j = 1$ in the previous 
proposition.
\end{remark}

\subsection{Conclusions}
\label{subsec:conclusions_rat}
Now that we have shown that hypotheses \ref{H0}--\ref{H4} are satisfied for 
torsion-free central filtrations (in the context of central filtrations), 
we can draw some conclusions. The first one 
(Lemma~\ref{min_filtrations_are_N-series}) allows us to recover Massuyeau's result.

\begin{proposition}[{\cite[Lem.~4.4]{Massuyeau}}]  
\label{prop:sqrt-N}
Let $G_*$ be an $N$-series. Then $\ssqrt{G_*}$ is the minimal $N_0$-series 
containing $G_*$.
\end{proposition}

The second one (Lemma~\ref{psd_of_min_filtrations}) recovers the decomposition theorem:
\begin{proposition}
\label{psd_of_min_ss_torsion}
Let $K_*$ act on $H_*$. Then $\ssqrt{H_* \rtimes K_*} = \ssqrt{H_*} \rtimes \ssqrt{K_*}$.
\end{proposition}

We can in fact use a shortcut to prove this Proposition directly from the results 
of~\cref{subsec:isolator}. Namely, we can use:
\begin{lemma}
\label{isolators_in_psd}
Let $K$ act on $H$, let $V$ be a subgroup of $K$ and $U$ be a subgroup 
of $H$ stable under the action of $V$ (so that $U \rtimes V$ is a subgroup 
of $H \rtimes K$). Suppose that $\ssqrt{U \rtimes V}$ is a subgroup of 
$H \rtimes K$. Then $\ssqrt{U \rtimes V} = \ssqrt{U} \rtimes \ssqrt{V}$.
\end{lemma}

\begin{proof}
The subgroup $\ssqrt{U \rtimes V}$ clearly contains $\ssqrt{U}$ and $\ssqrt{V}$, 
so it contains their product $\ssqrt{U} \rtimes \ssqrt{V}$. Conversely, 
let $(h,k) \in \ssqrt{U \rtimes V}$. Then there exists $n \geq 1$ such that 
$(h,k)^n \in U \rtimes V$. This implies $k^n \in V$, so $(1,k) \in \ssqrt{U \rtimes V}$. 
As a consequence, $(h, 1) = (h,k) (1,k)^{-1} \in \ssqrt{U \rtimes V}$, 
hence $h \in \ssqrt{U}$. Finally, $(h,k) \in \ssqrt{U} \rtimes \ssqrt{V}$, 
and the proof is complete.
\end{proof}

\begin{proof}[Another proof of Proposition~\ref{psd_of_min_ss_torsion}]
We already know (by applying Lemma~\ref{torsion_subgroup} to the quotient 
$H \rtimes K/H_i \rtimes K_i$) that $\ssqrt{H_i \rtimes K_i}$ is a subgroup 
of $H \rtimes K$, so we can apply Lemma~\ref{isolators_in_psd}. 
\end{proof}

\subsection{The rational lower central series}
\label{subsec:lcs-rat}
The \emph{rational lower central series}\/ $\gamma^{\rat}(G)$ of a 
group $G$ is the minimal torsion-free central filtration on it. 
We denote by $\gr_*^\Q(G)$ its associated graded. Since $\gamma^{\rat}(G)$ 
is central, it has to contain $\gamma_*(G)$, so it is also the 
minimal torsion-free central filtration containing $\gamma_*(G)$, 
which is $\ssqrt{\gamma_*(G)}$. It is also $\gamma_*^{\mathcal P}(G)$ 
if $\mathcal P$ is ``being torsion-free'' (in the context 
of central filtrations). We will apply Theorem~\ref{dec_of_P-LCS} 
to decompose the rational lower central series of a semi-direct product. 
Before doing that, let us remark that 
$\gr_1^\Q(G) = G/\ssqrt{\gamma_2(G)}$ is the quotient $G^{\abf}$ 
of $G^{\ab} = G/\gamma_2(G)$ by its torsion subgroup, and record 
some basic results about $\gr^\Q(G)$.

\begin{remark}
The topologically-minded reader may want to think of $G^{\ab} \otimes \Q$ 
as $H_1(G, \Q)$.
\end{remark}

\begin{proposition}
\label{engdeg1_rat} 
The rational Lie algebra $\gr^\Q_*(G) \otimes \Q$ is \emph{generated in degree $1$}. 
Precisely, it is generated (as a Lie ring) by 
$\gr_1^\Q(G) = G^{\ab} \otimes \Q \cong G^{\abf} \otimes \Q$. 
\end{proposition}

\begin{proof}
Since $\ssqrt{\gamma_*(G)} = \gamma_*^\Q(G)$, the isomorphism~\eqref{gr(G*)_rat} 
applied to $G_* = \gamma_*(G)$ becomes:
\[
\gr_*(G) \otimes \Q \cong \gr_*^\Q(G) \otimes \Q.
\]
The result then follows from the fact that $\gr_*(G)$ is generated 
in degree one as a Lie algebra over $\Z$ (Proposition~\ref{engdeg1}).
\end{proof}

\begin{corollary}
\label{rat_gr_i_are_fg}
If $G$ is finitely generated, then the groups $\gr_i^{\Q}(G)$ are 
finitely generated, torsion-free abelian groups.
\end{corollary}

\begin{proof}
If $G$ is finitely generated, then $G^{\ab} \otimes \Q$ is a finite-dimensional 
$\Q$-vector space. By Proposition~\ref{engdeg1_rat}, for every $i \geq 1$, 
a generating family for $\gr_i(G_*) \otimes \Q$, is given by the finitely 
many brackets of length $i$ of elements of a basis of $G^{\ab} \otimes \Q$. 
Now, for every $i \geq 1$, the abelian group $\gr_i^\Q(G)$ is torsion-free, 
so it injects into $\gr_i^\Q(G) \otimes \Q$, hence its rank is at most 
$\dim_\Q(\gr_i^\Q(G) \otimes \Q)$. Whence the result.
\end{proof}

\begin{corollary}
\label{cor_engdeg1_rat}
Let a group $K$ act on a group $H$ by automorphisms. If  $K$ acts trivially 
on $H^{\abf}$ (or, equivalently, on $H^{\ab} \otimes \Q$), then $K$ acts trivially 
on $\gr_*^\Q(H)$. 
\end{corollary}

\begin{proof}
The action of $K$ on $H$ is by group automorphisms on $H$, so it acts 
by Lie algebra automorphisms on  $\gr_*^\Q(H) \otimes \Q$. 
As a consequence, if its action on $H^{\abf} \otimes \Q$ is trivial 
then, by Proposition~\ref{engdeg1_rat}, its action on 
the whole of $\gr_*^\Q(H) \otimes \Q$ is also trivial. Finally, recall that 
$\gr_*^\Q(H)$ is torsion-free, so it injects into $\gr_*^\Q(H) \otimes \Q$; 
thus, $K$ must act trivially on it, too.
\end{proof}

We now apply Theorem~\ref{dec_of_P-LCS} to the case of torsion-free 
filtrations, using Corollary~\ref{cor_engdeg1_rat} to weaken the hypothesis 
of the third statement, to get:

\begin{theorem}
\label{dec_of_rat_LCS}
Let $G$ be a group. Then:
\begin{enumerate}[itemsep = 2pt]
\item \label{dr1}
$\gamma_*^\Q(G) = \ssqrt{\gamma_*(G)}$ is an $N$-series.
\end{enumerate}
Moreover, if a group $K$ acts on a group $H$ by automorphisms, then:
\begin{enumerate}[itemsep = 2pt]
\setcounter{enumi}{1}
    \item \label{dr2}
    $\gamma_*^{\Q}(H \rtimes K) = \ssqrt{\gamma_*^K(H)} \rtimes \gamma_*^{\Q}(K)$;
    \item \label{dr3}
    If $K$ acts trivially on $H^{\abf}$, then $\ssqrt{\gamma_*^K(H)} = \gamma_*^{\Q}(H)$.
\end{enumerate}  
In particular, if $K$ acts trivially on $H^{\abf}$, then we have a canonical isomorphism 
of graded Lie rings: $\gr_*^\Q(H \rtimes K) \cong \gr_*^\Q(H) \rtimes \gr_*^\Q(K)$.
\end{theorem}

\begin{remark}
The shortcuts indicated in Remark~\ref{rem:decomp-rat-lcs} and after 
Proposition~\ref{psd_of_min_ss_torsion} may be used to write down a 
quite straightforward deduction of \eqref{dr2} and \eqref{dr3} 
(in Theorem~\ref{dec_of_rat_LCS}) from Proposition~\ref{lcs_of_sdp}, 
if we assume \eqref{dr1} known.
\end{remark}

\begin{example}
\label{ex:milnor}
Let $f$ be a complex, homogeneous polynomial of degree $n$ in 
$d$ variables, and let $M=\{z \in \C^d : f(z)\ne 0\}$ be the 
complement of the zero set of $f$. 
As shown by Milnor, the restriction $f\colon M\to \C^*$ is a smooth 
fibration, with fiber $F=f^{-1}(1)$ and monodromy $h\colon F\to F$, 
$h(z)=e^{2\pi \ii/n} z$. If the map $h_*\colon H_1(F;\Q)\to H_1(F;\Q)$ 
is the identity, it follows from the last part of Theorem~\ref{dec_of_rat_LCS} 
that $\gr_{\ge 2}^\Q(\pi_1(F))\cong \gr_{\ge 2}^\Q(\pi_1(M))$.
When $f$ is the product of $n$ distinct linear factors, 
$M$ can be be viewed as the complement of an arrangement 
of $n$ hyperplanes in $\C^d$, defined by the factors of $f$.
For conditions that insure that $h_*=\id$ in that case, we refer 
to \cite{Suciu-roum2024, Suciu-decomp} and references therein.
\end{example}

\subsection{A guide to the literature}
\label{subsec:lit-guide}
The rational lower central series $\gamma^{\rat}(G)$ was introduced by 
Stallings in \cite{Stallings}. He used an inductive definition, by setting 
$\gamma^{\rat}_{1} (G)=G$ and 
\begin{equation}
\label{eq:gamma-q-filtration}
\gamma^{\rat}_{n+1}(G)=\ssqrt{[G,\gamma^{\rat}_{n}(G)]} .
\end{equation}  
This is the same as our formula~\eqref{induction_formula_for_min_torsion-free}, 
taking into account the remark following it, and the fact that 
$\gamma_{n+1}(G) = [G, \gamma_n(G)]$ is contained in $[G,\gamma^{\rat}_{n}(G)]$. 

Variants of this definition have appeared elsewhere in the literature. For instance, our 
formula~\eqref{induction_formula_for_min_torsion-free} is a direct generalization of Bass 
and Lubotzky's description of $\gamma^{\rat}_*(G)$, that they call $\gamma^{\zz} (G)$ 
\cite[Definition 7.1]{Bass-Lubotzky}. Stallings' definition is also straightforwardly 
equivalent to Higman's definition of $G_*^{\rat}$ \cite[\S 10.4]{Hillman}: $G_1^{\rat}=G$ 
and $G_{n+1}^{\rat}$ is the preimage in $G$ of $\Tors(G/[G,G_{n}^{\rat}])$). 

The alternative description of $\gamma^{\rat}_*(G)$ as $\ssqrt{\gamma_*(G)}$ (or, equivalently, 
of $\gamma^{\rat}_i(G)/\gamma_i(G)$ as the torsion subgroup $\Tors(G/\gamma_i(G))$) 
can also be found in several different places, such as \cite[Proposition 7.2(b)]{Bass-Lubotzky} 
or \cite[Appendix A]{FHT}.  This alternative description has also been used as a definition 
by Koberda~\cite{Koberda}, who called it $\gamma^T_*(G)$. 

Finally, an important feature of the filtration $\gamma_*^\Q(G) = \ssqrt{\gamma_*(G)}$ is 
that it coincides with the \emph{rational dimension series} of $G$, defined by 
$D_k^\Q(G) = G \cap (1 + (I_\Q G)^k)$, where $I_\Q G$ is the augmentation 
ideal in the group algebra $\Q G$ (see \cite[Thm.~IV.1.5]{Passi} 
or \cite[Thm.~XI.1.10]{Passman}). An important result related to this 
is Quillen's theorem from~\cite{Quillen}, which says that the associated 
graded $\gr((I_\Q G)^*)$ to the filtration of $\Q G$ by powers of 
$I_\Q G$ is the universal enveloping algebra of $\gr_*^\Q(G) \otimes \Q$.

\subsection{About $\pi$-free filtrations} 
\label{sec:pi-free}
Let $\pi$ be a set of prime numbers. Following \cite[\S0.5]{Baumslag}, we call 
\emph{$\pi$-numbers}\/ the positive integers whose prime factors are in $\pi$. 
A group is called \emph{$\pi$-free} if none of its elements has order a $\pi$-number. 
A subgroup $H$ of a group $G$ is called \emph{$\pi$-isolated in $G$}\/ 
if for every $g \in G$ and every $\pi$-number $n$, $g \notin H$ implies $g^n \notin H$. 
The $\pi$-isolator of a subgroup $H$ of a group $G$ is the set of elements $g \in G$ 
such that $g^n \in H$ for some $\pi$-number $n$. We denote it by $\ssqrt[\pi, G]{H}$, 
or just by $\ssqrt[\pi]{H}$ when $G$ is clear from the context.

\begin{definition}
\label{def:pi-free}
A filtration $G_*$ is called \emph{$\pi$-free}\/ if all the subgroups $G_i$ 
are $\pi$-isolated in $G$.
\end{definition}

If $G_{i+1} \triangleleft G_{i}$ for all $i \geq 1$, this is equivalent 
(by a direct adaptation of the proof of Lemma~\ref{G/G_i_torsion-free}) 
to the quotients $G_i/G_{i+1}$ being $\pi$-free.

Observe that if $\pi$ is the set of all prime numbers, then $\pi$-freeness is just 
torsion-freeness, and $\pi$-isolators are just isolators. We claim that everything 
that we have done in this section can be adapted to $\pi$-free filtrations in a 
straightforward fashion. The key point is that Mal'cev's result holds in fact 
for $\pi$-free nilpotent groups:

\begin{lemma}[{\cite[\S2.2,~Cor.~2.11]{Baumslag}}]
\label{G/Z_pi-free}
If $G$ is a $\pi$-free nilpotent group, then the quotient $G/\mathcal Z(G)$ 
by its center is $\pi$-free. 
\end{lemma}

We now explain how to adapt our results for torsion-free filtrations to this 
more general context. First, let us check that $\pi$-isolators behave like 
isolators. We denote by $\Tors_{\pi}(G)$ the set of $\pi$-torsion elements 
of a group $G$, that is, the set of elements of $G$ having order a $\pi$-number; 
in other words, $\Tors_{\pi}(G) = \ssqrt[\pi, G]{\{1\}}$. A straightforward 
adaptation of the proof of Lemma~\ref{torsion_subgroup} shows that if $G$ is 
nilpotent, then $\Tors_{\pi}(G)$ is a normal subgroup of $G$. Since, in general, 
for $N \triangleleft G$, we have $\ssqrt[\pi, G]{N}= q^{-1}(\Tors_{\pi}(G/N))$, 
where $q\colon G\surj G/N$ is the canonical projection, 
Proposition~\ref{isolator_subgroup} holds for $\pi$-isolators: 

\begin{proposition}
\label{pi-isolator_subgroup}
Let $G$ be a nilpotent group and $N\triangleleft G$ be a normal subgroup. 
Then the $\pi$-isolator $\ssqrt[\pi,G]{N}$ is a normal subgroup of $G$.
\end{proposition}

Let us turn to the construction of minimal $\pi$-free central filtration 
$K_*$ containing a given central filtration $G_*$. The construction by 
induction from~\cref{par_inductive_construction} adapts straightforwardly, 
by replacing the maximal torsion-free quotient of an abelian group by the 
maximal $\pi$-free quotient, which is the quotient by its $\pi$-torsion. We get: 
\begin{equation}
\label{induction_formula_for_min_pi-free}
K_{n+1} = \ssqrt[\pi]{G_{n+1} [G, K_n]}.
\end{equation}

We can also extend Massuyeau's construction to this context, to get the alternative 
description given by just taking $\pi$-isolators.

\begin{proposition}
\label{min_pi-free_filtrations}
Let $G_*$ be a central filtration on a group $G$. Then the $\pi$-isolators 
$\ssqrt[\pi]{G_i}$ form a $\pi$-free central filtration  $\ssqrt[\pi]{G_*}$ on 
$G = G_1$, which is the smallest one containing $G_*$.
\end{proposition}

\begin{proof}
The proof is the same as that of Proposition~\ref{min_torsion-free_filtrations}, 
using the general version of Mal'cev's Lemma, quoted as Lemma~\ref{G/Z_pi-free}.
\end{proof}

The results of~\cref{par_torsion-free_graded} also adapt easily. Namely, 
the morphism from $\gr(G_*)$ to $\gr(\ssqrt[\pi]{G_*})$ admits the same 
description of its kernel and of its cokernel as in the torsion-free 
case (with $\pi$-isolators instead of isolators), so they must be of 
$\pi$-torsion. Then tensoring with the subring $\Z[\pi^{-1}]$ of $\Q$ 
(which is torsion-free, hence flat as a $\Z$-module), we get a canonical 
isomorphism of graded abelian groups,
\begin{equation}
\label{gr(G*)_pi-free}
\gr(G_*) \otimes \Z[\pi^{-1}] \cong \gr(\ssqrt[\pi]{G_*}) \otimes \Z[\pi^{-1}].
\end{equation}

The rest of the proof of \ref{H0}--\ref{H4} for $\pi$-filtrations (in the 
context of central filtrations) is exactly the same as for torsion-free 
filtrations, replacing ``torsion-free'' by ``$\pi$-free'', isolators by 
$\pi$-isolators, and $\Q$ by $\Z[\pi^{-1}]$ everywhere. We spell out 
some consequences, for further reference, namely the generalizations of 
Proposition~\ref{prop:sqrt-N}, Proposition~\ref{psd_of_min_ss_torsion}, and 
Theorem~\ref{dec_of_rat_LCS} to this context: 

\begin{proposition}
\label{prop:pi-sqrt-N}
Let $G_*$ be an $N$-series. Then $\ssqrt[\pi]{G_*}$ is the minimal 
$\pi$-free $N$-series containing $G_*$.
\end{proposition}

\begin{proposition}
\label{psd_of_min_ss_pi-torsion}
Let $K_*$ act on $H_*$. Then $\ssqrt[\pi]{H_* \rtimes K_*} = 
\ssqrt[\pi]{H_*} \rtimes \ssqrt[\pi]{K_*}$.
\end{proposition}

\begin{remark}
\label{rem:shortcut-pi-isolators}
The same shortcut as the one mentioned after Proposition~\ref{psd_of_min_ss_torsion} 
works here, with a straightforward adaptation of Lemma~\ref{isolators_in_psd} 
to $\pi$-isolators.    
\end{remark}

The \emph{$\pi$-free lower central series}\/ $\gamma_*^{\pi^{-1}}(G)$ of a 
group $G$ is the minimal $\pi$-free central filtration on it 
(called ``$\gamma_*^\pi(G)$" in~\cite{Aschenbrenner-Friedl}). We denote by 
$\gr_*^{\pi^{-1}}(G)$ its associated graded.  We have $\gamma_*^{\pi^{-1}}(G) = 
\ssqrt[\pi]{\gamma_*(G)}$, and it is also $\gamma_*^{\mathcal P}(G)$ 
if $\mathcal P$ is ``being $\pi$-free'' (in the context of central filtrations). 
Notice that $\gr_1^{\pi^{-1}}(G) = G/\ssqrt[\pi]{\gamma_2(G)}$ is the quotient 
$G^{\ab}_{\pi^{-1}}$ of $G^{\ab} = G/\gamma_2(G)$ by its $\pi$-torsion subgroup.

\begin{theorem}
\label{dec_of_pi-free_LCS}
Let $G$ be a group. Then:
\begin{enumerate}[itemsep = 2pt]
\item $\gamma_*^{\pi^{-1}}(G) =  \ssqrt[\pi]{\gamma_*(G)}$ is an $N$-series.
\end{enumerate}
Moreover, if a group $K$ acts on a group $H$ by automorphisms, then:
\begin{enumerate}[itemsep = 2pt]
\setcounter{enumi}{1}
    \item $\gamma_*^{\pi^{-1}}(H \rtimes K) = \ssqrt[\pi]{\gamma_*^K(H)} 
    \rtimes \gamma_*^{\pi^{-1}}(K)$;
    \item If $K$ acts trivially on $H^{\ab}_{\pi^{-1}}$, then 
    $\ssqrt[\pi]{\gamma_*^K(H)} = \gamma_*^{\pi^{-1}}(H)$.
\end{enumerate}  
In particular, if $K$ acts trivially on $H^{\ab}_{\pi^{-1}}$, then we have a 
canonical isomorphism of graded Lie rings: $\gr_*^{\pi^{-1}}(H \rtimes K) 
\cong \gr_*^{\pi^{-1}}(H) \rtimes \gr_*^{\pi^{-1}}(K)$.
\end{theorem}

\section{Construction of filtrations}
\label{sec_Dark}

The goal of this section is to describe both the minimal $q$-torsion central 
filtration $G_*^{\qt}$ containing a given central filtration $G_*$, and the 
minimal $N_p$-series $G_*^{\pr}$ containing a given $N$-series $G_*$. 
Here $q \geq 1$ is any integer, whereas $p$ is a prime number. We first give 
constructions by induction (\cref{par_inductive_construction_q_and_p}), 
which are easy to do but do not give a precise enough description of the 
filtrations thus obtained. We then go on to give another construction, 
which uses more advanced tools of commutator calculus (namely, 
Dark's theorem~\cite{Dark, Passi}) to get a more precise 
description of these minimal filtrations. This goal is achieved in 
Corollary~\ref{explicit_min_q} for $q$-torsion central series and in 
Corollary~\ref{explicit_min_p} for $N_p$-series. This generalizes 
Passi's descriptions of the $p$-restricted LCS \cite{Passi}. 
In the $p$-restricted case, the possibility of such a generalization 
had already been observed by Massuyeau~\cite[Lem.~4.6]{Massuyeau}.

These descriptions will be used in the next sections (\cref{sec_q-torsion}, 
resp.~\cref{sec_p-restricted}) to show that $q$-torsion central filtrations 
(resp. $N_p$-series) satisfy the hypotheses from~\cref{sect:P-filtrations}, 
in the context of central filtrations (resp.~of $N$-series). In particular, 
this will lead to decomposition theorems for the minimal $q$-torsion central 
filtration $\gamma_*^{\qt}(G)$ and for the minimal $N_p$-series $\gamma_*^{\pr}(G)$ 
on a group $G$, along a decomposition of $G$ into a semi-direct product.  

\subsection{Constructions by induction}
\label{par_inductive_construction_q_and_p}
We can adapt the construction of~\cref{par_inductive_construction} to build the smallest 
$q$-torsion central filtration $K_*$ containing a given central filtration $G_*$. Namely, 
$K_1$ must equal $G = G_1$. Then if $K_n$ is constructed for some $n \geq 1$, the subgroup 
$K_{n+1}$ must contain $G_{n+1}[G, K_n]$, so $K_n/K_{n+1}$ must be a quotient of $q$-torsion 
of the abelian group $A \coloneqq K_n/G_{n+1}[G, K_n]$. In order for $K_{n+1}$ to be minimal, 
$K_n/K_{n+1}$ must be the maximal $q$-torsion quotient of $A$, which is $A/qA$. In other 
words, we get:
\begin{equation}
\label{induction_formula_for_min_q-torsion}
K_{n+1} = K_n^q G_{n+1}[G, K_n].
\end{equation}
By construction, $K_*$ is a $q$-torsion central filtration containing $G_*$, and it is the 
smallest such filtration.

\begin{remark}
\label{rem:subgrou-subset}
Recall that $K_n^q$ denotes the \emph{subgroup}\/ generated by the elements of the 
form $g^q$ with $g \in K_n$, but here we can in fact take only the \emph{subset}\/ 
of such elements, for $K_{n+1}$ to be the inverse image of $qA$ by $K_n \surj A$ 
(hence a subgroup of $G$).
\end{remark}

Notice that if we apply this construction to $G_* = \gamma_*(G)$, then $K_*$ is the minimal 
$q$-torsion central filtration containing $\gamma_*(G)$, but every central filtration contains 
$\gamma_*(G)$, so it is just the minimal $q$-torsion central filtration on $G$, denoted 
$\gamma^{\qt}_*(G)$. Moreover, $\gamma_n \subseteq K_n$ implies $\gamma_{n+1} \subseteq 
[G, K_n] \subseteq K_{n+1}$, so in this case, the factor $G_{n+1}$ in 
\eqref{induction_formula_for_min_q-torsion} is redundant, and we get the 
inductive description:
\begin{equation}
\label{induction_formula_for_LCS_q-torsion}
\gamma^{\qt}_{n+1} = (\gamma^{\qt}_n)^q [G, \gamma^{\qt}_n].
\end{equation}

We can also devise an inductive construction of the smallest $p$-restricted $N$-series $K_*$ 
containing a given $N$-series $G_*$, which looks somewhat like the previous ones, albeit a 
little bit more complex. Namely, $K_1$ must be equal to $G = G_1$. Then, let us assume that 
$K_1, \ldots, K_n$ are constructed, for some $n \geq 1$. Then $K_{n+1}$ must contain 
$G_{n+1}$, the subgroups $[K_k, K_{n+1-k}]$ for $1 \leq k \leq n$, and all the subgroups 
$K_i^{p^j}$ such that $i p^j \geq n+1$. These are all normal subgroups of $G = G_1$ 
(by induction hypothesis), and they are contained in $K_n$, so their product is too. 
Hence, we can define $K_{n+1}$ by:
\begin{equation}
\label{induction_formula_for_min_Np-series}
K_{n+1} = G_{n+1} \cdot \prod\limits_{k = 1}^n [K_k, K_{n+1-k}] \cdot 
\prod_{\substack{i p^j \geq n + 1 \\ i \leq n}} K_i^{p^j} .
\end{equation}
This infinite product makes sense (as the union of all the products 
of a finite number of factors), but since $K_i^{p^{j+1}} \subseteq K_i^{p^j}$, 
we can as well see it as a finite product, by considering only, for each 
$i \leq n$, the minimal $j$ such that $i p^j \geq n + 1$. 
By construction, $K_*$ is a $p$-restricted $N$-series (that is, an 
$N_p$-series) containing $G_*$, and it is the smallest such filtration. 

\subsection{Two consequences of Dark's theorem}
\label{subsec-Dark}
Dark's theorem is a very useful result concerning formal interactions of exponents 
with commutator calculus. The theorem and its proof may be found 
in~\cite[IV, Thm.~1.11]{Passi}. We list two consequences.

\begin{proposition}[{\cite[Cor.~IV.1.16]{Passi}}]
\label{Dark_for_commutators} 
There exists a unique map $\theta \colon (\N^*)^2 \to F_2 = \langle x,y \rangle$ 
such that
\[
[x^\alpha,y^\beta] = \prod\limits_{r,s \geq 1} \theta(r,s)^{\binom{\alpha}{r}\binom{\beta}{s}}
\]
for all $\alpha, \beta \in \N$, and such that each $\theta(r,s)$ is a product of 
$\{x^{\pm 1},y^{\pm 1}\}$-commutators in which $x^{\pm 1}$ appears at least 
$r$ times and $y^{\pm 1}$ appears at least $s$ times.
\end{proposition}

\begin{remark}
\label{theta(1,1)}
The map $\theta$ is independent of $\alpha$ and $\beta$, and we can use that to compute it. 
For instance, $\theta(1,1) = [x,y]$, as one sees by taking $\alpha = \beta = 1$.
\end{remark}

The second consequence is a formula originally due to P.~Hall \cite{Hall34}: 
\begin{proposition}
\label{Dark_for_products} 
There exists a unique map $\theta \colon \N \to F_2 = \langle x,y \rangle$ 
such that
\[
x^\alpha y^\alpha = \prod\limits_{r \geq 0} \theta(r)^{\binom{\alpha}{r}}
\]
for all $\alpha \in \N$, and such that each $\theta(r)$ is a product of 
$\{x^{\pm 1},y^{\pm 1}\}$-commutators of length at least $r$.
\end{proposition}

\begin{remark}
\label{theta(1)}
Again, the map $\theta$ is independent of $\alpha$, and we can use that to compute it. 
For instance, $\theta(0) = 1$, and $\theta(1) = xy$.
\end{remark}

Of course, these formal formulas, which hold in the free group, are meant to be applied 
to elements $x,y \in G$, for any group $G$ (using the evaluation map $F_2 \to G$ 
sending $x$ to $x$ and $y$ to $y$).

\subsection{Minimal $q$-torsion central filtrations}
\label{subsec:min-q-torsion}
Let $F \colon \N^* \times \N \rightarrow \N \cup \{\infty\}$ 
be a fixed map. Given a filtration $G_*$ by normal subgroups, let us put:
\[
G_r^F \coloneqq \prod\limits_{F(i,\alpha) \geq r} G_i^\alpha,
\]
where the infinite product is the (filtered) union of all its finite subproducts. 
Recall that if $G$ is a group, $G^\alpha$ is the \emph{subgroup} generated by the 
$g^\alpha$ for $g \in G$. If $G$ is normal in some larger group, then clearly 
$G^\alpha$ is too. Recall also that for normal subgroups $H$ and $K$ of a group 
$G$, $HK$ is a normal subgroup and $HK = KH$, so the order in the above product 
does not matter. Notice further that $G_{r+1}^F \subseteq G_r^F$, even if $G_*$ 
is not even a filtration, but any family of subgroups indexed by the integers.

We now generalize \cite[Cor.~IV.1.18]{Passi}, whose construction works in a much 
more general setting than just for $G_* = \gamma_*(G)$. We even have a result 
valid for every central filtration:

\begin{proposition}
\label{Central_filtrations_Dark}
Assume that $F$ satisfies
\begin{equation}
\label{centrality_of_G^F}
F \left(i + r , \tbinom{\alpha}{r}\right) \geq 1 + F(i,\alpha),\ 
\text{ for all $i, r, \alpha \geq 1$ such that $1 \leq r \leq \alpha$.}
\end{equation}
Then for any central filtration $G_*$, the filtration $G_*^F$ is central too. 
Moreover, if there is an integer $q \geq 1$ such that $F(i,q\alpha) 
\geq F(i,\alpha) + 1$ for every $i, \alpha \geq 1$, then $G_*^F$ 
is of $q$-torsion.
\end{proposition}

\begin{proof}
In order to prove the first statement, it is enough to show that if $x \in G$ 
and $y \in G_i$ with $F(i, \alpha) \geq n$, then  $[x, y^\alpha] \in G_{n+1}^F$. 
Indeed, by repeated use of the formula $[a, bc] = [a,b] ({}^b[a,c])$, and by 
normality of $G_{n+1}^F$, we then deduce that 
$[x, G_i^\alpha] \subseteq G_{n+1}^F$ as soon as $F(i, \alpha) \geq n$, and 
then that $[x, G_n^F] \subseteq G_{n+1}^F$, whence the centrality of $G_*^F$. 
Let us apply the formula from Proposition~\ref{Dark_for_commutators} to get: 
\[
[x,y^\alpha] =  \prod\limits_{r \geq 1} \theta(1,r)^{\binom{\alpha}{r}}.
\]
The centrality of $G_*$ implies that $\theta(1,r) \in G_{i + r}$. Indeed, 
Corollary~\ref{LCS_acts_on_central} ensures that any commutator of length 
at least $r+1$ with at least one occurrence of an element of $G_i$ is in 
$G_{i + r}$, and $\theta(1,r)$ is a product of such elements. Now, our 
hypothesis on $F$ says that $G_{i+r}^{\binom{\alpha}{r}} \subseteq G_{n+1}^F$ 
as soon as $F(i, \alpha) \geq n$. Thus, $[x,y^\alpha] \in  G_{n+1}^F$, and the 
centrality of $G_*^F$ is proved.

The abelian group $G_m^F/G_{m+1}^F$ is generated by the classes $\overline{x^\alpha}$, 
with $x \in G_i$ and  $F(i,\alpha) \geq m$. If $F(i,q\alpha) \geq F(i,\alpha) + 1$ 
(for some fixed integer $q \geq 1$), then $x^{q\alpha} \in G_{m+1}^F$, which means 
that $q\overline{x^\alpha} = 0$. Hence $G_m^F/G_{m+1}^F$ is of $q$-torsion, as claimed.
\end{proof}

\begin{corollary}
\label{explicit_min_q}
Let $q$ by an integer and $G_*$ be a central filtration. Then the minimal 
$q$-torsion central filtration containing $G_*$ is given by:
\[
G^{\qt}_* = \prod\limits_{i+j \geq *}G_i^{q^j}.
\]
\end{corollary}

\begin{proof}
Let us define the \emph{$q$-adic valuation} $\nu_q(n)$ of an integer $n \geq 1$ as the 
greatest integer $v$ such that $q^v$ divides $n$, and let us pose:
\[
F(i, \alpha) \coloneqq i + \nu_q(n).
\]
Then the filtration of the statement is $G_*^F$. Thus, we want to show that 
$G_*^F$ is central and of $q$-torsion; its minimality is then clear: if $K_*$ 
is a $q$-torsion central filtration, then $K_i^{q^j} \subseteq K_{i+j}$, 
hence $K_*$ contains $G_*^F$ if it contains $G_*$. Let us check the hypotheses 
of Proposition~\ref{Central_filtrations_Dark} for this particular $F$. 
In order to do this, we first remark that:
\begin{equation}
\label{inégalité_binôme}
\nu_q \bigl( \tbinom{\alpha}{r} \bigr) = 
\nu_q \bigl( \tfrac{\alpha}{r} \tbinom{\alpha-1}{r-1} \bigr) \geq 
(\nu_q(\alpha) - \nu_q(r))_+,
\end{equation}
where $n_+ = \max(n,0)$. As a consequence, $F \left(i + r , \binom{\alpha}{r}\right) 
\geq  i + r + \nu_q(\alpha) - \nu_q(r) \geq i + \nu_q(\alpha) + 1 =  F(i,\alpha) +1$. 
Since $F(i, q\alpha) = F(i, \alpha) + 1$, our conclusion follows from 
Proposition~\ref{Central_filtrations_Dark}.
\end{proof}

\subsection{Minimal $p$-restricted filtrations}
\label{subsec:min-p-rest}
The previous construction does not carry over to the $p$-restricted 
case without change. Or course, given a prime $p$, we can pose:
$F(i, \alpha) \coloneqq i p^{\nu_p(\alpha)},$
and we see immediately that $G_*^F$ has to be contains in 
any $p$-restricted filtration containing $G_*$. However, $F$ 
does not satisfy the hypothesis~\eqref{centrality_of_G^F}: 
\begin{example}
Suppose that $p$ divides $\alpha$, i.e., that  $\nu \coloneqq \nu_p(\alpha) \geq 1$. Then:
\[
F \bigl(i + p , \tbinom{\alpha}{p}\bigr) = 
(i + p)p^{\nu_p(\tbinom{\alpha}{p})} = (i+p)p^{\nu - 1},
\]
where the last equality is obtained through an easy count of $p$-factors in 
the binomial coefficient. For $i$ large enough, this will be less than 
$F(i,\alpha) = ip^\nu$, so $F$ does not satisfy~\eqref{centrality_of_G^F}.
\end{example}

However, the construction will work for $N$-series, by an application of the next result, 
which generalizes~\cite[Cor~IV.1.18]{Passi} to all $N$-series. 
\begin{proposition}
\label{N-series_Dark}
Let $F \colon \N^* \times \N \to \N \cup \{\infty\}$ be a map. Suppose that
\[
F\bigl(ri+sj, \tbinom{\alpha}{r}\tbinom{\beta}{s} \bigr) \geq F(i,\alpha) + F(j,\beta),
\]
for all $i,j,r,s, \alpha, \beta \geq 1$ such that $1 \leq r \leq \alpha$ and 
$1 \leq s \leq \beta$. Then, for every $N$-series $G_*$, the filtration $G_*^F$ is 
an $N$-series. Moreover, if there is a prime $p$ such that $F(i,p\alpha) \geq p F(i,\alpha)$ 
for every $i, \alpha \geq 1$, then $G_*^F$ is $p$-restricted.
\end{proposition}

\begin{proof}
As in the proof of Proposition~\ref{Central_filtrations_Dark}, 
we first remark that $[G_m^F, G_n^F] \subseteq G_{n+m}^F$ will follow from 
$[x^\alpha,y^\beta] \in G_{m+n}^F$, for $x \in G_i$, $y  \in G_j$, 
and $\alpha, \beta$ such that $F(i,\alpha) \geq m$ and $F(j,\beta) \geq n$. 
We can now use the formula from Proposition~\ref{Dark_for_commutators}, 
and we notice that, since $G_*$ is an $N$-series, we have 
$\theta(r,s) \in G_{rm+sn}$. Then, the hypothesis on $F$ 
implies immediately:
\[
\theta(r,s)^{\binom{\alpha}{r}\binom{\beta}{s}} \in G_{m+n}^F.
\]
Thus, $[x^\alpha,y^\beta] \in G_{m+n}^F$, as claimed.

Now, let us assume that $F(i,p\alpha) \geq p F(i,\alpha)$ for every 
$i, \alpha \geq 1$, for some prime number $p$. We want to show that 
$G_{pm}^F$ contains $(G_m^F)^p$ (for all $m \geq 1$). In other words, 
we want to show that $G_m^F/G_{pm}^F$ has exponent $p$. Like in the proof of 
Proposition~\ref{Central_filtrations_Dark}, the hypothesis implies that 
$x^{p\alpha} \in G_{pm}^F$, as soon as $x \in G_i$ with $F(i,\alpha) \geq m$, 
so this quotient is generated by elements of order $p$. It is not abelian, 
but it is $(p-1)$-nilpotent, whence the conclusion, thanks to 
Lemma~\ref{lem_Passman} below.
\end{proof}

\begin{lemma}[{\cite[Lem.~XI.1.15]{Passman}}]
\label{lem_Passman} 
Let $p$ be a prime number, and let $G$ be a $(p-1)$-nilpotent group. 
If $G$ is generated by elements of order $p$, then $G$ has exponent 
$p$ (i.e., $G^p=\{1\}$).
\end{lemma}

\begin{proof}
We reason by induction on the nilpotency class $c$ of $G$. If $c = 1$, that is, 
if $G$ is abelian, the result is clear. Let us fix $c < p$, and let us assume 
the conclusion known for groups whose nilpotency class does not exceed $c - 1$. 
Let $G$ be $c$-nilpotent, generated by a subset $S$ of elements of order $p$, 
and let $x,y \in S$. The subgroup $L \coloneqq \langle x, [x,y] \rangle$ 
satisfies $\gamma_2 L \subseteq \gamma_3 G$, since $\gamma_2(L)$ is normally 
generated by $[x,[x,y]]$. An easy induction shows that 
$\gamma_i L \subseteq \gamma_{i+1}G$ if $i \geq 2$, so $L$ is nilpotent 
of class at most $(c-1)$. Since it is generated by $x$ and $x^y$ ($= [x,y]^{-1}x$) ,
both having order $p$, the induction hypothesis gives $L^p =\{1\}$. 
In particular, $[x,y]^p = 1$ holds for every $x,y \in S$. Thus, 
the subgroup $\gamma_2 G$, which is generated by the conjugates of the 
commutators $[x,y]$ with $x,y \in S$ (it is in fact normally generated by 
these elements) satisfies the induction hypothesis, hence $(\gamma_2 G)^p = \{1\}$. 

Now, let us apply Proposition~\ref{Dark_for_products} to two elements $x,y \in G$ 
such that $x^p = y^p = 1$. If $2 \leq r < p$, we have $\theta(r) \in \gamma_2 G$, 
so  $\theta(r)^p = 1$. But $p$ is prime, so it divides $\binom{p}{r}$. Recalling that 
$\theta(0) = 1$ and $\theta(1) = xy$ (Remark~\ref{theta(1)}), we see that the formula 
is reduced to $1 = x^p y^p = (xy)^p \theta(p)$. The element $\theta(p)$ belongs 
to $\gamma_pG$, which is trivial, $G$ being nilpotent of class at most $p-1$. 
As a consequence, $(xy)^p = 1$. We have shown that the set of elements $g \in G$ 
satisfying $g^p = 1$ is a subgroup, whence our claim.
\end{proof}

\begin{corollary}
\label{explicit_min_p}
Let $p$ by a prime number, and let $G_*$ be an $N$-series. Then the minimal 
$N_p$-series containing $G_*$ exists and is given by:
\[
G^{\pr}_* = \prod\limits_{i p^j \geq *}G_i^{p^j}.
\]
\end{corollary}

\begin{proof}
This formula defines a filtration $G_*^F$ corresponding to the function $F$ 
defined above by $F(i, \alpha) \coloneqq i p^{\nu_p(\alpha)}$. We need to 
check that this function satisfies the hypotheses of Proposition~\ref{N-series_Dark}. 
The second hypothesis is clearly satisfied. As for the first one, we can use the 
inequality~\eqref{inégalité_binôme} for $q = p$. Let us denote by $a$, $b$, $u$, 
and $v$ the $p$-adic valuations of, respectively, $\alpha$, $\beta$, $r$, 
and $s$, with $r = r'p^u$, and $s = s'p^v$. We obtain:
\begin{align*}
F \bigl(ri+sj, \tbinom{\alpha}{r}\tbinom{\beta}{s} \bigr) 
\geq (ri +sj)p^{(a-u)_+}p^{(b-v)_+} &\geq r'p^a \cdot 
ip^{(b-v)_+} + s'p^b \cdot j p^{(a-u)_+} \\
  &\geq ip^a + jp^b = F(i,\alpha) + F(j,\beta).
\end{align*}
Thus, we can apply Proposition~\ref{N-series_Dark} to conclude that $G_*^F$ is an 
$N_p$-series, which is clearly minimal among the ones containing $G_*$. 
\end{proof}

\begin{remark}
More generally, it $t \geq 1$ is an integer, the function $F$ defined by 
$F(i,\alpha) = i p^{\lfloor\frac{v_p(\alpha)}{t}\rfloor}$ satisfies the 
first hypothesis of Proposition~\ref{N-series_Dark}, but only 
$F(i,p^t\alpha) \geq pF(i,\alpha)$. Therefore, $G^F_*$ is the smallest 
$N$-series $K_*$ containing $G_*$ and satisfying $K_i^{p^t} \subseteq K_{pi}$.
\end{remark}

\section{Filtrations of \texorpdfstring{$q$}{q}-torsion}
\label{sec_q-torsion}

Let $q \geq 1$ be any integer (not necessarily a prime number). We prove 
a decomposition theorem for the $q$-torsion LCS of semi-direct products 
(Theorem~\ref{dec_of_q-LCS}). The $q$-torsion LCS $\gamma_*^{\qt}$ is the 
minimal $q$-torsion central filtration on a group, so it is $\gamma_*^{\mathcal P}$, 
if the property $\mathcal P$ is ``being $q$-torsion central'', in the context 
of central filtrations. Thus, Theorem~\ref{dec_of_q-LCS} will follow from 
Theorem~\ref{dec_of_P-LCS} once we show that our hypotheses \ref{H0}--\ref{H4} 
are satisfied for this property. 
We already know that \ref{H0}--\ref{H2} hold (Lemma~\ref{easy_hyp_torsion}). 
Moreover, in the previous section we have constructed in two different ways 
the minimal $q$-torsion central filtration $G_*^{\qt}$ containing a given 
central filtration $G_*$, so \ref{H3} is satisfied too. 

It remains  to verify that hypothesis \ref{H4} is also satisfied. 
In order to do that, we use the explicit description of the filtration 
$G_*^{\qt}$ from Corollary~\ref{explicit_min_q}. We do this by showing 
that the linear map $\iota_* \colon \gr(G_*) \to \gr(G_*^{\qt})$ 
induced by the inclusion of $G_*$ into $G_*^{\qt}$ is ``surjective enough'' 
(in the sense outlined at the end of \cref{subsec:checking_hyp}). 
More precisely, we show that $\gr(G_*^{\qt})$ is generated by 
$\iota_*(\gr(G_*))$ under the $q$-power operation.

\subsection{The $q$-power operation}
\label{par_q-torsion_graded}
If $G_*$ is a $q$-torsion central filtration, then $\gr(G_*)$ is endowed 
with a degree-one map induced by $(-)^q$. This is a well-defined map 
from $\gr_i(G_*)$ to $\gr_{i+1}(G_*)$: if $g \in G_i$ and $\eta \in G_{i+1}$, 
then $[g, \eta] \in G_{i+2}$ by centrality, so modulo $G_{i+2}$, the 
elements $g$ and $\eta$ commute, and 
$(g \eta)^q \equiv g^q \eta^q \equiv g^q$. We call this map the 
\emph{$q$-power operation}, and we still denote it by $(-)^q$.

\begin{lemma}
\label{gr(G)_gen_gr(Gq)}
Let $G_*$ be any central filtration. The smallest sub-$\Z$-module of 
$\gr(G_*^{\qt})$ containing $\iota_*(\gr(G_*))$ and stable under the 
$q$-power operation is $\gr(G_*^{\qt})$ itself.
\end{lemma}

\begin{proof}
This is a linearisation of the formula from Corollary~\ref{explicit_min_q}. 
Precisely, let $M$ be the submodule of $\gr(G_*^{\qt})$ defined in the 
statement, and let $k \geq 1$, and let $v \in \gr_k(G_*^{\qt})$ be the 
class of some $g \in G_k^{\qt}$. By Corollary~\ref{explicit_min_q}, 
$g$ is a (finite) product of elements of the form $h^{q^j}$, with 
$h \in G_i$ and $i + j \geq k$. Hence $v = \overline g \in \gr_k(G_*^{\qt})$ 
is a sum of elements of the form 
\[
\overline {h^{q^j}} = \overline h^{q^j} \in \iota_i(\gr_i(G_*))^{q^j},
\]
with $i+j = k$ (the factors with $i + j \geq k+1$ give elements of 
$G_{k+1}^{\qt}$, which are trivial in $\gr_k(G_*^{\qt})$). This means 
that $v \in M$, whence the conclusion.
\end{proof}

\subsection{Proving \ref{H4}}
\label{par_(H4)_for_q-torsion}
We can now prove that hypothesis \ref{H4} is satisfied for $q$-torsion 
filtrations among central filtrations:

\begin{proposition}
\label{(H4)_for_q-torsion}
Let $G_*$ be a central filtration. Then $\mathcal A_j(G_*) 
\subseteq \mathcal A_j(G_*^{\qt})$ for all $j \geq 0$.
\end{proposition}

\begin{proof}
It follows directly from the induction 
formula~\eqref{induction_formula_for_min_q-torsion} or from 
Corollary~\ref{explicit_min_q} that an automorphism of 
$G$ preserving the subgroups $G_i$ also has to preserve 
the subgroups $G_i^{\qt}$, whence the conclusion for $j = 0$. 

If $\varphi$ is an automorphism in $\mathcal A_1(G_*)$, 
then in particular it preserves $G_*$, hence also $G_*^{\qt}$. 
Therefore, $\varphi$ induces linear automorphisms $\varphi_*$ 
of both $\gr(G_*)$ and $\gr(G_*^{\qt})$; the hypothesis 
$\varphi \in \mathcal A_1(G_*)$ then says that the first one is the 
identity of $\gr(G_*)$. Let $A$ be the abelian subgroup of elements 
fixed by $\varphi_*$ in $\gr(G_*^{\qt})$. Since 
$\varphi_* \circ \iota_* =\iota_* \circ \varphi_* = \iota_*$, the group 
$A$ contains $\iota_*(\gr(G_*))$. Moreover, $\varphi(g^q) = \varphi(g)^q$ 
implies that $\varphi_*$ commutes with $(-)^q$, so $A$ is 
stable under $(-)^q$. Thus, we can apply Lemma~\ref{gr(G)_gen_gr(Gq)} 
to conclude that $A$ is the whole of $\gr_*(G_*^{\qt})$, which means 
that $\varphi$ acts trivially on $\gr_*(G_*^{\qt})$. In other words, 
$\varphi \in \mathcal A_1(G_*^{\qt})$, and the case $j = 1$ is proved. 

Let us now assume that the conclusion holds for some fixed $j \geq 1$. Let 
$\varphi \in \mathcal A_{j+1}(G_*)$. Then $\varphi \in \mathcal A_j(G_*) \subseteq 
\mathcal A_j(G_*^{\qt})$; Lemma~\ref{general_Johnson(phi)} shows that 
$[\varphi, -]$ induces a well-defined endomorphism of degree $j$ of 
both $\gr(G_*)$ and $\gr(G_*^{\qt})$. As in the proof 
of Proposition~\ref{(H4)_for_torsion-free}, we obtain 
commutative squares:
\[
\begin{tikzcd}
\gr_i(G_*) \ar[r, "\iota_*"] \ar[d, "{[\varphi, -]}"] 
&\gr_i(G_*^{\qt})  \ar[d, "{[\varphi, -]}"] \\
\gr_{i+j}(G_*) \ar[r, "\iota_*"]
&\gr_{i+j}(G_*^{\qt}).
\end{tikzcd}
\]
The fact that $\varphi \in \mathcal A_{j+1}(G_*)$ means that 
$[\varphi, G_i] \subseteq [\varphi, G_{i+j+1}]$, which means 
exactly that the vertical map on the left is zero. 
Therefore, the subgroup $A \coloneqq \ker([\varphi, -])$ 
of $\gr(G_*^{\qt})$ contains $\iota_*(\gr(G_*))$. As in 
the previous case, we are going to show that $A$ is stable 
under the $q$-power operation, which will allow us to apply 
Lemma~\ref{gr(G)_gen_gr(Gq)} to conclude that $A$ is equal to 
$\gr(G_*^{\qt})$, so the vertical map is trivial, and 
$\varphi \in \mathcal A_{j+1}(G_*^{\qt})$ by 
Lemma~\ref{general_Johnson(phi)}. 

Let $g \in G_i^{\qt} \setminus G_{i+1}^{\qt}$ be such that $\overline g \in A$. 
By repeated applications of the formula $[x,yz] = [x,y]\cdot {}^y [x,z]$, 
we get: 
\begin{equation}
\label{[x,y^q]}
[\varphi, g^q] = [\varphi, g] \cdot {}^g [\varphi, g]\ \cdots\ 
{}^{g^{q-1}}[\varphi, g].
\end{equation}
Since $\overline g \in A$, we have $[\varphi, g] = 0 \in \gr_{i+j}(G_*^{\qt})$, 
which means that $[\varphi, g] \in G_{i+j+1}^{\qt}$. By centrality of $G_*^{\qt}$, 
the elements $g^k$ act trivially on $[\varphi, g]$ modulo $G_{i+j+2}^{\qt}$ so, 
modulo $G_{i+j+2}^{\qt}$, we get 
$[\varphi, g^q] \equiv [\varphi,g]^q$. But $[\varphi, g] \in G_{i+j+1}^{\qt}$ 
implies $[\varphi, g]^q \in G_{i+j+2}^{\qt}$, so the class of $[\varphi, g^q]$ 
in  $\gr_{i+j+1}(G_*^{\qt})$ is trivial, hence $\overline g^q \in A$. 
This proves that $A$ is stable under $(-)^q$, whence the conclusion.
\end{proof}

\begin{remark}
\label{rem:decomp-q-lcs}
As in Remark~\ref{rem:decomp-rat-lcs}, the reader interested only in proving 
Theorem~\ref{dec_of_q-LCS} only needs the case $j = 1$ in the previous 
proposition.
\end{remark}

\subsection{Conclusions}
\label{subsec:conclusions_q}
Now that we have shown that hypotheses \ref{H0}--\ref{H4} are satisfied for 
$q$-torsion central filtrations (in the context of central filtrations), 
we can draw some conclusions. Firstly, Lemma~\ref{min_filtrations_are_N-series} gives:
\begin{proposition}  
\label{prop:q-torsion_N-series}
Let $G_*$ be an $N$-series. Then $G_*^{\qt}$ is the minimal $q$-torsion $N$-series 
containing $G_*$.
\end{proposition}
Secondly, Lemma~\ref{psd_of_min_filtrations} gives the decomposition 
theorem for general central filtrations:
\begin{proposition}
\label{psd_of_min_q-torsion}
Let $K_*$ act on $H_*$. Then $(H_* \rtimes K_*)^{\qt} = 
H_*^{\qt} \rtimes K_*^{\qt}$.
\end{proposition}

\subsection{The $q$-torsion lower central series}
\label{subsec:lcs-q}
The \emph{$q$-torsion lower central series}\/ $\gamma_*^{\qt}(G)$ of 
a group $G$ is the minimal $q$-torsion central filtration on it. 
We denote by $\gr_*^{\qt}(G)$ its associated graded. 
Since $\gamma^{\qt}(G)$ is central, it has to contain $\gamma_*(G)$, 
so it is also the minimal $q$-torsion central filtration containing 
$\gamma_*(G)$, which is $\gamma_*(G)^{\qt}$. In other words, it is 
$\gamma_*^{\mathcal P}(G)$ if $\mathcal P$ is ``being $q$-torsion'' 
(in the context of central filtrations). We are going to apply 
Theorem~\ref{dec_of_P-LCS} to decompose the $q$-torsion lower 
central series of a semi-direct product. Before that, let us remark 
that $\gr_1^{\qt}(G) = G/(G^q\gamma_2(G)) = G^{\ab}/qG^{\ab} = 
G^{\ab} \otimes (\Z/q) = H_1(G, \Z/q)$. 

\begin{lemma}
\label{engdeg1_q}
Let $G$ be a group. Then $\gr_*^{\qt}(G)$ is generated, as a Lie algebra together with the 
$q$-power operation, by its degree $1$ piece, which is $G^{\ab} \otimes (\Z/q)$.
\end{lemma}

\begin{proof}
The inclusion $\gamma_*(G) \subseteq \gamma_*^{\qt}(G)$ induces a morphism
$\iota_* \colon \gr(G) \to \gr_*^{\qt}(G)$. By Lemma~\ref{gr(G)_gen_gr(Gq)}, its image 
generates $\gr_*^{\qt}(G)$ under the $q$-power operation. By Proposition~\ref{engdeg1}, 
$\gr(G)$ is generated, as a Lie ring, by its degree~$1$ piece. Since the image of 
$\gr_1(G) = G^{\ab}$ by $\iota_*$ is $\gr_1^{\qt}(G) = G^{\ab} \otimes (\Z/q)$, 
the conclusion follows.
\end{proof}

As an immediate consequence, we get:
\begin{corollary}\label{triv_action_q}
Let a group $K$ act on a group $H$ by automorphisms. Suppose that $K$ acts trivially 
on $H^{\ab} \otimes (\Z/q)$. Then $K$ acts trivially on $\gr_*^{\qt}(H)$. 
\end{corollary}

We now apply Theorem~\ref{dec_of_P-LCS} to the case of $q$-torsion filtrations, 
using Lemma~\ref{engdeg1_q} to weaken the hypothesis of the third statement. 
We obtain the following theorem, whose last statement is a mod-$q$ version 
of the Falk--Randell Theorem, generalizing to not necessarily prime integers 
the mod-$p$ version of Bellingeri and Gervais~\cite[Thm.~3.3]{Bellingeri-Gervais}, 
in a more conceptual framework.

\begin{theorem}
\label{dec_of_q-LCS}
Let $G$ be a group. Then:
\begin{enumerate}[itemsep = 2pt]
\item $\gamma_*^{\qt}(G) = \gamma_*(G)^{\qt}$ is an $N$-series.
\end{enumerate}
Moreover, if a group $K$ acts on a group $H$ by automorphisms, then:
\begin{enumerate}[itemsep = 2pt]
\setcounter{enumi}{1}
    \item $\gamma_*^{\qt}(H \rtimes K) = \gamma_*^K(H)^{\qt} \rtimes 
    \gamma_*^{\qt}(K)$;
    \item If $K$ acts trivially on $H^{\ab} \otimes (\Z/q)$, then 
    $\gamma_*^K(H)^{\qt} = \gamma_*^{\qt}(H)$.
\end{enumerate}  
In particular, if $K$ acts trivially on $H^{\ab} \otimes (\Z/q)$ (which is also 
$H_1(H, \Z/q)$), that is, if $H \rtimes K$ is a $q$-almost direct product, then 
we have a canonical isomorphism of graded Lie rings: 
$\gr_*^{\qt}(H \rtimes K) \cong \gr_*^{\qt}(H) \rtimes \gr_*^{\qt}(K)$, 
which also preserves the $q$-power operation.
\end{theorem}

Variants of this theorem can be obtained using the same techniques. 
For instance, we can recover, with a simpler proof:
\begin{theorem}[{\cite[Thm.~4.1]{Bardakov-Bryukhanov-Neshchadim}}]
\label{BBN_thm}
Let a group $K$ act on a group $H$ by automorphisms. Assume that $K$ acts 
trivially on $H^{\ab} \otimes (\Z/q)$ for some $q \geq 2$. Then 
$\gamma_*(H \rtimes K) \subseteq \gamma_*^{\qt}(H) \rtimes \gamma_*(K)$.
\end{theorem}

\begin{proof}
Under our hypothesis, Theorem~\ref{dec_of_q-LCS} ensures that 
$\gamma_*^{\qt}(H) \rtimes \gamma_*^{\qt}(K) = \gamma_*^{\qt}(H \rtimes K)$ 
is an $N$-series. In other words, $\gamma_*^{\qt}(K)$ acts on $\gamma_*^{\qt}(H)$ 
(see Proposition-Definition~\ref{psd_of_N-series}). Since $\gamma_*(K)$ is the 
minimal $N$-series on $K$, we have, thanks to the statements already quoted, 
$\gamma_*(K) \subseteq \gamma_*^{\qt}(K) \subseteq \mathcal A_*(K, \gamma_*^{\qt}(H))$. 
As a consequence, $\gamma_*^{\qt}(H) \rtimes \gamma_*(K)$ is an $N$-series on 
$H \rtimes K$, which thus contains the minimal $N$-series $\gamma_*(H \rtimes K)$ 
on this group.
\end{proof}

\begin{remark}
\label{rk_on_BBN}
The converse inclusion is clearly not true in general. Indeed, we have 
$\gamma_2(\Z \times \Z) = \{0\} \not\supseteq \gamma_2^{\qt}(\Z) \times 
\gamma_2(\Z) = q\Z \times \{0\}$.
Moreover, the hypothesis is sharp: the result does not hold if the action 
of $K$ on $H^{\ab} \otimes (\Z/q) = H/\gamma_2^{\qt}(H)$ is not trivial. 
In fact, for $k \in K$ and $h \in H$, having $k \cdot h \neq h$ modulo 
$\gamma_2^{\qt}(H)$ means that $[k,h] \notin \gamma_2^{\qt}(H)$, where 
$[k,h]$ is computed in $H \rtimes K$ (that is, $[k,h] = (k \cdot h)h^{-1}$). 
But $[k,h]$ always belongs to $\gamma_2(H \rtimes K)$, so if $K$ does not 
act trivially on $H^{\ab} \otimes (\Z/q)$, then $\gamma_2(H \rtimes K) 
\not\subseteq \gamma_2^{\qt}(H) \rtimes \gamma_2(K)$.
\end{remark}

\subsection{A guide to the literature}
\label{subsec:guide-lit_q}

The $q$-torsion lower central series $\gamma^{\qt}(G)$ was first introduced by Stallings 
in~\cite{Stallings}, where $q$ is assumed prime, but a large part of the paper works 
without this hypothesis (which is not actually used before his Lemma~4.2 and the 
second part of his Theorem~4.3). The assumption that $q$ is prime is also made 
in~\cite{Paris} and in~\cite{Bellingeri-Gervais}, though it is also 
often unnecessary in that context. 

The link between the lower central series and filtrations of algebras 
alluded to in~\cref{subsec:lit-guide} does not generalize to the $q$-torsion 
lower central series, but rather to the $p$-restricted ones 
(see \cref{subsec:guide-lit_p}). Moreover, unlike the $p$-restricted 
Lie algebra structure in the $p$-restricted case, the structure induced 
by $(-)^q$ on the associated graded is not well-understood, even 
when $q$ is prime. Motivations for studying $q$-torsion filtrations 
(and, in particular, $\gamma^{\qt}$) are to be found elsewhere. 
Our three main ones are the following:
\begin{itemize}[itemsep = 1.5pt, topsep = 2pt, leftmargin = 1em]
\item If $G$ is a group with $qG^{\ab} = 0$, then, since $\gr(G)$ 
is generated by $G^{\ab}$ as a Lie ring (see Prop.~\ref{engdeg1}), we must have $q \gr(G) = 0$, 
which means that $\gamma_*(G)$ is of $q$-torsion. Then $\gamma_*(G)$ (the minimal 
central filtration on $G_*$) is equal to $\gamma_*^{\qt}(G)$ (the minimal $q$-torsion 
central filtration on $G_*$), so in this case, studying $\gamma_*^{\qt}(G)$ is actually 
the same as studying $\gamma_*(G)$. This happens for instance when $q = 2$ and $G$ is
a Coxeter group. 
\item As mentioned above, $\gamma^{\qt}$ is linked with group homology~\cite{Stallings}.
\item  $\gamma^{\qt[p]}$ governs the same residual property as the $p$-restricted 
filtration $\gamma^{\pr}$, but it is easier to handle (see Remark~\ref{rem:p-cofinal} 
below; we will call RTN$_p$ this residual property, which coincides with being residually 
a finite $p$-group for finitely generated groups). In particular, its nice inductive 
description \eqref{induction_formula_for_LCS_q-torsion} allows one to apply induction 
arguments; see for instance Theorem~\ref{Residual_prop_of_sdp} and its proof.
\end{itemize}

\section{\texorpdfstring{$p$}{p}-restricted filtrations}
\label{sec_p-restricted}
This section is devoted to the proof of a decomposition theorem for the $p$-restricted 
LCS of semi-direct products (Theorem~\ref{dec_of_p-LCS}). As for the rational 
LCS (\cref{sect:lcs-rational}) and for the $q$-torsion LCS (\cref{sec_q-torsion}), 
the decomposition theorem will follow from Theorem~\ref{dec_of_P-LCS} if we show 
that it can be applied. Namely, we show that our hypotheses \ref{H0}--\ref{H4} 
are satisfied for $N_p$-series (that is, $p$-restricted $N$-series), in the 
context of $N$-series. Our reason for restricting the context to $N$-series only 
(as opposed to central filtrations) is that $p$-restricted central filtrations 
seem to be very ill-behaved in general. For instance, there are examples of 
$p$-restricted central filtrations for which the $p$-power operation on the 
group does not induce a well-defined operation on the associated graded. 

\subsection{Proving \ref{H0}--\ref{H2}}
\label{par_(H0)-(H2)_for_p-restricted}
The hypothesis \ref{H0} is clear. Then \ref{H1} is given by the following:
\begin{proposition}[{\cite[Prop~8.5]{Habiro-Massuyeau}}]
\label{(H1)_p-restricted}
Let $G_*$ be an $N_p$-series. Then $\mathcal A_*(G_*)$ is an $N_p$-series.
\end{proposition}

\begin{proof}
Let $\varphi \in \mathcal A_j(G_*)$ and $g \in G_i$. Using Proposition~\ref{Dark_for_commutators}, 
we get:
\[
[\varphi^p, g] = \prod\limits_{k=1}^p \theta(k,1)^{\binom{p}{k}},
\]
with $\theta(k,1) \in G_{i+kj}$, for all $k$. If $k < p$, then $p$ 
divides $\binom{p}{k}$, hence $\theta(k,1)^{\binom{p}{k}} \in G_{p(i+kj)} 
\subseteq G_{i+pj+1}$. Since $\theta(p,1)$ is also in $G_{i+pj}$, we have:
\[
[\varphi^p, g] \in G_{i+pj}.
\]
This is true for every $g \in G_i$, for all $i$. Hence 
$\varphi^p \in \mathcal A_{pj}(G_*)$, as required.
\end{proof}

As for \ref{H2}, it is a straightforward generalization of~\cite[Prop.~3.10]{Darne1}:
\begin{proposition}
\label{(H2)_p-restricted}
Let $K_*$ and $H_*$ be $N_p$-series. If $K_*$ acts on $H_*$, then the filtration 
$H_* \rtimes K_*$ is an $N_p$-series.
\end{proposition}

\begin{proof}
Denote $H_* \rtimes K_*$ by $G_*$. We already know that it is an $N$-series 
(see Prop.-Def.~\ref{action_of_LCS}).
Let us show that it is $p$-restricted. An element of $G_j = H_j \rtimes K_j$ 
is a product $hk$, with $h \in H_j$ and $k \in K_j$. Using 
Proposition~\ref{Dark_for_products}, we get:
\begin{equation}
\label{eqn_Dark}
h^p k^p = \prod\limits_{i=1}^p \theta(i)^{\binom{p}{i}} = 
(hk)^p \cdot \theta(2)^{\binom{p}{2}} \cdots \theta(p-1)^p \cdot \theta(p),
\end{equation}
with $\theta(1) = hk$ and $\theta(i) \in G_{ij}$ for any $i$, since $G_*$ is 
an $N$-series. We use this formula to show, by induction on $d \leq p$, 
that the inclusion 
\[
G_j^p \subseteq G_{dj}
\]
holds for all $j\ge 1$. 
This is clearly true for $d = 1$. Let us assume that it holds for $d-1$.
Let $hk \in G_j$, as above. For the sake of clarity, let us rewrite 
the formula \eqref{eqn_Dark}:
\[
(hk)^p = h^p k^p \cdot \theta(p)^{-1} \cdot \prod\limits_{i=p-1}^2 
\theta(i)^{-\binom{p}{i}}.
\]
Using, respectively, that $H_*$ is $p$-restricted, that $K_*$ 
is too and that $G_* = H_* \rtimes K_*$ is an $N$-series, we get: 
\[h^p,\ k^p,\ \theta(p)\ \in G_{pj} \subseteq G_{dj},\]
where the inclusion comes from the inequality $d \leq p$.
If $2 \leq i < p$, then $\theta(i) \in G_{ij}$. Since $p$ 
divides $\binom{p}{i}$, the induction hypothesis implies:
\[
\theta(i)^{\binom{p}{i}} \in G_{ij}^p \subseteq G_{(d-1)ij} \subseteq G_{dj},
\]
because $(d-1)ij \geq dj$. Finally, we get what we were looking for:
\[
(hk)^p \in G_{dj}.
\]
This completes the induction step, and the proof of the inclusion 
$G_j^p \subseteq G_{pj}$.
\end{proof}

\subsection{The $p$-power operation}
\label{par_p-restricted_graded}
If $G_*$ is an $N_p$-series, then $\gr(G_*)$ is endowed with maps 
$\gr_i(G_*) \to \gr_{pi}(G_*)$ induced by the $p$-power operation $(-)^p$. 
We use the following result by Lazard, which is a consequence of a deep 
result saying that every $N_p$-series comes from an algebra filtration 
on the group algebra $\mathbb F_p G$:

\begin{proposition}[{\cite[Cor~6.8]{Lazard}}]
\label{gr(p-restr)_is_p-restr}
Let $G_*$ be an $N_p$-series. Then $\gr(G_*)$, endowed with the 
bracket coming from commutators together with the maps induced 
by $(-)^p$, is a $p$-restricted Lie algebra.
\end{proposition}

We use this to give a more precise version of Lemma~\ref{general_Johnson(phi)} 
for $N_p$-series, by adapting the construction of Remark~\ref{rk_H_psd_Z} 
to this context. Recall that a \emph{$p$-restricted derivation}\/ of a 
$p$-restricted Lie algebra $L$, in the sense of Jacobson~\cite{Jacobson} 
is a linear endomorphism $\partial\colon L\to L$ satisfying the identities:
\begin{equation}
\label{eq:p-res-derivation}
\partial([x,y]) = [\partial x, y] +[x, \partial y] 
\quad \text{and} \quad
\partial(x^p) = \ad_x^{p-1}(\partial x).
\end{equation}
In a $p$-restricted Lie algebra, the maps $[x,-]$ are $p$-restricted derivations.

\begin{lemma}
\label{p-restricted_Johnson(phi)}
Let $H_*$ be an $N_p$-series, and let $\varphi \in \mathcal A_j(H_*)$. 
Then the degree-$j$ derivation  $\tau(\varphi)$ of the Lie ring $\gr(H_*)$ 
induced by $[\varphi,-]$ (see Lemma~\ref{general_Johnson(phi)}) is a 
$p$-restricted derivation.
\end{lemma}

\begin{proof}
Let us consider, as in Remark~\ref{rk_H_psd_Z}, the group $K = \Z$ acting on $H$ through 
the powers of $\varphi$, that is, through the homomorphism $a \colon \Z \to \Aut(H)$ 
sending $1$ to $\varphi$. Let $K$ be endowed with the minimal $p$-restricted 
filtration on $\Z$ containing $1$ in degree $j$, which is given explicitly by 
$K_k  = \Z$ for $k \leq j$, and 
$K_k = p^\alpha \Z$ if $p^{\alpha - 1}j < k \leq p^\alpha j$ for $k > j$. 
Proposition~\ref{(H1)_p-restricted} says that $\mathcal A_*(H_*)$, hence 
also $\mathcal A_*(K, H_*) = a^{-1}(\mathcal A_*(H_*))$ is an 
$N_p$-series. Hence $K_*$ acts on $H_*$ (which means that 
$K_* \subseteq \mathcal A_*(K, H_*)$) if and only if 
$\varphi \in \mathcal A_j(H_*)$. 

Under this hypothesis, we can apply Proposition~\ref{(H2)_p-restricted} 
to conclude that $H_* \rtimes K_*$ is an $N_p$-series on $H \rtimes \Z$. 
Then Lazard's result quoted above as Proposition~\ref{gr(p-restr)_is_p-restr} 
implies that $\gr(H_* \rtimes K_*)$ is a $p$-restricted Lie algebra. 
But the latter identifies with  $\gr(H_*) \rtimes \gr(K_*)$, 
where $\gr(K_*)$ is (explicitly) the free $p$-restricted Lie algebra 
generated by the degree-$j$ element $\overline \varphi$. 
Then the degree-$j$ derivation of $\gr(H_*)$ induced by $[\varphi,-]$ 
is none other than the restriction of $[\overline\varphi,-]$, 
computed in  $\gr(H_* \rtimes K_*)$. This is of the form $[x,-]$ in a 
$p$-restricted Lie algebra, so it is a $p$-restricted derivation, as claimed.
\end{proof}

If $G_*$ is any $N$-series, let us consider the Lie morphism $\iota_* \colon \gr(G_*) \to 
\gr(G_*^{\pr})$ induced by the inclusion of $G_*$ into $G_*^{\pr}$. Using exactly the 
same proof as that of Lemma~\ref{gr(G)_gen_gr(Gq)}, we can deduce from 
Corollary~\ref{explicit_min_p} the $p$-restricted analogue:

\begin{lemma}
\label{gr(G)_gen_gr(Gp)}
Let $G_*$ be an $N$-series. The smallest sub-$\Z$-module of $\gr(G_*^{\pr})$ 
containing $\iota_*(\gr(G_*))$ and stable under the $p$-power operation is 
$\gr(G_*^{\pr})$ itself.
\end{lemma}

\subsection{Proving \ref{H4}} 
\label{subsec:proveH4}
Hypothesis \ref{H3} has already been proven in \cref{sec_Dark}, where the minimal 
$N_p$-series containing a given $N$-series has been constructed in two ways. We 
now use the explicit description from Corollary~\ref{explicit_min_p}, together with 
Proposition~\ref{(H2)_p-restricted}, to prove \ref{H4} for $N_p$-series in 
the context of $N$-series:
\begin{proposition}
\label{(H4)_for_p-restricted}
Let $G_*$ be an $N$-series. Then $\mathcal A_j(G_*) \subseteq \mathcal A_j(G_*^{\pr})$ 
for all $j \geq 0$.
\end{proposition}

\begin{proof}
The strategy of proof is the same as for Proposition~\ref{(H4)_for_q-torsion}, 
using the $p$-power operation instead of the $q$-power operation. The proof 
for $j = 0$ and $j = 1$ is exactly the same as in the proof of 
Proposition~\ref{(H4)_for_q-torsion}. The proof of the induction step 
(assuming that the conclusion holds for some fixed $j \geq 1$) also works 
in the same way: if $\varphi \in \mathcal A_{j+1}(G_*)$, we get a well-defined 
endomorphism of degree $j$ of both $\gr(G_*)$ and $\gr(G_*^{\pr})$, with the 
same commutative diagram, and we need to show that $A \coloneqq \ker([\varphi, -])$ 
is stable under the $p$-power operation. But this is a direct consequence of 
Lemma~\ref{general_Johnson(phi)}, which says that $[\varphi, -]$ is a $p$-restricted 
derivation. Indeed, if $[\varphi, x] = 0$, then  $[\varphi, x^p] = \ad_x^{p-1}
([\varphi, x]) = 0$. We can thus apply Lemma~\ref{gr(G)_gen_gr(Gp)} to conclude that 
$A$ is $\gr(G_*^{\pr})$, so that $[\varphi, -] = 0$, and then $\varphi \in \mathcal A_{j+1}
(G_*^{\pr})$ by Lemma~\ref{general_Johnson(phi)}, which is the required conclusion.
\end{proof}

\begin{remark}
\label{rem:decomp-p-lcs}
As for torsion-free filtrations (Remark~\ref{rem:decomp-rat-lcs}) and for $q$-torsion ones 
(Remark~\ref{rem:decomp-q-lcs}), the reader interested solely in the proof of 
Theorem~\ref{dec_of_p-LCS} only needs the case $j = 1$ of the previous proposition.
\end{remark}

The hypotheses \ref{H0}--\ref{H4} are now proved for $p$-restricted filtrations 
(in the context of $N$-series), so we can draw some conclusions. Before turning 
to the ones concerning the $p$-restricted lower central series, let us notice that 
Lemma~\ref{min_filtrations_are_N-series} is not relevant here, since we are already 
dealing only with $N$-series. However, Lemma~\ref{psd_of_min_filtrations} gives 
the following decomposition theorem, for general $N$-series:
\begin{proposition}
\label{psd_of_min_p-restricted}
Let $K_*$ act on $H_*$. Then $(H_* \rtimes K_*)^{\pr} = H_*^{\pr} \rtimes K_*^{\pr}$.
\end{proposition}

\subsection{The $p$-restricted lower central series}
\label{subsec:lcs-p}
For a prime $p$, the \emph{$p$-restricted lower central series} $\gamma_*^{\pr}(G)$ 
of a group $G$ is the minimal $N_p$-series on it. We denote by $\gr_*^{\pr}(G)$ its 
associated graded. Since the filtration $\gamma^{\pr}(G)$ is central, 
it has to contain $\gamma_*(G)$, 
so it is also the minimal $N_p$-series containing $\gamma_*(G)$, which is 
$\gamma_*(G)^{\pr}$. In other words, it is $\gamma_*^{\mathcal P}(G)$ if 
$\mathcal P$ is ``being  $p$-restricted'' (in the context of $N$-series). 
We are going to apply Theorem~\ref{dec_of_P-LCS} to decompose the $p$-restricted 
lower central series of a semi-direct product. Before that, let us remark that 
$\gr_1^{\pr}(G) = G/(G^p\gamma_2(G)) = G^{\ab}/pG^{\ab} = G^{\ab} \otimes \mathbb F_p$.
With the same proof as that of Lemma~\ref{engdeg1_q} (using Lemma~\ref{gr(G)_gen_gr(Gp)} 
instead of Lemma~\ref{gr(G)_gen_gr(Gq)}), we get:

\begin{lemma}
\label{engdeg1_p}
Let $G$ be a group. Then $\gr_*^{\pr}(G)$ is generated, as a $p$-restricted 
Lie algebra, by its degree $1$ piece, which is $G^{\ab} \otimes \mathbb F_p$.
\end{lemma}

We directly deduce:
\begin{corollary}
\label{triv_action_p}
Let $K$ be a group acting by automorphisms on a group $H$. Suppose that 
$K$ acts trivially on $H^{\ab} \otimes \mathbb{F}_p$. Then $K$ acts 
trivially on $\gr_*^{\pr}(H)$. 
\end{corollary}

We now apply Theorem~\ref{dec_of_P-LCS} to the case of $N_p$-series, 
using Corollary~\ref{triv_action_p} to weaken the hypothesis of the third statement. We get:
\begin{theorem}\label{dec_of_p-LCS}
Let a group $K$ act on a group $H$ by automorphisms. Then:
\begin{enumerate}[itemsep = 2pt]
    \item $\gamma_*^{\pr}(H \rtimes K) = \gamma_*^K(H)^{\pr} \rtimes \gamma_*^{\pr}(K)$;
    \item If $K$ acts trivially on $H^{\ab} \otimes \mathbb{F}_p$, then 
    $\gamma_*^K(H)^{\pr} = \gamma_*^{\pr}(H)$.
\end{enumerate}  
In particular, if $K$ acts trivially on $H^{\ab} \otimes \mathbb F_p$, then 
we have a canonical isomorphism of graded $p$-restricted Lie rings: 
$\gr_*^{\pr}(H \rtimes K) \cong \gr_*^{\pr}(H) \rtimes \gr_*^{\pr}(K)$.
\end{theorem}

\subsection{A guide to the literature}
\label{subsec:guide-lit_p}

In general, $p$-restricted filtrations come from filtrations on associative 
algebras over $\F_p$. In fact, it is a classical theorem that every $p$-restricted 
filtration on a group $G$ comes from a filtration on the group algebra 
$\mathbb F_pG$ (see for instance \cite[Chap.~3,~Thm.~1.7]{Passi} 
and the reference therein). As a matter of fact, the $p$-restricted 
lower central series $\gamma^{\pr}(G)$ was introduced by Zassenhaus 
in~\cite{Zassenhaus} as the \emph{mod-$p$ dimension series}\/ of $G$, 
that is, as coming from the filtration by the powers of the augmentation 
ideal of $\mathbb F_pG$ (see also~\cite[Thm.~IV.1.5]{Passi} 
and~\cite[Thm.~XI.1.19]{Passman}). The main result of~\cite{Zassenhaus} 
is actually to show that the mod-$p$ dimension series can also be 
described by the combinatorial formula that we obtain by applying 
Corollary~\ref{explicit_min_p} to the lower central series:
\[
\gamma^{\pr}_* = \prod\limits_{i p^j \geq *}\gamma_i^{p^j}.
\]
This also appears as~\cite[Thm.~5.6]{Lazard} in the more general setting of 
mod-$p^h$ dimension series; see also~\cite[Thm.~XI.1.20]{Passman}. 

All this makes the $p$-restricted lower central series the right analogue 
of the rational lower central series in finite characteristic, which 
explains why it is the filtration appearing in the version of Quillen's 
theorem over $\mathbb F_p$~\cite{Quillen}: the associated graded 
$\gr((I_{\mathbb F_p}G)^*)$ to the filtration of $\mathbb F_p G$ 
by powers of the augmentation ideal $I_{\mathbb F_p} G$ is the 
universal enveloping algebra of the $p$-restricted algebra $\gr_*^{\pr}(G)$.

\section{Residual properties of semi-direct products}
\label{sect:res-semi}
One of the many reasons for studying the lower central series and its rational 
and modular variants comes from the fact that these series control the corresponding 
residual properties of a group $G$. Namely, $G$ is residually nilpotent, respectively, 
residually torsion-free nilpotent, or residually $p$ if and only if the intersection 
of the terms of the series $\gamma_*(G)$, respectively $\gamma_*^{\rat}(G)$, or 
$\gamma_*^{\qt[p]}(G)$ is trivial. 

\begin{remark}
\label{rem:p-cofinal}
If $p$ is prime, the filtrations $\gamma_*^{\qt[p]}$ and $\gamma_*^{\pr}$ are \emph{cofinal}: 
$\gamma_*^{\pr}$ is of $p$-torsion, so $\gamma_*^{\qt[p]} \subseteq \gamma_*^{\pr}$, while 
$\gamma_{\lfloor \log_p(*) \rfloor}^{\qt[p]}$ is an $N_p$-series, so it contains 
$\gamma_*^{\pr}$. As a consequence, although their associated graded behave very 
differently, these two filtrations play the same role when it comes to residual properties.
\end{remark}

Our main goal in this section is to use the machinery introduced above to generalize 
the following easy corollary of Propositions~\ref{lcs_of_sdp} and~\ref{lcs_of_adp}
(partially spelled out by Falk and Randell as~\cite[Thm.~2.6]{FR88}):
\begin{corollary}
\label{cor:res-nilp}
An almost-direct product $H \rtimes K$ is residually nilpotent if 
and only if $H$ and $K$ are.
\end{corollary}

Such a generalization (to residually $\mathcal R$-groups, where $\mathcal R$ is 
nilpotency with some additional conditions) will be obtained as Theorem~\ref{thm:res-P}. 
When the additional conditions concern torsion, we give a further generalization in 
Theorem~\ref{Residual_prop_of_sdp}. These generalizations will be applied mainly 
to residually torsion free nilpotent (RTFN) groups, and to residually $p$-groups.

\subsection{Residual properties and filtrations: a general framework} 
\label{subsec:residual-filtrations}
Let $\mathcal R$ be a class of groups. A group $G$ is said to 
be \emph{residually $\mathcal R$}\/ if for any $g \in G$, $g \ne 1$, 
there exists a group $Q \in \mathcal R$ and an epimorphism 
$\psi \colon G \surj Q$ such that $\psi(g) \ne 1$. We denote 
by $\res(\mathcal R)$ the class of residually $\mathcal R$ groups. 
We now describe a setting in which a given residual property can be 
interpreted nicely in terms of filtrations on a group. Precisely, 
\textbf{for the rest of this section, $\mathcal P$ is a property 
of filtrations that satisfies the following hypotheses:}
\begin{enumerate}[itemsep = 2pt, label=(HR\arabic*) , 
wide=1em,  leftmargin=*, start=0]
\item \label{HR0}
If $f\colon K \rightarrow G$ is a group morphism and $G_*$ is a $\mathcal P$-filtration 
on $G$, then $f^{-1}(G_*)$ is a $\mathcal P$-filtration on $K$;
\item \label{HR1}
If $G_*$ is a $\mathcal P$-filtration, then for any $n \geq 1$, $G_*/G_n$ is a 
$\mathcal P$-filtration on $G_1/G_n$;
\item \label{HR2}
For every group $G$, there is a smallest $\mathcal P$-filtration 
$\gamma_*^{\mathcal P}(G)$ such that $\gamma_1^{\mathcal P}(G) = G$. 
\end{enumerate}

Notice that~\ref{HR0} is exactly~\ref{H0} from \cref{sect:P-filtrations}, and 
that~\ref{HR2} is a weak variant of~\ref{H3}, implied by it  under mild assumptions 
(see the very beginning of \cref{subsec:p-lcs}). Only~\ref{HR1} is really new. 
However, since $\gr(G_*/G_n)$ is a truncation of $\gr(G_*)$ (killing degrees 
$\geq n$), \ref{HR1} is clearly satisfied when $\mathcal P$ is really a condition 
on $\gr(G_*)$ unaffected by truncating, in particular a torsion condition (an 
argument similar to the ones leading to Lemma~\ref{easy_hyp_torsion}). 
Thus we have proved:
\begin{lemma}
\label{checking_(HRi)}
Torsion-free (resp.~$q$-torsion) central filtrations satisfy \ref{HR0}--\ref{HR2}.
\end{lemma}

The associated residual property will be linked to the existence of separating 
$\mathcal P$ filtra\-tions.
\begin{definition}
\label{def_separating}
A filtration $G_*$ is called~\emph{separating} if $\bigcap_{n \geq 1} G_n = \{1\}$.
\end{definition}

Let us define $\mathcal R_{\mathcal P}$ as the class of groups $G$ satisfying 
$\gamma_n^{\mathcal P}(G) = \{1\}$ for some $n \geq 1$. Under the above hypotheses, 
we have:

\begin{lemma}
\label{res_prop_and_separation}
Let $G$ be a group. The following properties are equivalent:
\begin{itemize}[itemsep = 1.5pt]
\item The group $G$ is residually $\mathcal R_{\mathcal P}$;
\item The filtration $\gamma_*^{\mathcal P}(G)$ is separating, that is, 
$\bigcap_{n \geq 1} \gamma_n^{\mathcal P}(G) = \{1\}$;
\item There exists a separating $\mathcal P$-filtration $G_*$ such that $G = G_1$.
\end{itemize}
\end{lemma}

\begin{proof}
The equivalence between the last two statements is clear. Let us proof the equivalence 
between the first two.

Let $x$ be a non-trivial element of $G$, and let $n \geq 1$. 
The statement follows readily from the following claim: $x \notin \gamma_n^{\mathcal P}(G)$ 
if and only if there exists a surjective morphism $\pi\colon G \surj H$ such 
that $\gamma_n^{\mathcal P}(H) = \{1\}$ and $\pi(x) \neq 1$. 

To prove this claim, 
let us first assume that there is such a morphism. Hypothesis \ref{HR0} implies 
that $\pi^{-1}(\gamma_*^{\mathcal P}(H)) \supseteq \gamma_*^{\mathcal P}(G)$. Thus, 
$x \in \gamma_n^{\mathcal P}(G)$ would imply $\pi(x) \in \gamma_n^{\mathcal P}(H) = \{1\}$, 
a contradiction. As a consequence, $x \notin \gamma_n^{\mathcal P}(G)$, and one implication 
is proved. 

Conversely, suppose that $x \notin \gamma_n^{\mathcal P}(G)$. This means that the 
class of $x$ in the quotient $H \coloneqq G/\gamma_n^{\mathcal P}(G)$ is not trivial. 
Hypothesis \ref{HR1} says that $H_* \coloneqq \gamma_*^{\mathcal P}(G)/\gamma_n^{\mathcal P}(G)$ 
is a $\mathcal P$-filtration on $H$, so that $\gamma_n^{\mathcal P}(H) \subseteq H_n = \{1\}$, 
proving our claim, and the Lemma.
\end{proof}

\begin{remark}
Applied to inclusions, hypothesis~\ref{HR0} implies that $\res(\mathcal R_{\mathcal P})$ 
is stable under taking subgroups. In fact, it implies directly that 
$\mathcal R_{\mathcal P}$ itself is closed under taking subgroups: 
if $H \subset G$, then $\gamma_*^{\mathcal P}(G) \cap H$ is a 
$\mathcal P$-filtration, so it contains $\gamma_*^{\mathcal P}(H)$. 
\end{remark}

\begin{example}
\label{ex:res-finite}
Although the above terminology is not introduced to deal with residual finiteness, 
let us point out that residually finite groups do fit in this framework. 
Indeed, if $G$ is a group, we can introduce, following~\cite{DFPR}, the 
intersection $N_n$ of all the subgroups of index at most~$n$ of $G$. 
Then obviously $N_1 = G$ and $N_{n+1} \subseteq N_n$, so we have a 
filtration on $G$, and we can show that it is $\gamma_*^{\mathcal P}(G)$ 
for some property $\mathcal P$ satisfying our hypotheses. Namely, say that 
$G_*$ satisfies $\mathcal P$ if for every $g \notin G_n$, there is a subgroup 
$H$ of index at most $n$ containing $G_n$ but not $g$. It is easy to check that 
this works; moreover, since $\bigcap_{i\ge 1} N_i$ is the intersection of all finite-index 
subgroups of $G$, Lemma~\ref{res_prop_and_separation} says that $\res(\mathcal R_{\mathcal P})$ 
is the class of residually finite groups. However, notice that the class 
$\mathcal R_{\mathcal P}$ is larger than the class of finite groups, 
unless we restrict to finitely generated groups. Indeed, if $G$ is 
finitely generated, then all subgroups $N_n$ have finite index 
in $G$ (see~\cite{DFPR}), so one of them is trivial if and only 
if $G$ is finite; but this does not hold in general.
\end{example}

In what follows we focus on central filtrations. When $\mathcal P$ implies 
centrality, then $\mathcal R_{\mathcal P}$ is contained in the class of 
nilpotent groups, so $\res(\mathcal R_{\mathcal P})$ is a subclass of 
residually nilpotent groups.

\begin{example}[Residually nilpotent groups] 
\label{ex:res-nilp}
If $\mathcal P$ is ``being central'', then $\gamma_*^{\mathcal P}$ is the 
lower central series, and $\mathcal R_{\mathcal P}$ is (by definition) 
the class of nilpotent groups.
\end{example}

\begin{example}
\label{ex:bs-mn}
For each $0<m\le \abs{n}$, let $G=\BS(m,n)$ be the \emph{Baumslag--Solitar group}, 
defined by the  presentation $\langle t,a\mid ta^mt^{-1}=a^n\rangle$. 
Then $G$ is residually finite if and only if $m=1$ or $m=\abs{n}$, 
and it is residually nilpotent if and only if 
$m=1$ and $n\ne 2$, or $m=\abs{n}=p^r$ for some prime $p$ 
and some $r>0$; see \cite{Bardakov-Neshchadim-BS,KMP} 
and references therein (see also Example~\ref{ex:baumslag-solitar}, 
where the lower central series of $\BS(1,n)$ is computed). In particular, 
$\BS(1,2)$ is residually finite but not residually nilpotent, whereas 
the Klein bottle group $\BS(1,-1)$ is residually nilpotent.
\end{example}

\begin{example}[RTFN groups]
\label{ex:rtfn-groups}
When $\mathcal P$ is ``being torsion-free central'', then $\gamma_*^{\mathcal P}$ is 
$\ssqrt{\gamma_*}$. For a group $G$, having $\ssqrt{\gamma_n} = \{1\}$ is equivalent 
to having $\gamma_n = \{1\}$ and $G = G/\gamma_n$ being torsion-free, so 
$\mathcal R_{\mathcal P}$ is the class of torsion-free nilpotent groups, 
and $\res(\mathcal R_{\mathcal P})$ is the class of \emph{residually torsion 
free nilpotent}\/ (RTFN) groups.
\end{example}

\begin{example}[Residually $\pi$-free nilpotent groups]
\label{ex:res-pi-free-nilp}
As a straightforward generalization of the previous example, if $\pi$ is 
a set of primes, and $\mathcal P$ is ``being $\pi$-free central'', then 
$\mathcal R_{\mathcal P}$ is the class of $\pi$-free nilpotent groups.
\end{example}

\begin{example}[Residually $p$ groups]
\label{ex:res-p-groups}
Let $q \geq 2$ be an integer, and let $\mathcal P$ be ``being 
$q$-torsion central''. Then, as we shall see (Proposition~\ref{res(R_q)}), 
$\mathcal R_q \coloneqq \mathcal R_{\mathcal P}$ is the class of nilpotent 
groups whose exponent divides a power of $q$. If $q = p$ is prime, 
and \emph{for finitely generated groups}, being residually $\mathcal R_p$ 
is the same as being residually a finite $p$-group, or 
\emph{residually $p$}\/ (Corollary~\ref{res(R_q)_fg}).
\end{example}

\subsection{Residual properties and semi-direct products}
\label{subsec:residual_prop_sd_prod}

We now turn to the heart of our study of residual properties of semi-direct products. 
Precisely, we are looking for conditions to decide when a semi-direct product is 
residually $\mathcal R_{\mathcal P}$, for some property $\mathcal P$ of filtrations. 
In the sequel, we will deduce conditions for a semi-direct product to be residually 
torsion-free nilpotent (RTFN), or residually a finite $p$-group.

Combining Lemma~\ref{res_prop_and_separation} with 
Theorem~\ref{dec_of_P-LCS}, we obtain the following generalization 
of a result of Falk and Randell~\cite{FR88} to variants of the 
lower central series:

\begin{theorem}
\label{thm:res-P}
Under hypotheses \ref{H0}--\ref{H4} and \ref{HR1}, 
if a group $K$ acts on a group $H$ such that 
the induced action on $\gr_*^{\mathcal P}(H)$ is trivial, 
then $H \rtimes K$ is residually $\mathcal R_{\mathcal P}$ 
if and only if $H$ and $K$ are.
\end{theorem}

\begin{proof}
By Lemma~\ref{res_prop_and_separation}, $H \rtimes K$ is residually 
$\mathcal R_{\mathcal P}$ if and only if $\gamma_*^{\mathcal P}(H \rtimes K)$ 
is separating. By Theorem~\ref{dec_of_P-LCS}, this is equivalent to both 
$\gamma_*^K(H)^{\mathcal P}$ and $\gamma_*^{\mathcal P}(K)$ being separating. 
The latter condition says precisely that $K$ is residually $\mathcal R_{\mathcal P}$. 
As for the former, the hypothesis that $K$ acts trivially on $\gr_*^{\mathcal P}(H)$ 
implies, by the third statement of Theorem~\ref{dec_of_P-LCS}, that 
$\gamma_*^K(H)^{\mathcal P} = \gamma_*^{\mathcal P}(H)$. Hence the separability 
of $\gamma_*^K(H)^{\mathcal P}$ is equivalent to $H$ being residually 
$\mathcal R_{\mathcal P}$, and the conclusion follows.
\end{proof}

\begin{remark}
The condition on the action can be relaxed: in order to get the conclusion, 
we do not need $\gamma_*^K(H)^{\mathcal P}$ and $\gamma_*^{\mathcal P}(H)$ 
to be equal, but only to be \emph{cofinal} (see~\cref{subsec:residual_prop_sd_prod_unipotent}).
\end{remark}

\begin{remark}
In \cite[Thm.~5.7]{Metaftsis-Papistas}, Metaftsis and Papistas show 
that certain groups of the form $G=H\rtimes K$ are residually nilpotent, provided $H$ is 
residually $p$, $K$ is residually nilpotent, and $K$ acts trivially on $H^{\ab} \otimes \F_p$ 
(a direct consequence of Theorem~\ref{BBN_thm}, which is also 
\cite[Thm.~4.1]{Bardakov-Bryukhanov-Neshchadim}). This is not always the case 
if this condition does not hold. For instance, take 
$G = \Z/p \rtimes \Z =\langle t, s \mid s^p=tst^{-1}s^{-2}=1 \rangle$, 
for some prime $p \ne 2$. Then $s = [t,s] = [t, [t,s]] = [t, [t,[t,s]]] = \cdots$, 
and so  $s\in \gamma_\omega(G)$, showing that $G$ is not residually nilpotent.
\end{remark}

\subsection{RTFN groups}
\label{subsec:rtfn}

A group $G$ is said to be {\em residually torsion-free nilpotent}\/ (RTFN) 
if every non-trivial element can be detected in a torsion-free nilpotent 
quotient. Clearly, such groups are residually nilpotent, but the 
converse does not hold. For instance, every finite nilpotent group 
is residually nilpotent, but not RTFN.

By \cite[Ch.~VI, Theorem~2.26]{Passi}, the group $G$ is RTFN precisely 
when the group-algebra $\Q[G]$ is residually nilpotent, that is, 
$\bigcap_{n\geq 1}I^{n}=\{0\}$, where $I$ is the augmentation ideal.  
When $G$ is finitely generated, the RTFN condition is equivalent to 
the injectivity of the canonical map $G\to \mathfrak{M}(G)$ to 
its prounipotent (or, Mal'cev) completion, where $\mathfrak{M}(G)$ 
is the set of group-like elements in the Hopf algebra obtained by 
completing $\Q[G]$ with respect to the $I$-adic filtration (see for 
instance \cite{Suciu-Wang} and references therein).

Finitely generated RTFN groups are torsion-free and bi-orderable. 
If $G$ is residually nilpotent and $\gr_n (G)$ is torsion-free for 
$n\ge 1$, then $G$ is residually torsion-free nilpotent. 
Examples of RTFN groups include free groups \cite{Magnus35}, 
fundamental groups of surfaces having a handle \cite{Baumslag2010}, 
right-angled Artin groups \cite{Duchamp-Krob}, and the pure braid 
groups \cite{FR88}.

\medskip

As noticed above in Example~\ref{ex:rtfn-groups}, RTFN groups correspond to 
the case when $\mathcal P$ is ``being torsion-free central'' in the previous 
machinery. Then our hypotheses are satisfied (see Lemma~\ref{checking_(HRi)}), 
we have $\gamma_*^{\mathcal P} = \ssqrt{\gamma_*}$ (see 
Proposition~\ref{min_torsion-free_filtrations}), and 
Lemma~\ref{res_prop_and_separation} recovers the well-known result 
(compare~\cite[Proposition 7.2(e)]{Bass-Lubotzky}):

\begin{corollary}
\label{cor:rtfn}
Let $G$ be a group. The following properties are equivalent:
\begin{itemize}[itemsep = 1.5pt]
\item The group $G$ is RTFN;
\item The filtration $\gamma_*^{\rat}(G)$ is separating, that is, 
$\bigcap_{n \geq 1} \gamma_n^{\rat}(G) = \{1\}$;
\item There exists a separating torsion-free central filtration $G_*$ such that $G = G_1$.
\end{itemize}
\end{corollary}

Using Corollary~\ref{cor_engdeg1_rat} to weaken its hypothesis, 
Theorem~\ref{thm:res-P} becomes, in this context:
\begin{corollary}
\label{cor:rtfn-ext}
If a group $K$ acts on a group $H$ such that 
the induced action on $H^{\ab} \otimes \Q$ (or, equivalently, on $H^{\abf}$) is trivial, 
then $H \rtimes K$ is RTFN if and only if $H$ and $K$ are.
\end{corollary}

This has the following topological consequence:
\begin{corollary}
\label{cor:RTFN_pi1_in_fibration} 
Let $p\colon E \surj B$ be a Serre fibration between connected, 
locally path-connected spaces, with fiber $F$. Suppose that the 
corresponding monodromy action of $\pi_1(B)$ on $H_1(F, \Q)$ is 
trivial, and that the fibration splits up to pointed homotopy 
(that is, there exists $s \colon B \to E$ preserving basepoints 
such that $p \circ s$ is homotopic to $\id_B$ via a pointed homotopy). 
Then $\pi_1(E)$ is RTFN if and only if $\pi_1(F)$ and $\pi_1(B)$ are.
\end{corollary}

\begin{proof}
The existence of a splitting $s$ implies that the long exact sequence 
of homotopy groups of the fibration breaks into split short exact 
sequences. At the level of fundamental groups, this yields   
a semi-direct product decomposition $\pi_1(E) \cong \pi_1(F) 
\rtimes \pi_1(B)$. The corresponding action of $\pi_1(B)$ on $\pi_1(F)$, 
which is the restriction of the conjugation action of $\pi_1(E)$ on itself, 
induces the monodromy action on $\pi_1(F)^{\ab} = H_1(F)$. The other 
hypothesis says that the induced action on $\pi_1(F)^{\ab} \otimes \Q 
\cong H_1(F,\Q)$ is trivial, and so we can apply Corollary~\ref{cor:rtfn-ext}.
\end{proof}

Given a space $X$ and a map $f\colon X\to X$, we let 
$T_f=X\times [0,1]/(x,0)\sim (f(x),1)$ be the mapping torus of $f$. 

\begin{corollary}
\label{cor:circle-fibration} 
Let $X$ be a connected, locally path-connected space, and  
let $f\colon X\to X$ be a map inducing the identity on $H_1(X;\Q)$. 
Then $\pi_1(T_f)$ is RTFN if and only if $\pi_1(X)$ is.
\end{corollary}

\begin{proof}
Apply Corollary~\ref{cor:RTFN_pi1_in_fibration} to the split fibration 
$T_f \to S^1$ (with fiber $X$), which is a locally trivial fibration, 
and in particular a Serre fibration.
\end{proof}

\begin{example}
\label{ex:torelli}
Let $\Sigma$ be a closed, orientable surface, 
let $f\colon \Sigma\to\Sigma$ be a smooth map that 
belongs to the rational Torelli group (i.e., $f$ induces the 
identity on $H_1(\Sigma;\Q)$), and let $M=T_f$ be the 
$3$-manifold which fibers over the circle with fiber $\Sigma$ 
and monodromy $f$. As mentioned previously, the 
surface group $\pi_1(\Sigma)$ is known to be RTFN. 
By the above corollary, then, the $3$-manifold group $\pi_1(M)$ 
is also RTFN.
\end{example}

\begin{example}[Braids on surfaces]
\label{ex:braids}
Using the Fadell-Neuwirth fibrations to decompose them as semi-direct products, we will 
show in \cref{sec:braids} that pure braid groups on orientable surfaces different from the 
sphere are RTFN. This will actually use a bit more that Corollary~\ref{cor:rtfn-ext}, which 
can be applied only to planar surfaces (in particular for the usual pure braid group -- see 
Example~\ref{ex:pure braids}): the more subtle Theorem~\ref{Residual_prop_of_sdp} 
will be needed to handle the general case.
\end{example}

\begin{remark}
\label{rk:res_pi-free_nilp}
All the results of this section (\cref{subsec:rtfn}) generalize easily to residually 
$\pi$-free nilpotent groups, where $\pi$ is a set of primes (that is, to the 
context described in~\cref{sec:pi-free}; see also Example~\ref{ex:res-pi-free-nilp}). 
Then RTFN groups correspond to the case where $\pi$ is the set of all primes, 
and residually nilpotent groups to the case where $\pi = \varnothing$. 
\end{remark}

\subsection{RTN$_q$ groups} 
\label{par_RTNq}

We now investigate the residual property associated with $q$-torsion central 
filtration. Namely, we apply the above machinery to the case where our property 
$\mathcal P$ of filtrations is ``being central of $q$-torsion", for some integer $q$. 
We first show that the corresponding class $\mathcal R_{\mathcal P}$ is the 
class $\mathcal R_q$ of \emph{nilpotent $q$-groups} (see the next definition). 
Then $\res(\mathcal R_q)$ is the class of \emph{residually nilpotent $q$-groups}, 
which we call RTN$_q$ groups, for short. These are groups whose non-trivial 
elements are detected in nilpotent quotients which are $q$-groups, in the following sense: 

\begin{definition}
\label{def:q-group}
A \emph{$q$-group}\/ is a group $G$ having some power of $q$ as exponent, 
that is: there exists an integer $t \geq 1$ such that every element $g \in G$ 
satisfies $g^{q^t} = 1$.
\end{definition}

If $q = p$ is prime and $G$ is finite, this definition coincides with 
the usual definition of a (finite) $p$-group. We will in fact show 
(Proposition~\ref{prop:gf_q-grps_are_finite}) that finitely 
generated $q$-groups are finite.

\begin{proposition}
\label{res(R_q)}
Let $G$ be a group. Then there exists $n \geq 1$ such that 
$\gamma_{n+1}^{\qt}(G) = 1$ if and only if $G$ is a nilpotent $q$-group.
\end{proposition}

\begin{proof}
If $g \in G$, we clearly have $g^{q^{k-1}} \in \gamma_k^{\qt}(G)$, 
for every $k \geq 1$. Thus, $\gamma_{n+1}^{\qt}(G) = 1$ implies $G^{q^n} = 1$. 
Moreover, $\gamma_*^{\qt}$ contains $\gamma_*$, so it also implies that 
$G$ is $n$-nilpotent. Let us prove the converse: let us assume that 
$\gamma_{s+1}(G) = 1$ and that $G^{q^t} = 1$ for some $s,t \geq 1$, and 
let us prove that $\gamma_{st + 1}^{\qt}(G) = 1$. We need to find a 
$q$-torsion filtration $N_*$ on $G$ with $N_{st +1} = 1$. For every 
integer $\alpha \geq 0$ and every $r \in \{0, \dots , q - 1\}$, we 
define $N_{1 + \alpha q + r}$ to be 
$\gamma_{\alpha+1}^{q^r} \gamma_{\alpha +2}$. Since $G^{q^t} = 1$, 
we have $N_{1 + \alpha t} = \gamma_{\alpha + 1}$. Thus, $N_*$ is a 
filtration on $G$, and it is a refinement of its lower central series, 
so it is central. By construction, $N_*$ is of $q$-torsion, and 
$N_{1+st} = \gamma_{s+1}(G) = 1$, whence our claim.
\end{proof}

In this context, our hypotheses are satisfied (Lemma~\ref{checking_(HRi)}), 
and Lemma~\ref{res_prop_and_separation} becomes:
\begin{corollary}
\label{cor:rtnq}
Let $G$ be a group. The following properties are equivalent:
\begin{itemize}[itemsep = 1.5pt]
\item The group $G$ is residually $\mathcal R_{\mathcal P}$;
\item The filtration $\gamma_*^{\qt}(G)$ is separating, that is, $\bigcap_{n \geq 1} \gamma_n^{\qt}(G) = \{1\}$;
\item There exists a separating $q$-torsion central filtration $G_*$ such that $G = G_1$.
\end{itemize}
\end{corollary}

Using Corollary~\ref{triv_action_q} to weaken its hypothesis, Theorem~\ref{thm:res-P} 
becomes, in this context:
\begin{corollary}
\label{cor:rtnq-ext}
If a group $K$ acts on a group $H$ such that 
the induced action on $H^{\ab} \otimes (\Z/q)$ is trivial, 
then $H \rtimes K$ is RTN$_q$ 
if and only if $H$ and $K$ are.
\end{corollary}

We will see later on (in~\cref{subsec:residually-p}) that the previous two 
results are in fact generalizations of results of Paris~\cite[Prop. 2.3]{Paris} 
and Bellingeri and Gervais~\cite[Cor.~3.5]{Bellingeri-Gervais} to any integer 
(not necessarily prime) and any group (not necessarily finitely generated). 
First, we need to study nilpotent $q$-group and investigate the link with 
(finite) $p$-groups.

\subsection{Nilpotent $q$-groups} 

We now investigate in more detail the structure of nilpotent $q$-groups. We first 
show how such groups decompose as a product of (possibly infinite) $p$-groups, 
allowing us to reformulate the definition of RTN$_q$ groups in terms of primes 
dividing $q$. Then we show that that finitely generated $q$-groups are in fact 
finite (Proposition~\ref{prop:gf_q-grps_are_finite}), showing in particular 
that is $q = p$ is prime, finitely generated RTN$_q$ group are just residually 
(finite) $p$-groups.

For a group $G$, let $\Tors_p(G)$ denote the set of \emph{$p$-elements}\/ 
of $G$ (the ones whose order is a power of $p$). We recall the following 
Sylow-type theorem for nilpotent groups (see \cite[\S 5.2]{Robinson}):

\begin{proposition}
\label{torsion_as_prod_of_p-torsion}
Let $G$ be a nilpotent group whose elements have finite order. Then $G$ is the 
restricted direct product of the subgroups $\Tors_p(G)$, for all prime numbers $p$.
\end{proposition} 

\begin{proof}
That the sets $\Tors_p(G)$ are subgroups of $G$ is proved in the same way as 
Lemma~\ref{torsion_subgroup}. That elements of coprime order commute follows 
from the usual formula:
\begin{equation}
[x, y^q] = [x, y] \cdot {}^y [x, y]\ \cdots\ {}^{y^{q-1}}[x, y].
\end{equation}
Indeed, let $x$ (resp.~$y$) be elements of order $\alpha$ (resp.~$\beta$), 
with $\alpha$ and $\beta$ coprime. We can then show by induction that 
$[x,y] \in \gamma_k$ for every $k \geq 1$. This clearly holds for 
$k = 1$, and if  $[x,y] \in \gamma_k$ for some $k$, then modulo 
$\gamma_{k+1}$, $[x,y]$ is central, so the above formula gives 
$[x, y^q] \equiv [x, y]^q$ for all integer $q$. Thus $[x,y]^\alpha \equiv 
[x,y]^\beta \equiv 1$, so $[x,y] \equiv 1$, which means that $[x,y] \in 
\gamma_{k+1}$, and the induction is complete. Since $G$ is nilpotent, 
we must have $[x,y] = 1$.

Finally, let us show that any $g \in G$ can be written uniquely as a finite 
product of elements $g_p \in \Tors_p(G)$ indexed by prime numbers $p$. 
By hypothesis, $g$ is of finite order $n = \prod p_i^{\alpha_i}$, 
so this holds in $\langle g \rangle \cong \Z/n \cong \prod \Z/p_i^{\alpha_i}$, 
by the Chinese remainder theorem. But if $g = xy$ is a decomposition into 
a product of elements $x$ and $y$ of coprime order $\alpha$ and $\beta$, 
then for $u, v$ satisfying $u \alpha + v \beta = 1$, we have $x = g^{v\beta}$ 
and $y = g^{u \alpha}$. As a consequence, for any decomposition of $g$ into 
a product of elements $g_p \in \Tors_p(G)$, the $g_p$ must be powers of $g$: 
the unique such decomposition in $\langle g \rangle$ is also the unique one in $G$.
\end{proof}

We readily infer (using projection on factors) that a group $G$ is RTN$_q$ if 
and only is every non-trivial elements $g \in G$ is detected in a quotient 
which is a $p$-group for some prime factor $p$ of $q$ (that may depend on $g$). 
In other words:
\begin{corollary}
\label{res(R_q)_in_terms_of_divisors}
Let $p_1, \dots, p_l$ be the list of prime divisors of $q$. Then: 
\[
\res(\mathcal R_q) = \res\left(\mathcal R_{p_1} 
\cup \cdots \cup \mathcal R_{p_l}\right).
\]
\end{corollary}

In order to make the link with residually $p$-groups, we now turn 
to finitely generated groups.
\begin{proposition}
\label{prop:gf_q-grps_are_finite}
A finitely generated group $G$ is a nilpotent $q$-group if and only if 
it is a product of finite $p$-groups, indexed by prime divisors $p$ of $q$.
\end{proposition}

\begin{proof}
The ``only if'' part follows easily from the fact that finite $p$-groups 
are nilpotent. Let us prove the converse. Let $G$ be a finitely generated 
nilpotent group satisfying $G^{q^t} = 1$ for some $t$. Then we can apply 
Proposition~\ref{torsion_as_prod_of_p-torsion} to $G$. Since the orders 
of elements of $G$ divide $q^t$, the subgroup $\Tors_p(G)$ is trivial 
unless $p$ divides $q$. Moreover, the factors in the product must be finitely 
generated for their product to be so. We are left with showing that 
\emph{every finitely generated nilpotent $H$ having only $p$-elements is finite}. 
This follows rather easily from Proposition~\ref{engdeg1}. Indeed $H^{\ab}$ is a 
finitely generated abelian group consisting of $p$-elements. Let $p^m$ be the 
maximal order of its generators. Then  $p^m x = 0$ for every $x \in H^{\ab}$. 
But then all the groups $\gr_i(H)$ must be finitely generated abelian groups 
whose elements $x$ satisfy $p^m x = 0$, that is, finitely generated 
$\Z/p^m$-modules; hence, they are finite. By induction using the 
decomposition of $H/\gamma_{i+1}$ as an extension of 
$H/\gamma_i$ by $\gr_i(H)$, the quotients $H/\gamma_i$ are finite. Since $H$ 
is nilpotent, $H/\gamma_i = H$ for $i$ big enough; therefore, $H$ is finite.
\end{proof}

\begin{corollary}
\label{res(R_q)_fg}
Let $p_1, \dots, p_l$ be the list of prime divisors of $q$. Then: 
\[
\res(\mathcal R_q) \cap \{\text{f.g. groups}\} = \res\left(\{p_1\text{-groups}\} 
\cup \cdots \cup \{p_l\text{-groups}\}\right) \cap \{\text{f.g. groups}\}.
\]
\end{corollary}

In particular, if $q = p$ is prime, then a finitely generated group is RTN$_p$ 
if and only if it is residually a (finite) $p$-group (see the more precise Corollary~\ref{res_p_vs_RTN_p}).

\subsection{Residually $p$ groups}
\label{subsec:residually-p}

If $p$ is prime, a group $G$ is said to be {\em residually $p$}\/ if every 
nontrivial element of $G$ can be detected in a \emph{finite} $p$-group quotient. 
Since $p$-groups are nilpotent, residually $p$ groups are residually nilpotent, 
but the converse does not hold. 

\begin{example}
\label{ex:bs-p}
The Baumslag--Solitar group $\BS(m,n)$ from Example \ref{ex:bs-mn} 
is residually $p$ if and only if $m=1$ and $n\equiv 1 \!\pmod{p}$; or 
$m=n=p^r>1$; or $m=-n=2^r>1$ and $p=2$; see \cite{Bardakov-Neshchadim-BS,KMP}.  
It follows that the group $\BS(1,p+1)$ is residually $p$, while the group 
$\BS(1,p)$ is residually nilpotent (if $p\ne 2$) but not residually $p$.
\end{example}

Let us spell out the direct consequence of Corollaries~\ref{res(R_q)} 
and~\ref{res(R_q)_fg} that allows us to apply the results of~\cref{par_RTNq} 
to residually~$p$ groups:
\begin{corollary}
\label{res_p_vs_RTN_p}
Residually~$p$ groups are RTN$_p$. Conversely, finitely generated RTN$_p$ 
groups are residually~$p$.
\end{corollary}

The converse implication does not hold for infinitely generated groups. Indeed, 
there are nilpotent $p$-groups (in our sense) that are not residually finite, 
hence not residually~$p$.
\begin{example}[see~{\cite[p.~299]{Wehrfritz}}]
\label{ex:nilp_p-gp_not_res_fini}
Let $H$ be the group generated by $x_i, y_i$ (for $i \in \N$) and $z$, submitted to the relations 
$[x_i, y_i] = z$, $[x_i,z] = [y_i,z] = 1$ and $x_i^p = y_i^p = z^p = 1$ (for all $i \in \N$), 
and  $[x_i,x_j] = [y_i,y_j] = [x_i,y_j] = 1$ for $i \neq j$. This group can be seen as the 
restricted central product of countably many copies of the Heisenberg group $H_p$ 
over $\F_p$, which is the group of upper-triangular unitary matrices with coefficients 
in $\F_p$. One can see that $H$ is a central extension of $(\Z/p)^{(\N)}$ by $\Z/p$, 
so it is $2$-nilpotent, and its exponent divides $p^2$. In fact, it has exponent $p$, 
since $H_p$ does. As a consequence, it is a nilpotent $p$-group, in our sense. 
However, in any finite quotient of $H$, some relation of the form $x_i = x_j$ for 
$i \neq j$ must become true, which implies $z = [x_i, y_i] = [x_j, y_i] = 1$, so 
$z$ is killed in any such quotient. Thus, $H$ is not residually finite, and in 
particular not residually~$p$.
\end{example}

Thanks to Corollary~\ref{res_p_vs_RTN_p}, Corollary~\ref{cor:rtnq} specializes to Paris' 
\cite[Prop.~2.3]{Paris}:

\begin{lemma}[\cite{Paris}]
\label{lem:paris-bis}
If a group $G$ is residually $p$, then $\bigcap_{n\ge 1} \gamma^{\qt[p]}_{n}(G)=
\{1\}$. The converse holds if $G$ is finitely generated.
\end{lemma}

As for Corollary~\ref{cor:rtnq-ext}, it specializes to the following result
of Bellingeri and Gervais:

\begin{corollary}[{\cite[Cor.~3.5]{Bellingeri-Gervais}}]
\label{cor:rtnp-ext}
Let $K$ and $H$ be finitely generated groups. If $K$ acts on a group $H$ 
such that the induced action on $H^{\ab} \otimes \F_p$ is trivial, 
then $H \rtimes K$ is residually $p$ if and only if $H$ and $K$ are.
\end{corollary} 

\begin{remark}
\label{rem:BG-discuss}
The statement of \cite[Cor.~3.5]{Bellingeri-Gervais} assumes that $H \rtimes K$ 
is a $p$-almost direct product, which is equivalent to our hypothesis, since 
$H^{\ab} \otimes \F_p \cong H_1(H, \mathbb F_p)$. Those authors neglected to spell 
out the finite generation hypothesis in their statement, although this 
hypothesis is essential in their use of Paris' result. Finally, notice 
that this statement being a corollary of their Theorem~3.3 is a 
particular case of our Theorem~\ref{thm:res-P} being a rather 
direct consequence of Theorem~\ref{dec_of_P-LCS} (which, 
in this context, specializes as Theorem~\ref{dec_of_q-LCS}).
\end{remark}

\subsection{Finitely generated RTFN groups are residually $p$-groups}
\label{subsec:fg_RTFN_implies_res_p}
This statement is a classical theorem by Gruenberg~\cite[Thm.~2.1]{Gruenberg}, 
who gave a very elegant proof. We give a different proof with our language, 
which allows for a greater control on the filtrations involved: albeit 
more technical, ours is more constructive, in that it does not rely 
on any statement of existence of a maximal subgroup satisfying certain 
properties. Our main tool is a generalization of~Proposition~\ref{Dark_for_products}, 
which is another consequence of Dark's theorem:
\begin{proposition}
\label{Dark_for_n_products} 
Let $n \geq 1$ be an integer. There exists a unique map 
$\theta_n \colon \N \rightarrow F_n = \langle x_1, \dots, x_n \rangle$ 
such that
\[
x_1^\alpha \cdots x_n^\alpha = (x_1 \cdots x_n)^\alpha\prod\limits_{r \geq 2} 
\theta_n(r)^{\binom{\alpha}{r}}
\]
for all $\alpha \in \N$, and such that each $\theta_n(r)$ is a product of 
$\{x_1^{\pm 1}, \dots, x_n^{\pm 1}\}$-commutators of length at least $r$.
\end{proposition}

\begin{remark}
Here we have included in our statement the analogue of Remark~\ref{theta(1)}: 
the map $\theta_n$ given by Dark's theorem is independent of $\alpha$, and 
we can use that fact to compute it. In particular, $\theta_n(0) = 1$ and 
$\theta_n(1) = x_1 \cdots x_n$.
\end{remark}

Our main technical statement is the following:
\begin{lemma}
\label{lem_interaction_G*_and_G*^(q)}
Let $q \geq 2$ and let $G_*$ be a central filtration 
without $q$-torsion. Then:
\[
G_j^{\qt} \cap G_i \subseteq  G_i^{q^{j-i}}G_{i + 1},
\]
for all $j > i \geq 1$.
\end{lemma}

\begin{proof}
Let $K_* \coloneqq G_*^{\qt}$, and let $j > i \geq 1$. We will show by induction 
on $k \in \{1, \dots, i\}$ that $K_j \cap G_i \subseteq H_k G_{i+1}$,
where $H_k$ is the product of the $G_s^{q^t}$ for $s+t = j$ and $k \leq s \leq i$. 
For $k = 1$, this is a direct consequence of the explicit description of 
$K_j$ from Corollary~\ref{explicit_min_q}: $K_j = \prod_{s+t \geq j}G_s^{q^t}$. 
The case $k = i$ will give the announced result. Let us then assume that 
the inclusion $K_j \cap G_i \subseteq H_k G_{i+1}$ has been proven 
for some $k \in \{1, \dots, i-1\}$, and let us prove that 
$K_j \cap G_i \subseteq H_{k+1} G_{i+1}$. 

Let $g \in K_j \cap G_i$. The induction hypothesis says that $g \in H_k G_{i+1}$. 
Since $H_k = G_k^{q^{j-k}}H_{k+1}$, if we pose $\alpha \coloneqq q^{j-k}$, 
there are elements $x_1, \dots, x_n \in G_k$ such that 
$g \in x_1^\alpha \cdots x_n^\alpha H_{k+1}G_{i+1}$. 
Now we can apply the formula from Proposition~\ref{Dark_for_n_products}. 
We are going to show that all the factors involved in the resulting 
expression of $x_1^\alpha \cdots x_n^\alpha$ are in $H_{k+1}$ or $G_{i+1}$. 
This will imply that $g \in H_{k+1}G_{i+1}$, as expected.

We first inspect the formula modulo $G_{k+1}$. Recall that $k < i$ and 
that $g \in G_i$, so $g \in G_{k+1}$. Moreover, $G_{i+1} \subseteq G_{k+1}$. 
We also have that $H_{k+1} \subseteq G_{k+1}$, by definition. As a consequence, 
$x_1^\alpha \cdots x_n^\alpha \in g G_{i+1}H_{k+1} \subseteq G_{k+1}$. 
Furthermore, if $r \geq 2$, then $\theta_n(r)$ is a product of 
commutators of length at least $2$ of the $x_i^{\pm 1} \in G_k$, 
so by centrality of $G_*$, they belong to $G_{k+1}$. Finally, we get 
$(x_1 \cdots x_n)^\alpha \in G_{k+1}$. But $x_1 \cdots x_n \in G_k$, 
$\alpha = q^{j-k}$, and our hypothesis on $G_*$ says that $G_k/G_{k+1}$ 
does not have non-trivial element $u$ satisfying $u^q = 1$. 
We conclude that $x_1 \cdots x_n \in G_{k+1}$. As a consequence, 
$(x_1 \cdots x_n)^\alpha \in 
G_{k+1}^\alpha \subseteq G_{k+1}^{q^{j-k-1}}$, so 
$ (x_1 \cdots x_n)^\alpha \in H_{k+1}$.

If $r \geq i+1$, then $\theta_n(r)$ is a product of commutators of length at least $i+1$ 
of elements of $G_1$ (in fact of $G_r$), so by centrality of $G_*$, we have 
$\theta_n(r) \in G_{i+1}$ (we do not need a sharper bound here).

If $2 \leq r \leq i$, then $\theta_n(r)$ is a product of commutators 
of length at least $r$ of elements of $G_k$, so by centrality of $G_*$, 
we have $\theta_n(r) \in G_{k+r-1}$. In order 
to show that $\theta_n(r)^{\binom{\alpha}{r}} \in H_{k+1}$, we need to investigate 
$v_q\left(\binom{\alpha}{r}\right)$. We apply Lemma~\ref{inégalité_binôme}, recording 
that $2 \leq r \leq i$ and $\alpha = q^{k-i}$:
\[
v_q\bigl(\tbinom{\alpha}{r}\bigr) \geq (j-k) - v_q(r) \geq (j-k) - (r-1).
\]
Finally, $\theta_n(r)^{\binom{\alpha}{r}} \in 
G_{k+r-1}^{q^{j-k-r+1}} \subseteq H_{k+1}$, 
since $k+r-1 \geq k+1$.As announced, we have shown that all 
the factors of the formula for $x_1^\alpha \cdots x_n^\alpha$ 
are in $H_{k+1}G_{i+1}$, hence $g \in H_{k+1}G_{i+1}$, 
which finishes the induction step and the proof.
\end{proof}

\begin{corollary}[{\cite[Thm.~2.1]{Gruenberg}}]
\label{Gruenberg_thm}
Finitely generated RTFN groups are residually (finite) $p$-groups, for every prime $p$.
\end{corollary}

\begin{proof}
Let $p \geq 2$ be a prime number. If $G$ is finitely generated RTFN, then 
$G_* \coloneqq \gamma^\Q_*(G)$ is a torsion-free 
filtration of $G$, which is separating (Cor.~\ref{cor:rtfn}). Moreover, the groups 
$\gr_i(G_*)$ are finitely generated torsion-free abelian (Cor.~\ref{rat_gr_i_are_fg}), 
so for every $i \geq 1$, we have $\bigcap_{j \geq 1} p^j \gr_i(G_*) = \{0\}$. 
Let $g \in G$, and let $i \geq 1$ such that $g \in G_i \setminus G_{i+1}$, so that 
$\bar g \in \gr_i(G_*)$ is not trivial. If $g \in G_j^{\qt[p]}$ with $j > i$, 
then by Lemma~\ref{lem_interaction_G*_and_G*^(q)}, we have $\bar g \in p^{j-i} \gr_i(G_*)$. 
This cannot hold for every $j > i$, so there is some $j > i$ such that $g \notin G_j^{\qt[p]}$. 
This shows that $G_*^{\qt[p]}$ is a separating $p$-torsion filtration, whence the 
result by Corollaries~\ref{cor:rtnq} and~\ref{res_p_vs_RTN_p}. 
\end{proof}

\begin{remark}
Because of Corollary~\ref{res(R_q)_fg} and~\ref{res_p_vs_RTN_p}, 
for finitely generated groups, being a residually (finite) $p$-group 
for every prime $p$ is the same as being RTN$_q$ for every $q \geq 2$.
\end{remark}

\begin{remark}
The finite generation hypothesis is essential here, as the case of the 
additive group $\Q$ shows: this group is RTFN, but $p\Q = \Q$, so it 
does not have any non-trivial finite $p$-group as a quotient. 
Moreover, the hypothesis that $G$ is RTFN is used in full in the 
proof, even for showing that $G$ is residually $p$ for any particular 
choice of $p$: our argument fails if there is torsion in the filtration 
considered. Indeed, Lemma~\ref{lem_interaction_G*_and_G*^(q)} can be 
applied only for a filtration with no $p$-torsion, but if there is 
torsion of order prime to $p$ in a filtration $G_*$, then for some 
$i \geq 1$, we have $\bigcap_{j \geq 1} p^j \gr_i(G_*) \neq \{0\}$, 
so the argument of the previous proof fails.
\end{remark}

\begin{remark}
\label{rem:p-not-rt}
Examples of finitely generated 
groups that are residually $p$ for all primes $p$ yet are not RTFN 
were given by Hartley \cite{Hartley}.
\end{remark}

\begin{example}
\label{ex:BS1n}
The Baumslag--Solitar groups $\BS(1,n)$ with $n>1$ are not RTFN, 
since they may be residually $p$ for some primes $p$, but not for all 
primes (see~Example~\ref{ex:bs-p}). 
\end{example}

\begin{example}
\label{ex:free_groups_are_RTFN}
Finitely generated free groups are known to be RTFN (\cite{Magnus35}; see also~\cite{Baumslag2010}). 
As a consequence, they are residually $p$ for every prime $p$ (a result that can also be obtained 
directly from a mod-$p$ version of the Magnus expansion argument used to show that they are RTFN). 
From this, we deduce that all free groups are RTFN and residually $p$ for all $p$. In fact, if finitely 
generated free groups are residually $\mathcal R$, then all free groups are. Indeed, suppose that 
finitely generated free groups are residually $\mathcal R$, and let $X$ be any set. If $w \neq 1$ 
belongs to the free group $F[X]$ on $X$, the subset $I \subseteq X$ of letters appearing in the 
word $w$ is finite, and then the canonical projection $F[X] \surj F[I]$ does not kill $w$, since it 
sends it to the same word $w \neq 1$. By assumption, there is a quotient of $F[I]$ (hence 
of $F[X]$) satisfying $\mathcal R$ where $w$ survives. This shows that $F[X]$ is residually $\mathcal R$.
\end{example}

\begin{example}
We will see in \cref{sec:braids} that the group of pure braids on a surface different from the sphere is 
always $RTN_2$, hence it is residually $2$ if it is finitely generated, which is the case for surfaces of 
finite type. However, it is RTFN only for orientable surfaces. For non-orientable surfaces of finite type, 
it is residually $2$ but not residually $p$ for any prime $p \neq 2$.
\end{example}

Recall that if $\pi$ is a set of primes, a $\pi$-number is a product of elements of $\pi$, and a 
$\pi$-element of a group is an element whose order is a $\pi$-number. A group is called a 
$\pi$-group or a $\pi$-torsion group (resp.~a $\pi$-free group) if all (resp.~none) of its 
non-trivial elements are $\pi$-elements. A nice consequence of Gruenberg's theorem 
(also a refinement thereof) is the following result, whose proof occupies the rest of the 
present section:
\begin{theorem}
\label{thm:res_pi-torsion_vs_res_pi'-free}
Let $\pi$ be a nonempty set of primes and let $\pi'$ be its complement. A finitely generated 
group $G$ is residually nilpotent $\pi$-torsion if and only if it is residually nilpotent $\pi'$-free.
\end{theorem}

\begin{remark}
A finitely generated nilpotent $\pi$-torsion group is always finite, by a straightforward 
generalization of Proposition~\ref{prop:gf_q-grps_are_finite}.
\end{remark}

\begin{corollary}
Let $q \geq 2$ and let $\pi_q$ be the set of primes not dividing $q$. A finitely generated group 
$G$ is residually a nilpotent $\pi_q$-free group if and only if it is residually~$q$.
\end{corollary}

In particular, this says that for finitely generated groups, the filtrations $\gamma^{\qt}$ 
and $\gamma^{\pi_q^{-1}}$ govern the same residual property. In particular:
\begin{corollary}
A finitely generated group is residually $p$ if and only if $\gamma_*^{\pi_p^{-1}}(G)$ is separating.
\end{corollary}

Before proving Theorem~\ref{thm:res_pi-torsion_vs_res_pi'-free}, we need to recall some 
classical facts about finitely generated nilpotent groups.

\begin{lemma}\label{subgroups_of_index_n}
If $G$ is a finitely generated group and $n \geq 1$ is an integer, there are finitely 
many subgroups of $G$ of index $n$.
\end{lemma}

\begin{proof}
Every subgroup $H$ of $G$ of index $n$ is the stabilizer of $H$ in the action of $G$ on 
$G/H$ by (left) translation. Since there is a bijection $G/H \cong \{1, \dots, n\}$, the 
subgroup $H$ is a stabilizer for an action of $G$ on $\{1, \dots, n\}$. Such an action 
is given by a morphism from $G$ to the symmetric group $\Sym_n$. Since $G$ is 
finitely generated, there are finitely many such morphisms, hence finitely many 
subgroups of index $n$. 
\end{proof}

\begin{proposition}
Polycyclic groups are virtually torsion-free.
\end{proposition}

\begin{proof}
We reason by induction. For cyclic groups, the result is trivial. Now, we need to show that, 
given an extension $H \inj G \surj C$, if $H$ is finitely generated and virtually torsion-free 
and $C$ is cyclic, then $G$ is virtually torsion-free. Assume that $H$ has a torsion-free 
subgroup $F$ of finite index $n$. If $C$ is finite, then $F$ is a subgroup of finite index 
of $G$, whence the result. So we are left with the case where $C$ is infinite cyclic, 
generated by $t$. Then the extension must be split (any lift of $t$ to $G$ gives a splitting), 
so $G \cong H \rtimes C$. Since all the $t^iFt^{-i}$ have index $n$ in $G$, 
Lemma~\ref{subgroups_of_index_n} ensures that $S = \bigcap_i t^iFt^{-i}$ is 
a finite intersection, so it is of finite index in $H$. Then $S \rtimes C$ is of the 
same finite index in $G$, and it is torsion-free.
\end{proof}

Since a finitely generated nilpotent group $G$ is obtained as a tower of extensions 
from the finitely generated abelian groups $\gamma_iG/\gamma_{i+1}G$, it is polycyclic. 
Hence the:
\begin{corollary}
\label{fg_nilp_virt_TF}
Finitely generated nilpotent groups are virtually torsion-free.
\end{corollary}

\begin{proof}[Proof of Theorem~\ref{thm:res_pi-torsion_vs_res_pi'-free}]
A $\pi$-group is certainly $\pi'$-free, so residually nilpotent $\pi$-groups are clearly residually 
nilpotent $\pi'$-free. Conversely, suppose that $G$ is residually nilpotent $\pi'$-free. 
If $g \in G \setminus \{1\}$, there is a nilpotent $\pi'$-free quotient $N$ of $G$ 
where $\overline g \neq 1$. Let $T \coloneqq \Tors(N)$ be its torsion subgroup 
(see Lemma~\ref{torsion_subgroup}). Let us construct a quotient of $N$ which is 
a nilpotent $\pi$-group where the image of $\overline g$ is not trivial.

If $\overline g \notin T$, then its image in $N/T$ is not trivial. Since $N/T$ is torsion-free, 
Gruenberg's theorem ensures that it is residually $p$ for any $p \in \pi$. In particular, 
there is a quotient of $N/T$ (hence of $G$) which is a finite $p$-group (hence a 
nilpotent $\pi$-group) where the image of $\overline g$ (hence of $g$) is not trivial, 
whence the result in this case.

If $\overline g \in T$, then by Corollary~\ref{fg_nilp_virt_TF}, there exists $F$ a torsion-free 
subgroup of $N$ of finite index, which can be assumed normal in $N$, up to replacing it 
with the (finite) intersection of its conjugates. Then $\overline g \notin F$ hence its 
image $[g]$ in $N/F$ is not trivial. By applying Proposition~\ref{torsion_as_prod_of_p-torsion} 
to the finite nilpotent group $N/F$, we get the decomposition $N/F = \Tors_\pi(N/F) \times 
\Tors_{\pi'}(N/F)$ (where $\Tors_\pi(G)$ is the set of $\pi$-elements of $G$, which is also, by 
Proposition~\ref{torsion_as_prod_of_p-torsion}, the restricted direct product of the 
$\Tors_p(G)$, for $g \in \pi$). Since $N$ is $\pi'$-free, $g$ is a $\pi$-element in $N$, 
hence $[g] \in \Tors_\pi(N/F)$. Hence $g$ has a non-trivial image in the quotient 
$\Tors_\pi(N/F)$ of $N/F$ (hence of $G$), which is a nilpotent $\pi$-group. This 
gives the conclusion in this case, and concludes the proof.
\end{proof}

\subsection{Torsion conditions and unipotent actions}
\label{subsec:residual_prop_sd_prod_unipotent}

Let $R$ be a ring. An $R$-linear action of a group $K$ on an $R$-module $V$ 
is called \emph{unipotent}\/ if there is a finite filtration by 
sub-$R$-modules $V = V_1 \supseteq V_2 \supseteq \cdots 
\supseteq V_{l+1} = \{0\}$ stable under the action of $K$, 
and such that the induced action of $K$ on the each quotient 
module $V_i/V_{i+1}$ is trivial. 

Let $R$ be either the ring $\Z[\pi^{-1}]$ for a set $\pi$ of primes, 
or the ring $\Z/q$ for some integer $q \geq 2$. Let $\mathcal P = \mathcal P_R$ 
be, respectively, the property ``being $\pi$-free central'' or 
``being $q$-torsion central'', which, in both cases, means that 
in addition to the centrality of $G_*$, the natural map 
$\gr(G_*)\to \gr(G_*) \otimes R$ is injective. Recall that 
$\mathcal P_R$ satisfies \ref{H0}--\ref{H4} (see~\cref{sect:lcs-rational}, 
in particular~\cref{sec:pi-free} for $R = \Z[\pi^{-1}]$, 
and~\cref{sec_q-torsion} for $R =\Z/q$). For every group $G$, 
let us denote $\gamma_*^{\mathcal P}(G)$ by $\gamma_*^R(G)$ 
and by $\gr_*^R(G)$ its associated graded. Let us also denote 
by $\mathcal R_R \coloneqq \mathcal R_{\mathcal P}$ the 
class of groups satisfying $\gamma_n^R = 1$ for for $n \geq 1$. 
Recall that residually $\mathcal R_R$ groups are:
\begin{itemize}[itemsep = 1.5pt]
\item residually nilpotent groups if $R = \Z$; 
\item RTFN groups if $R = \Q$; 
\item RTN$_q$ groups if $R = \Z/q$;
\item Residually $p$ groups if $R = \F_p$ and the groups in question 
are finitely generated.
\item residually $\pi$-free nilpotent groups if $R = \Z[\pi^{-1}]$ 
(see Remark~\ref{rk:res_pi-free_nilp});
\end{itemize}

The following result generalizes of one implication of Corollaries~\ref{cor:rtfn-ext} and~\ref{cor:rtnq-ext}: 

\begin{theorem}
\label{Residual_prop_of_sdp}
Let a group $K$ act on a group $H$ by automorphisms. If $H$ and $K$ 
are residually $\mathcal R_R$ and the action of $K$ on $H^{\ab} \otimes R$ 
is unipotent, then the semi-direct product $H \rtimes K$ is residually 
$\mathcal R_R$. 
\end{theorem}

\begin{remark}
Theorem~\ref{Residual_prop_of_sdp} and its proof generalize a result of Hall \cite[Lem.~7]{Hall58}, 
which is an important ingredient of~\cite{Baumslag2010}, where it is used to show that some 
groups of the form $N \rtimes \Z$ are RTFN (using the fact that the LCS of $N$ is torsion-free 
to avoid working over $\mathbb Q$, which our theorem would allow). 
\end{remark}

\begin{proof}[Proof of Theorem~\ref{Residual_prop_of_sdp}]
We are going to construct a separating filtration 
on $H \rtimes K$ satisfying $\mathcal P_R$. Then it will contain 
$\gamma_*^R(H \rtimes K)$, whence the result. This filtration 
will be of the form $H_* \rtimes \gamma_*^R(K)$, where $H_*$ is 
a subdivision of $\gamma_*^R(H)$ on which $\gamma_*^R(K)$ acts. 
Then the fact that $H$ and $K$ are residually $\mathcal R_R$ 
will clearly imply that $H_* \rtimes \gamma_*^R(K)$ is separating, 
since its intersection is contained in the intersection 
(over $i,j \geq 1$) of the subgroups $\gamma_i^R(H) \rtimes \gamma_j^R(K)$.

The main ingredient in the construction of such a filtration $H_*$ is 
Lemma~\ref{Unipotent_action_on_gr}, which implies that under our hypothesis, 
there is, for each $i \geq 1$, a finite filtration of $R$-modules 
$\gr_i^R(H) \otimes R = W_{i1} \supseteq W_{i2} \supseteq \cdots 
\supseteq W_{i, l(i)} = \{0\}$, stable under the action of $K$, 
and such that the action of $K$ on the quotient $W_{ij}/W_{i, j+1}$ is trivial. 
Then let $H_{ij}$ be the inverse image of $W_{ij}$ under the canonical 
map $\eta_i\colon \gamma_i^R(H) \surj \gr_i^R(H) 
\inj \gr_i^R(H) \otimes R$. Finally, let $H_*$ 
denote the filtration:
\begin{align*}
H = H_{11} \supseteq H_{12} \supseteq 
&\cdots \supseteq H_{1, l(1)} = \gamma_2^R(H) = H_{21} \supseteq H_{22} \supseteq \cdots \\
&\cdots \supseteq H_{i, l(i)} = \gamma_{i+1}^R(H) = H_{i+1,1} \supseteq H_{i+1, 2} 
\supseteq \cdots. 
\end{align*}
Clearly, $H_*$ is a refinement of the central filtration $\gamma_*^R(H)$, 
so it is central. For each $i,j \geq 1$, the map $\eta_i$ induces 
an injection $\eta_{ij}$ of $H_{ij}/H_{i, j+1}$ into $W_{ij}/W_{i, j+1}$. 
These are $\Z$-linear injections into $R$-modules, so $\gr(H_*)$ 
injects into $\gr(H_*) \otimes R$. Thus $H_*$ satisfies $\mathcal P_R$. 
Moreover, the maps $\eta_i$ are $K$-equivariant, so $H_*$ is stable under 
the action of $K$. In addition, the $K$-equivariance of the induced 
injection $\eta_{ij}$ implies that the action of $K$ on $\gr(H_*)$ 
is trivial. In other words, $K = \mathcal A_1(K, H_*)$, so (because 
of \ref{H0} and \ref{H1} for $\mathcal P_R$), $\gamma_*^R(K) \subseteq 
\mathcal A_*(K, H_*)$, which means that $\gamma_*^R(K)$ acts on $H_*$. 
Then (because of \ref{H2}) the filtration $H_* \rtimes \gamma_*^R(K)$ 
satisfies $\mathcal P_R$, and it is separating, as required.
\end{proof}

\begin{lemma}
\label{Unipotent_action_on_gr}
If the action of $K$ on $H^{\ab} \otimes R$ is unipotent, then its action 
on each $\gr_i^R(H) \otimes R$ is also unipotent.
\end{lemma}

\begin{proof}
The proof is somewhat different for the two cases (when $R$ is a localisation 
or a quotient of $\Z$), but it uses the same idea of generation in degree one. 
In both cases, we reason by induction on the degree $i$, the result being true 
for $\gr_1^R(H) \otimes R \cong H^{\ab} \otimes R$ by hypothesis. Let us denote by 
$H^{\ab} \otimes R = V_1 \supseteq V_2 \supseteq \cdots \supseteq V_{l+1} = \{0\}$ 
the filtration witnessing the unipotence of the action on $H^{\ab} \otimes R$. 
Let us now assume that we already have such a filtration in degree $i$ (for some 
fixed $i \geq 1$), and denote it by 
$\gr_i^R(H) \otimes R = W_1 \supseteq W_2 \supseteq \cdots \supseteq W_{k+1} = \{0\}$. 

If $R = \Z[\pi^{-1}]$, let us recall that $\gamma_*^R = \ssqrt[\pi]{\gamma_*}$, 
so \eqref{gr(G*)_pi-free}, together with Proposition~\ref{engdeg1}, implies that 
$L_* \coloneqq \gr_*^R(H) \otimes R \cong \gr(H) \otimes R$ is generated in degree 
one as a Lie algebra over $R$. Thus the Lie bracket $L_1 \otimes_R L_i \to L_{i+1}$ 
is surjective, so if we pose:
\[
U_j \coloneqq \sum\limits_{s+t-1 = j} [V_s, W_t],
\]
we have a filtration on the $R$-module $U_1 = L_{i+1}$ satisfying $U_{l+k+1} = 0$. 
Moreover, let $x \in V_s$ and $y \in W_t$ (so that $[x,y] \in U_{s+t-1}$), and let 
$k \in K$. By definition of $V_*$ and $W_*$, we have $k \cdot x \in x + V_{s+1}$ and 
$k \cdot y \in y + W_{t+1}$. Since $K$ acts by Lie algebra automorphisms, we have:
\[
k \cdot [x,y] \in [x + V_{s+1}, y + W_{t+1}] \subseteq [x,y] + U_{s+t}.
\]
It follows that for every $j \geq 1$ and every $z \in U_j$, $k \cdot z \in z + U_{j+1}$, 
which finishes the proof of the induction step, and of our statement for 
$R = \Z[\pi^{-1}]$.

If $R = \Z/q$, the same argument works, but we have to take $q$-powers into account: 
Lie brackets are not enough to generate $L_* = \gr_*^R(H) \otimes R$ (which is just 
$\gr_*^R(H)$ in this case) from $L_1$. Namely, the inductive description of 
$\gamma_*^R = \gamma_*^{\qt}$ \eqref{induction_formula_for_LCS_q-torsion} gives 
that $L_{n+1} = L_n^q + [L_1, L_n]$, where $L_n^q$ is the abelian subgroup of 
$L_{n+1}$ generated by the elements $x^q$, with $x \in L_n$. Let us define:
\[
U_j \coloneqq W_j^q + \sum\limits_{s+t-1 = j} [V_s, W_t].
\]
In order to show that $k \cdot z \in z + U_{j+1}$ for every $k \in K$, every 
$j \geq 1$, and every $z \in U_j$, we need to decompose $z$ as a sum of 
brackets and $q$-powers. Brackets are treated as in the previous case, 
so we are left with showing that:
\[
k \cdot x^q \in x^q + U_{j+1}, \:  
\text{for all $x \in W_j$ and $k \in K$.}
\] 
This is a relation in the central subgroup $L_{n+1}  = \gr^{\qt}_{n+1}(H)$ of $H/\gamma^{\qt}_{n+2}$. 
Let us consider the canonical projection $\pi$ from $H/\gamma^{\qt}_{n+2}$ onto 
$H/\gamma^{\qt}_{n+1}$, of which $L_n$ is a central subgroup. Let us choose 
a lift $\tilde x \in \pi^{-1}(x)$ of $x$ modulo $\gamma^{\qt}_{n+2}$. 
By definition of the $q$-power operation, the element $x^q \in L_{n+1}$ 
is none other that $\tilde x^q$ (a power computed in $H/\gamma^{\qt}_{n+2}$). 
Then $k \cdot \tilde x^q = (k \cdot \tilde x)^q$, and since $k \cdot x \in 
x + W_{j+1}$ in $L_n$, we have $k \cdot \tilde x = \tilde x w$ for some 
$w \in \pi^{-1}(W_{j+1})$. 
Now, notice that $U_{j+1}$, viewed as a subgroup of $H/\gamma^{\qt}_{n+2}$, 
contains the commutator $[H/\gamma^{\qt}_{n+2}, \pi^{-1}(W_j)]$, which 
corresponds to the bracket $[V_1, W_j]$ in the associated graded. 
Thus, $[\tilde x,w] \in U_{j+1}$, which means that $\tilde x$ 
and $w$ commute modulo $U_{j+1}$, hence $(\tilde x w)^q \equiv
\tilde x^q w^q$ modulo $U_{j+1}$. Since $w^q \in U_{j+1}$, 
this finishes the proof that $k \cdot \tilde x^q \equiv \tilde x^q$ 
modulo $U_{j+1}$, that is, that $k \cdot x^q \equiv x^q$ modulo 
$U_{j+1}$ in $L_{n+1}$.
\end{proof}

\section{Application to pure braids on surfaces}
\label{sec:braids}

We now apply the tools and methods introduced above to pure braids on surfaces. 
Most of these groups admit a decomposition as an iterated semi-direct products 
of free groups and fundamental groups of closed surfaces. We use this to prove 
the following general result, which is a combination of Theorem~\ref{Pn(S)_RTFN} 
(for non-closed orientable surfaces), Theorem~\ref{Pn(S)_res-2} and 
Corollary~\ref{cor:nor-orientable-non-res} (for non-closed non-orientable surfaces),  
and Theorem~\ref{Pn(closed_S)} (for closed surfaces):
\begin{theorem}
\label{res_prop_of_Pn(S)}
Let $S$ be a connected surface different from the sphere and the projective plane, 
and let $n \geq 2$. Then the group $P_n(S)$ of pure braids on $S$ is:
\begin{itemize}[itemsep = 2pt]
    \item RTFN and RTN$_q\ $ for all $q \geq 2$ if $S$ is orientable.
    \item RTN$_2$ but not RTFN nor RTN$_q\ $ for any odd $q \geq 3$ if $S$ is not orientable.
\end{itemize}
\end{theorem}

In order to shed some light on the meaning of this statement, recall that 
if a prime $p$ divides an integer $q$, we have the implications:
\[
\text{residually } p \ \ \Rightarrow\ \ \mathrm{RTN}_p\ (= \text{residually }\mathcal R_{\F_p}) \ \ \Rightarrow\ \ \mathrm{RTN}_q\ (= \text{residually }\mathcal R_{\Z/q}) .
\]
The first one is Corollary~\ref{res_p_vs_RTN_p}, which also says that the converse 
holds for finitely generated groups. The second one is clear, and it implies that being 
RTN$_q$ for all $q \geq 2$ is the same as being RTN$_p$ for every prime~$p$. 
The link with properties $\mathcal R_R$ from~\cref{subsec:residual_prop_sd_prod_unipotent} 
is important, since Theorem~\ref{Residual_prop_of_sdp} will be the main ingredient in the 
proof of Theorem~\ref{res_prop_of_Pn(S)}.

For infinitely generated groups, the above implications are not equivalences, as 
Example~\ref{ex:nilp_p-gp_not_res_fini} shows. In particular, we cannot infer 
that $P_n(S)$ is always residually $p$ if $S$ is orientable, nor that $P_n(S)$ 
is always residually $2$ if $S$ is not orientable. It may be the case 
(note that infinitely generated free groups are residually $p$ -- see 
Example~\ref{ex:free_groups_are_RTFN}), and maybe a refinement of our 
techniques could be devised in order to show it, but we were unable to do so.

However, for surfaces of finite type (that is, surfaces homotopy equivalent 
to compact surfaces or, equivalently, surfaces whose fundamental group is 
finitely generated), the groups $P_n(S)$ are finitely generated, so, thanks 
to Corollary~\ref{res_p_vs_RTN_p}, Theorem~\ref{res_prop_of_Pn(S)} takes 
the more familiar form:
\begin{corollary}
\label{res_prop_of_Pn(S)_fg}
Let $S$ be a connected surface of finite type different from the sphere 
and the projective plane, and let $n \geq 2$. Then the group $P_n(S)$ 
of pure braids on $S$ is:
\begin{itemize}[itemsep = 2pt]
    \item RTFN (hence residually~$p$ for all prime~$p$) if $S$ is orientable.
    \item residually~$2$ but not residually~$p$ for any prime $p \neq 2$ 
    if $S$ is not orientable.
\end{itemize}
\end{corollary}

This result was already known, except for the negative part, which answers 
a question of \cite{Bellingeri-Gervais}. However, our proof is new, and it 
does not need any explicit computations with presentations. Precisely, 
since a surface of finite type is isotopy equivalent to some $\Sigma_{g,p}$, 
this is a combination of:
\begin{itemize}[itemsep = 2pt]
    \item Bellingeri, Gervais, and Guaschi's \cite[Thm.~4]{Bellingeri-Gervais-Guaschi} for
     orientable non-closed surfaces. Their proof used an explicit embedding into the 
    Torelli group of an orientable surface, known to be RTFN by a result of Hain.
    \item Bardakov and Bellingeri's \cite[Cor.~7]{Bardakov-Bellingeri} for orientable closed surfaces. 
    We borrow their method to deduce it from the previous one, but we are able to circumvent their calculations.
    \item Bellingeri and Gervais' \cite[Thm.~1.1]{Bellingeri-Gervais} for non-orientable surfaces. 
    They adapted the method from the previous papers about orientable surfaces, using a result 
    of Paris about the mod-$2$ Torelli group of non-orientable surfaces in place of the result of Hain.
\end{itemize}

\begin{remark}
In Corollary~\ref{res_prop_of_Pn(S)_fg}, Gruenberg's theorem (Cor.~\ref{Gruenberg_thm}) 
applies to the finitely generated group $P_n(S)$, which is RTFN, hence residually $p$ 
for all prime $p$. However, this conclusion can also be reached directly from the 
semi-direct product decomposition, without the use of Gruenberg's theorem. The latter 
approach is in fact necessary for the proof of the more general Theorem~\ref{res_prop_of_Pn(S)}.
\end{remark}

Theorem~\ref{res_prop_of_Pn(S)} is the main goal of this section. 
In particular, we focus on residual properties, and not on the 
computation of the associated graded of the corresponding filtrations, 
which would require the introduction of some more tools.

\begin{remark}
\label{pi_1(surfaces)}
For $n = 1$, that is, for $\pi_1(S)$, the statement is somewhat different: 
$\pi_1(S)$ is always RTFN, unless $S$ is the projective plane or the Klein 
bottle (in which case it is residually $2$ but not residually $p$ for any 
$p \neq 2$). Indeed, in all the other cases, either it is a free group, 
which is RTFN and residually $p$ for every prime $p$ 
(see Example~\ref{ex:free_groups_are_RTFN}), or $S$ is a closed surface 
having a handle, and then \cite{Baumslag2010} says that $\pi_1(S)$ is 
RTFN (hence residually $p$ for every prime $p$ by Gruenberg's theorem).
\end{remark}

\subsection{General results}
The reader is referred for instance to~\cite[Chap.~6]{DPS} for an introduction 
to braids on surfaces. In particular, we assume that $S$ is any (second-countable) 
connected surface (not necessarily orientable nor compact), and we consider surfaces 
up to isotopy equivalence, since an isotopy equivalence between $S$ and $S'$ induces 
an isomorphism between the corresponding braid groups. In particular, since $S$ and 
$S \setminus \partial S$ are isotopy equivalent, we can assume that 
$\partial S = \varnothing$, and use only the classification of surfaces 
without boundaries from~\cite{Richards}. 

Let us recall the computation of $P_n(S)^{\ab}$, which is the case $\lambda = (1, \dots, 1)$ 
in~\cite[Prop.~6.47]{DPS}. Let us first recall that if we fix an inclusion of a closed disk 
$\mathbb D$ into $S$ and a base configuration $\mathbf c = (c_1, \dots, c_n)$ in the interior 
of $\mathbb D$, we get an associated morphism from the usual pure braid group:
\[
\varphi \colon P_n = P_n(\D) \longrightarrow P_n(S).
\]
Moreover, if we associate to a pure braid the list of homotopy classes of its 
strands (which are loops on $S$), we get a morphism:
\[
\pi_S \colon P_n(S) \longrightarrow \pi_1(S)^n.
\]
Recall also that a surface is called~\emph{planar} if it injects into the plane, 
which is equivalent to it having no handle nor crosscap, and then it is homeomorphic 
to the plane minus a subset of a Cantor set, so we can think of it as a ``disk with 
(possibly infinitely many) holes" \cite{Richards}.
\begin{proposition}[{\cite[Prop.~6.47, for $\lambda = (1, \dots , 1)$]{DPS}}]
\label{Pn(S)^ab}
If $S$ is not isotopy equivalent to the sphere $\mathbb S^2$, then:
\[P_n(S)^{\ab} \cong
\begin{cases}
H_1(S)^n \oplus P_n^{\ab} &\text{if $S$ is planar}, \\
H_1(S)^n &\text{else.}
\end{cases}\]
Explicitly: in both cases, the projection $P_n(S)^{\ab}\surj H_1(S)^n = (\pi_1(S)^n)^{\ab}$ 
is induced by $\pi_S$, and for planar surfaces, the injection $P_n^{\ab} \inj P_n(S)^{\ab}$ 
is induced by $\varphi$.
\end{proposition}

Let $\mathbb D$ be a closed disk in $S$, and let 
$Q_m$ be a set of $m$ points in the interior of $\mathbb D$. The main tool 
for obtaining semi-direct product decompositions of pure braid groups on 
surfaces is the Fadell--Neuwirth exact sequence, which was introduced by 
Fadell and Neuwirth ~\cite[Thm.~3]{Fadell-Neuwirth} (see also 
\cite[Prop.~1.3]{Birman}, and~\cite[Prop.~6.15]{DPS} for the present statement):
\begin{proposition}[Fadell--Neuwirth exact sequences]
\label{Fadell-Neuwirth}
If $S$ is a connected surface different from the sphere and the projective plane, 
the following sequence of canonical maps is a short exact sequence, for all $m,n \geq 1$:
\begin{equation}
\tag{$FN_{m,n}$}
\label{FN}
\begin{tikzcd}[column sep=20pt]
1 \ar[r]& P_n(S \setminus Q_m) \ar[r]& P_{m+n}(S) \ar[r]& P_m(S) \ar[r]& 1.
\end{tikzcd}
\end{equation}
Furthermore, if $S$ is not closed, then the surjection $P_{m+n}(S) \surj P_m(S)$ splits.
\end{proposition}

The projections $\pi_S$ induces a morphism of extensions with source $(FN_{n-1,1})$, 
that is (recalling that $P_1(S) = \pi_1(S)$), a commutative diagram
\[
\begin{tikzcd}
\pi_1(S \setminus Q_{n-1}) \ar[r, hook] \ar[d, two heads, dotted, "\pi_S"]
&P_n(S) \ar[r, two heads] \ar[d, two heads, "\pi_S"]
&P_{n-1}(S) \ar[d, two heads, "\pi_S"] \\
\pi_1(S)  \ar[r, hook]
&\pi_1(S)^n  \ar[r, two heads]
&\pi_1(S)^{n-1}, \ar[l, bend right]
\end{tikzcd}
\]
where we still denote by $\pi_S$ the induced homomorphism between kernels. 
This implies that the map $p \coloneqq \pi_S^{\ab}$ from $\pi_1(S \setminus Q_{n-1})^{\ab} = 
H_1(S \setminus Q_{n-1})$ to $\pi_1(S)^{\ab} = H_1(S)$ is equivariant with 
respect to the action of $P_{n-1}(S)$ on its source, the action of 
$\pi_1(S)^{n-1}$ on its target, and the map $\pi_S \colon P_{n-1}(S) 
\surj \pi_1(S)^{n-1}$ (see~\cref{sec_ppd}). Moreover, the bottom 
extension corresponds to the trivial action of $\pi_1(S)^{n-1}$ on $\pi_1(S)$. 

The map $p \colon H_1(S \setminus Q_{n-1}) \to H_1(S)$ is induced by the 
quotient of $S \setminus Q_{n-1}$ by $\mathbb D \setminus Q_{n-1}$, so 
it is part of the long exact sequence associated with the pair 
$(S \setminus Q_{n-1}, \mathbb D \setminus Q_{n-1})$. Its kernel $K$ 
is thus the image of $H_1(\mathbb D \setminus Q_{n-1}) \cong \Z^{n-1}$ 
in $H_1(S \setminus Q_{n-1})$. The above can be reformulated by saying 
that the two-step filtration
\[
0 \subseteq K \subseteq H_1(S \setminus Q_{n-1})
\] 
is stable under the action of $P_{n-1}(S)$, which acts trivially on 
$H_1(S \setminus Q_{n-1})/K \cong H_1(S)$. Hence it will be a unipotent action 
(in the sense of~\cref{subsec:residual_prop_sd_prod_unipotent}) if we show 
that $P_{n-1}(S)$ acts trivially on $K$.

Let $x_i \in \pi_1(\mathbb D \setminus Q_{n-1})$ be the homotopy class of the 
loop $\eta_i C_i \eta_i^{-1}$, where $\eta_i$ is a path in $\mathbb D \setminus Q_{n-1}$ 
from the basepoint to a neighborhood of the $i$-th puncture $c_i$, and $C_i$ is an 
arbitrary small circle around the puncture. Then the group $\pi_1(\D \setminus Q_{n-1})$ 
is free on the $n-1$ generators $x_1, \dots, x_{n-1}$. 
Consider the morphism from $\pi_1(\mathbb D \setminus Q_{n-1})$ to 
$\pi_1(S \setminus Q_{n-1})$, which injects into $P_n(S)$ as above. 
We still denote by $x_i$ the image of $x_i$ in $P_n(S)$.

Let us call \emph{two-sided} (resp.~\emph{one-sided}) a loop on a surface 
along which orientation is preserved (resp.~reversed). Notice that such a 
terminology does not require any theorem about the existence of tubular 
neighbourhoods to make sense. 

\begin{lemma}
\label{basis-conj}
Let $\beta \in P_n(S)$. For $i = 1, \dots, n-1$, there exists $w_i \in \pi_1(\mathbb D 
\setminus Q_{n-1})$ such that $\beta x_i \beta^{-1} = w_i x_i^\epsilon w_i^{-1}$, 
where $\epsilon = +1$ if the $i$-th strand of 
$\beta$ is two-sided, else $\epsilon = -1$.
\end{lemma}

\begin{proof}
Let $I = [0,1]$ be the unit interval. For a pure braid $\alpha \in P_n(S)$, 
we also denote by $\alpha \colon I \to \Conf_n(S)$, $t \mapsto 
(\alpha_1(t), \dots, \alpha_n(t))$ a path representing the braid. 

Recall that if the group $\Homeo(S)$ of self-homeomorphisms of $S$ is endowed with 
the compact-open topology, then the evaluation map $\Homeo(S) \to \Conf_n(S)$ 
(sending $\varphi$ to $\varphi(\mathbf c)$) is a locally trivial fibration 
(see Birman~\cite[Thm.~4.1]{Birman}, whose proof works verbatim on any surface). 
The path $\beta \colon I \to \Conf_n(S)$ can then be lifted to a path 
$I \to \Homeo(S)$ ending at~$\id_S$. In other words, there exists a path 
$t \mapsto \varphi_t$ of homeomorphisms of $S$ such that $\varphi_1 = \id_S$ 
and $\varphi_t(c_i) = \beta_i(t)$ for every $t \in I$. We also denote by 
$\Phi$ the self-homeomorphism of $I \times S$ defined by $(t,x) \mapsto 
(t, \varphi_t(x))$. Then for every braid $\alpha$, we have 
$\beta \alpha \beta^{-1} \simeq \varphi_0(\alpha)$, 
via the following explicit homotopy $H \colon I^2 \to \Conf_n(S)$:
\[
H(s,t) = 
\begin{cases}
\beta(3st)  &\text{for}\ \ t \leq \frac s3,\\
\varphi_{s^2}\!\left(\alpha\!\left(\frac{3t-s}{3-2s}\right)\right) 
&\text{for}\ \ \frac s3 \leq t \leq 1 - \frac s3,\\
\beta(3s(1-t)) &\text{for}\ \ t \geq 1 - \frac s3.
\end{cases}
\]

Since $\beta$ is a pure braid, the homeomorphism $\varphi_0$ fixes 
the points $c_i$, so it sends a disk around it homeomorphically to 
a neighborhood of $c_i$. Up to a local isotopy fixing $c_i$, we 
can assume that $\varphi_0$ equals $\pm \id$ on a small disk around $c_i$. 
Moreover, we can assume that this small disk contains the arbitrarily 
small circle $C_i$ from above, so $\varphi_0(C_i) = C_i^{\pm 1}$. 
Now, with $x_i = \eta_i C_i \eta_i^{-1}$ introduced above, seen 
as a braid with only $c_1$ moving and $c_2, \dots, c_n$ staying 
put, we have  
\[
\beta x_i \beta^{-1} \simeq \varphi_0(x_i) = \varphi_0(\eta_i) 
\varphi_0(C_i) \varphi_0(\eta_i^{-1}) = w_i x_i^{\pm 1} w_i^{-1},
\]
if we pose $w_i \coloneqq \varphi_0(\eta_i)\eta_i^{-1} 
\in \pi_1(\D \setminus Q_{n-1})$. 

In order to finish the proof, we need to show that the sign is $+1$ 
if and only if orientation is preserved along $\beta_i$. This is quite 
intuitive, but we give a formal proof by considering the effect 
on orientation covers of the following commutative diagram:
\[
\begin{tikzcd}
&I \ar[d, "{(\id,\ \beta_i)}"] \ar[rd, bend left, "\beta_i"] 
\ar[ld, swap, "{(\id,\ c_i)}"]
&\\
I \times S \ar[r, "\Phi"]
&I \times S \ar[r, "\proj_2"]
&S.
\end{tikzcd}
\]

Let us denote by $\widetilde M$ the $2$-fold cover of a manifold $M$ obtained 
from local orientations (note that we need to consider manifolds with 
boundaries -- see for instance~\cite[Chap.~21, \S~4]{May} for technical 
details). We then have a canonical homeomorphism 
$\widetilde{I \times S} \cong I \times \widetilde S 
\cong \proj_2^*(\widetilde S)$. Recall that a self-homeomorphism $u$ 
of a manifold $M$ induces an automorphism $\tilde u$ of the covering 
$\widetilde M$. If $f$ is a morphism of coverings of a space $X$ and 
$x \in X$, then we denote by $f_x$ ($= f|_{\{x\}}$) the restriction 
of $f$ to the fiber over $x$. 

On the one hand, by construction, the sign in the 
previous equation is $+1$ if and only if $\widetilde{\varphi_0} \colon 
\widetilde S \to \widetilde S$ acts trivially on the fiber above $c_i$, 
that is, if and only if $(\widetilde{\varphi_0})_{c_i} = \id$. On the other 
hand, orientation is preserved by $\beta_i$ if and only if the unique 
covering morphism $f \colon \widetilde I \to \beta_i^*(\widetilde S)$ satisfying 
$f_1 = \id$ also satisfies $f_0 = \id$.

Let us show that the morphism $f$ is the one obtained by pulling back 
$\widetilde \Phi$ to coverings of $I$. Precisely, $\widetilde \Phi$ 
is pulled back to a morphism $g$ from $(\id, c_i)^*(\widetilde{I \times S})$ 
to $(\id, \beta_i)^*(\widetilde{I \times S})$. Using the identification 
$\widetilde{I \times S}\cong \proj_2^*(\widetilde S)$, we see that the 
source of $g$ is canonically identified with $\widetilde I = c_i^*(\widetilde S)$ 
(where $c_i$ denotes the constant map equal to $c_i$), and its target is identified to 
$\beta_i^*(\widetilde S)$. Since $g_1$ identifies with $\widetilde \Phi_{(1,c_i)}$, it
is the identity, because $\Phi|_{1 \times S} = \varphi_1 = \id_S$ (using 
once again the canonical identification between local orientations 
near the boundary and local orientations of the boundary 
\cite[Chap.~21, \S~4]{May}). Hence $g = f$, by unicity of $f$, as announced. 
But $g_0 = \widetilde \Phi_{(0,c_i)}$, which identifies canonically with 
$(\widetilde{\varphi_0})_{c_i}$. This shows that $f_0 = (\widetilde{\varphi_0})_{c_i}$, 
so one is the identity if and only if the other is, whence the conclusion.
\end{proof}

\begin{remark}
In the previous proof, the isotopy class of $\varphi_0$ is entirely determined 
by the braid $\beta$, since it is, by construction, its image by the connecting 
morphism $P_n(S) = \pi_1(\Conf_n(S)) \rightarrow \pi_0(\Homeo(S))$ in the long 
exact sequence in homotopy associated with the above fibration. This is the usual 
construction associating a mapping class to a braid.
\end{remark}

Lemma~\ref{basis-conj} will allow us to apply Theorem~\ref{Residual_prop_of_sdp} 
to the extension $(FN_{n-1, 1})$, whenever it splits, that is, whenever $P_n(S) \cong 
\pi_1(S \setminus Q_{n-1}) \rtimes P_{n-1}(S)$. As stated in Proposition~\ref{Fadell-Neuwirth}, 
this extension does split whenever $S$ is \emph{not} a closed surface.

\subsection{Non-closed orientable surfaces}

Assume that $S$ is orientable. Then Lemma~\ref{basis-conj} 
says that the action of an element of $P_n(S)$ by conjugation sends each $x_i$ 
to a conjugate of $x_i$ in $\pi_1(S \setminus Q_{n-1})$. This implies 
immediately that the action of $P_{n-1}(S)$ on $K$ fixes the images 
of the elements $x_i$. These elements 
generate $K$ (the image of $\pi_1(\mathbb D \setminus Q_{n-1})$ 
in $H_1(S \setminus Q_{n-1})$), so the action of $P_{n-1}(S)$ 
on $K$ is trivial. Thus (using the terminology 
from~\cref{subsec:residual_prop_sd_prod_unipotent}):
\begin{corollary}
\label{unipotent_action_orientable}
The action of $P_{n-1}(S)$ on $H_1(S \setminus Q_{n-1})$ is unipotent.
\end{corollary}

In fact, this still holds after tensoring with any ring $R$. This follows from 
the same reasoning over $R$, since $H_1(S \setminus Q_{n-1}) \otimes R \cong 
H_1(S \setminus Q_{n-1},R)$ by the universal coefficient theorem. Alternatively, 
by tensoring the short exact sequence $K \inj H_1(S \setminus Q_{n-1}) \surj H_1(S)$ by $R$, 
one gets the exact sequence $K \otimes R \to H_1(S \setminus Q_{n-1})\otimes R \surj 
H_1(S)\otimes R$, so if $\overline K$ is the image of the first map, we get directly 
a $2$-step filtration $0 \subseteq \overline K \subseteq H_1(S \setminus 
Q_{n-1}) \otimes R$ with respect to which the action is unipotent. 
Either way, we get:

\begin{corollary}
\label{unipotent_action_orientable_R}
The action of $P_{n-1}(S)$ on $H_1(S \setminus Q_{n-1}) \otimes R$ is unipotent, for 
every commutative ring $R$.
\end{corollary}

Assume further that $S$ is not closed. Then the extension $(FN_{n-1,1})$ splits, which 
allows us to use our methods to prove:
\begin{theorem}\label{Pn(S)_RTFN}
Let $S$ be a connected surface that is orientable, but not closed. Then $P_n(S)$ is 
RTFN, and RTN$_q\ $ for all $q \geq 2$.
\end{theorem}

\begin{proof}
We reason by induction on $n$. If $n = 1$, then $\pi_1(S) = P_1(S)$ is a free group, so it is 
RTFN and residually $p$ for every prime $p$ (even if it is not finitely generated -- see 
Example~\ref{ex:free_groups_are_RTFN}), hence it is RTN$_q\ $ for all $q \geq 2$. 
Let us now fix $n \geq 2$, assume that $P_{n-1}(S)$ is RTFN and RTN$_q\ $ for all 
$q \geq 2$, and show that $P_n(S)$ too. We want to apply Theorem~\ref{Residual_prop_of_sdp} 
to the above action of $P_{n-1}(S)$ on $\pi_1(S \setminus Q_{n-1})$. 
Corollary~\ref{unipotent_action_orientable_R} says that the action of $P_{n-1}(S)$ 
on $\pi_1(S \setminus Q_{n-1})^{\ab} \otimes R$ is unipotent. Moreover, 
$\pi_1(S \setminus Q_{n-1})$ is a free group, so it is RTFN and RTN$_q\ $ 
for all $q \geq 2$ (using again Example~\ref{ex:free_groups_are_RTFN}), 
so all the hypotheses of Theorem~\ref{Residual_prop_of_sdp} are satisfied, 
both for $R = \Q$ and for $R = \Z/q$, whence the conclusion.
\end{proof}

\begin{remark}
The use of Theorem~\ref{Residual_prop_of_sdp} does not give a lot of control on the 
filtrations involved. For \emph{planar surfaces}, the semi-direct product 
$P_n(S) \cong \pi_1(S \setminus Q_{n-1}) \rtimes P_{n-1}(S)$ is an almost-direct one, 
which allows us to describe completely the lower central series (see~\cref{par:planar}). 
However, the decomposition $P_n(S) \cong \pi_1(S \setminus Q_{n-1}) \rtimes P_{n-1}(S)$ 
is never an almost-direct product if $S$ is not planar (see Remark~\ref{not_adp_when_not_planar}); 
in the general case, getting such a description would require the introduction of more tools, 
so we leave this work for a later paper.
\end{remark}

\subsection{Non-closed non-orientable surfaces}

\emph{Assume that $S$ is not orientable}. Then the same holds over $\mathbb F_2$. 
Namely, Lemma~\ref{basis-conj} implies that $P_{n-1}(S)$ acts trivially on the 
image $\overline K$ of $K$ inside $H_1(S \setminus Q_{n-1}, \mathbb F_2) \cong 
H_1(S \setminus Q_{n-1}) \otimes \F_2$ (the latter isomorphism follows, 
as above, from the universal coefficient theorem). As a result:
\begin{corollary}
\label{unipotent_action_not_orientable}
The action of $P_{ n-1}(S)$ on $H_1(S \setminus Q_{n-1}) \otimes \F_2$ is unipotent.
\end{corollary}

Assume that $S$ is not orientable but not closed. Then the extension $(FN_{n-1,1})$ 
splits, but the decomposition $P_n(S) \cong \pi_1(S \setminus Q_{n-1}) \rtimes P_{n-1}(S)$ 
is never a $2$-almost-direct product. Nevertheless, we can apply 
Theorem~\ref{Residual_prop_of_sdp} to deduce from 
Corollary~\ref{unipotent_action_not_orientable} the following:
\begin{theorem}\label{Pn(S)_res-2}
Let $S$ be a connected surface that is not orientable and not closed. 
Then $P_n(S)$ is RTN$_2$.
\end{theorem}

We now show that $P_n(S)$ is not residually $p$ if $p \neq 2$. 
We need the following criterion:

\begin{proposition}
\label{prop:non-res_p}
Let $G$ be a group, and $q \neq 0$ be an integer. Suppose there 
exist $g, x \in G$ such that $x \neq 1$ and $[g,x]= x^q$. If $q$ is 
invertible in $R$, then $G$ is not residually $\mathcal R_R$. 
In particular, $G$ is neither RTFN nor RTN$_{q'}$ for any integer 
$q'$ prime to $q$.  Nor is it residually~$p$ for any prime $p$ not dividing $q$.
\end{proposition}

\begin{proof}
Let us first show by induction on $k \geq 1$ that the relation $[g,x]= x^q$ 
implies $x^{q^k} \in \gamma_{k+1}(G)$. The case $k = 1$ is clear. Suppose that 
$x^{q^k} \in \gamma_{k+1}(G)$ for some $k \geq 1$. The relation can be written 
$gxg^{-1} = x^{1+q}$, hence:
\[
[g,x^{q^k}] = gx^{q^k}g^{-1}x^{-q^k} = (gxg^{-1})^{q^k}x^{-q^k} = 
(x^{1+q})^{q^k}x^{-q^k} = x^{q^{k+1}}.
\]
Thus $x^{q^{k+1}} \in [G, \gamma_{k+1}(G)] = \gamma_{k+2}(G)$, 
which finishes the induction. 

Now, if $K$ is a nilpotent quotient of $G$, then $\gamma_{c+1}(K) = 1$ 
so $x^{q^c} = 1$ in $K$. As a consequence, if $\mathcal R$ is any 
class of nilpotent groups having no non-trivial elements whose 
order divides a power of $q$, then $G$ is not residually 
$\mathcal R$, since $x = 1$ in any quotient of $G$ belonging to $\mathcal R$. 
\end{proof}

\begin{corollary}
\label{cor:nor-orientable-non-res}
For any non-orientable surface $S$ and any $n \geq 2$, the pure 
braid group $P_n(S)$ is neither RTFN nor RTN$_q$ for any odd $q \geq 3$. 
In particular, it is not residually $p$ for any prime $p \neq 2$.
\end{corollary}

\begin{proof}
Let $\beta \in P_n(S)$ be a braid whose first strand is one-sided (for instance, 
the strands $2, \dots, n$ of $\beta$ can be constant paths at $c_2, \dots, c_n$, 
and the first strand can move along any one-sided loop on $S \setminus \{c_2, \dots, 
c_n\}$). Then Lemma~\ref{basis-conj} gives a relation of the form 
$\beta x_1 \beta^{-1} = w x_1^{-1}w^{-1}$ in $P_n(S)$. This relation 
can be written as $[w^{-1}\beta, x_1] = x_1^{-2}$, so Proposition~\ref{prop:non-res_p} 
applied to $g = w^{-1}\beta$, $x = x_1$, and $q = -2$ gives the desired conclusion.
\end{proof}

We can actually compute the biggest RTFN quotient of $P_n(S)$.

\begin{proposition}
\label{prop:non-orient-quotient}
For any non-orientable surface $S$ and any $n \geq 2$, if $2$ is invertible 
in the ring $R$, then $P_n(S)/\gamma_\infty^R \cong (\pi_1(S)/\gamma_\infty^R)^n$, 
which equals $\pi_1(S)^n$, unless $S$ is the projective plane (for which is it 
trivial) or the Klein bottle (for which $\pi_1(S)/\gamma_\infty^R \cong \Z$).
\end{proposition}

\begin{proof}
Notice that the proof of Proposition~\ref{prop:non-res_p} says a little more than 
its statement. Namely, it says that whenever $2$ is invertible in $R$, the 
element $x$ must belong to $\bigcap \gamma_i^R = \gamma_\infty^R$. The proof of 
Corollary~\ref{cor:nor-orientable-non-res} then implies that $x_1 \in \gamma_\infty^R(P_n(S))$. 
But $x_1$ is none other than the classical pure braid generator $A_{1n}$ (up to conjugation). 
By inverting the roles of the strands, we see that all the generator $A_{ij}$ of $P_n$ 
(or, rather, of its image $\varphi(P_n)$ in $P_n(S)$) must belong to $\gamma_\infty^R(P_n(S))$, 
so $\varphi(P_n) \subseteq \gamma_\infty^R(P_n(S))$. Since $\gamma_\infty^R(P_n(S))$ is 
normal in $P_n(S)$, it contains the normal closure of the image of $P_n$, which, 
by~\cite[Prop.~6.10]{DPS} (generalizing a result of Goldberg), is exactly the 
kernel of the projection $\pi_S \colon P_n(S) \rightarrow \pi_1(S)^n$. As a 
consequence, $P_n(S)/\gamma_\infty^R$ is canonically identified with the 
maximal residually $\mathcal R_R$ quotient of $\pi_1(S)^n$. But $\pi_1(S)$ 
is residually $\mathcal R_R$ unless $S$ is the projective space or the Klein 
bottle (see Remark~\ref{pi_1(surfaces)}), and in the latter case, $\pi_1(S)$ 
is either $\Z/2$ or $\Z \rtimes \Z$, which are easy to deal with explicitly. 
\end{proof}

\subsection{Closed surfaces}
\label{subsec:closed-surfaces}

For closed surfaces, the extension $(FN_{n-1,1})$ does not split 
in general (see the solution of the splitting problem 
from~\cite[Thm.~2]{Gonçalvez-Guaschi-2}). 
However, a very nice topological argument of Gon\c{c}alvez and Guaschi 
(see~\cite[\S~2,~Thm.~1]{Goncalves-Guaschi} and its proof, 
which works for any surface different from $\mathbb S^2$ 
and $\mathbb P^2$) gives, for $S \neq \mathbb S^2, \mathbb P^2$, 
a splitting of the extension $(FN_{1,n-1})$. In other words,  
there is a semi-direct product decomposition:
\[
P_n(S) \cong P_{n-1}(S \setminus \{\pt\}) \rtimes \pi_1(S).
\]
In order to apply Proposition-Definition~\ref{def_ppd}, we need to study the 
corresponding group morphism $P_{n-1}(S \setminus \{\pt\})^{\ab} \to P_n(S)^{\ab}$, 
which is induced by $(z_1, \dots , z_{n-1}) \mapsto (z_1, \dots, z_{n-1}, \pt)$ 
at the level of configuration spaces. Thanks to Proposition~\ref{Pn(S)^ab}, 
we identify it with the map $H_1(S \setminus \{\pt\})^{n-1} \to H_1(S)^n$ 
given by the product of $n-1$ copies of the canonical map 
$i_*\colon H_1(S \setminus \{\pt\}) \to H_1(S)$ followed 
by the canonical inclusion $H_1(S)^{n-1} \inj H_1(S)^{n-1} \times H_1(S)$. 
As above, there are two cases:
\begin{itemize}[itemsep = 2pt]
    \item \emph{If $S$ is orientable}, then $i_*$ is an isomorphism, 
    hence the semi-direct product is an almost-direct one.
    \item \emph{If $S$ is not orientable}, then $i_*$ is the quotient 
    of $\Z^g$, with basis $a_1, \dots, a_g$, by the subgroup generated by 
    $2(a_1 + \cdots + a_g)$. Then it becomes injective only when tensored 
    by $\F_2$, hence the semi-direct product is not almost-direct, nor 
    is it $p$-almost direct for $p \neq 2$, but it is a $2$-almost direct product. 
\end{itemize}
From this, we deduce the case of closed surfaces in Theorem~\ref{res_prop_of_Pn(S)} 
(recall that for finitely generated groups, such as pure braid groups on surfaces of 
finite type, being residually~$p$ is the same as being RTN$_p$):
\begin{theorem}
\label{Pn(closed_S)}
Let $S$ be a closed connected surface different from the sphere and the 
projective plane, and let $n \geq 2$. Then the group $P_n(S)$ of pure 
braids on $S$ is:
\begin{itemize}[itemsep = 2pt]
    \item RTFN (hence residually~$p$ for all prime~$p$) if $S$ is orientable.
    \item residually~$2$ but not RTN$_q$ for any odd $q \geq 3$ 
    (hence not residually~$p$ for any prime $p \neq 2$) if $S$ is not orientable.
\end{itemize}
\end{theorem}

\begin{proof}
Apply Corollaries~\ref{cor:rtfn-ext} (in the orientable case) 
and~\ref{cor:rtnq-ext} (for $q = 2$, in the non-orientable case) 
to the above decomposition, using the case of non-closed surfaces 
(Theorems~\ref{Pn(S)_RTFN} and~\ref{Pn(S)_res-2}) to deal 
with $P_{n-1}(S \setminus \{\pt\})$. For the negative statement, 
apply Corollary~\ref{cor:nor-orientable-non-res}.
\end{proof}

\subsection{The sphere and the projective plane} 
\label{subsec:S2-RP2}

It has been known for a long time that $P_n(S)$ has torsion if and 
only if $S$ is the sphere $\mathbb S^2$ (except for $P_2(\mathbb S^2) 
\cong P_1(\mathbb S^2) = 1$) or the projective plane 
$\mathbb P^2$ \cite{vanBuskirk1966}. In fact, as observed 
in \cite[Rem.~18]{Bellingeri-Gervais-Guaschi}, since $P_3(\mathbb S^2) 
\cong \Z/2$ and $P_n(\mathbb S^2) \cong P_{n-3}(\mathbb S^2 \setminus Q_3) 
\times \Z/2$ for $n \geq 4$ \cite{Goncalves-Guaschi04}, a decomposition 
as a direct product of a $2$ group and a residually $2$ group, the group 
$P_n(\mathbb S^2)$ is residually $2$ but not RTN$_q$ for any odd $q \geq 3$ 
(hence not residually $p$ for any prime $p \neq 2$).

As for $P_n(\mathbb P^2)$, it is known not to be RTN$_q$ for any 
odd $q \geq 3$ (hence not residually $p$ for any prime $p \neq 2$), 
which can be seen as part of Corollary~\ref{cor:nor-orientable-non-res}, 
or as a consequence of the existence of $2$-torsion. However, whether 
it is residually $2$ is still an open question 
for $n \geq 4$, as discussed at the end of~\cite{Bellingeri-Gervais}.

\subsection{Planar surfaces}
\label{par:planar}
Among orientable surfaces, we can distinguish the particular case of planar 
surfaces (also known as ``disks with holes"), for which the situation is 
particularly nice, allowing for a computation of the Lie ring associated with $P_n(S)$.

Let $S$ be a planar surface (notice that it is then orientable 
and not closed). Then $\pi_1(S \setminus Q_{n-1})$ is a free 
group and Proposition~\ref{Fadell-Neuwirth} says that 
$P_n(S) \cong \pi_1(S \setminus Q_{n-1}) \rtimes P_{n-1}(S)$, 
since planar surfaces without boundary are not compact. We now 
show that the latter group is an almost-direct product. This 
follows quite directly from Proposition~\ref{Pn(S)^ab} and 
Proposition-Definition~\ref{def_ppd}. 

\begin{lemma}
\label{adp_planar_surfaces}
If $S$ is a connected planar surface, then $P_n(S) \cong 
\pi_1(S \setminus Q_{n-1}) \rtimes P_{n-1}(S)$ is an 
almost-direct product.
\end{lemma}

\begin{proof}
We have $S \setminus Q_{n-1} \simeq S \vee (\D \setminus Q_{n-1})$ 
(a homotopy equivalence stemming from $\D \simeq \D \vee \D$, 
since $S$ is a disk with holes), 
so $\pi_1(S \setminus Q_{n-1}) \cong \pi_1(S) * F_{n-1}$, which implies 
$H_1(S \setminus Q_{n-1}) \cong H_1(S) \oplus \Z^{n-1}$. Then the map 
$H_1(S \setminus Q_{n-1}) \rightarrow P_n(S) ^{\ab}$ identifies, by way of 
Proposition~\ref{Pn(S)^ab}, with the map from $H_1(S) \oplus \Z^{n-1}$ to 
$H_1(S)^n \oplus P_n^{\ab}$ sending $H_1(S)$ isomorphically to the 
first factor, and sending each element $\bar x_i$ in the basis of 
$\Z^{n-1}$ to the class of the usual pure braid generator 
$\overline A_{i,n} \in P_n^{\ab}$. Since $P_n^{\ab}$ is free 
abelian on the elements $\overline A_{ij}$, this map is injective. 
The result follows, using Proposition-Definition~\ref{def_ppd}.
\end{proof}

\begin{remark}\label{not_adp_when_not_planar}
When $S$ is not planar, but not closed, and $n \geq 2$, then 
$P_n(S) \cong \pi_1(S \setminus Q_{n-1}) \rtimes P_{n-1}(S)$ 
is never an almost-direct product. Indeed, we still have 
$\pi_1(S \setminus Q_{n-1}) \cong H_1(S) \oplus \Z^{n-1}$, 
but the map to $P_n(S)^ab \cong H_1(S)^n$ kills $\Z^{n-1}$, 
so it is never injective.
\end{remark}

Finally, we get a generalization of the classical computation of $\gr(P_n)$ recalled in 
Example~\ref{ex:pure braids}, initially done in~\cite{FR, FR88}.
\begin{corollary}
If $S$ is a connected planar surface, then $P_n(S)$ is RTFN. Moreover, the filtration 
$\gamma_*(P_n(S))$ is torsion-free, and $\gr(P_n(S))$ decomposes as an iterated 
semi-direct product of free Lie rings, as follows:
\[
\gr(P_n(S)) \cong \gr(\pi_1(S \setminus Q_{n-1})) \rtimes 
\Bigl(\gr(\pi_1(S \setminus Q_{n-2})) \rtimes \Bigl(\cdots 
\rtimes  \gr(\pi_1(S)) \cdots \Bigr) \Bigr).\]
\end{corollary}

\begin{proof}
Apply the usual results of Falk and Randell (Proposition~\ref{lcs_of_adp} and 
Corollary~\ref{cor:rtfn-ext}) to the almost-direct product decomposition 
of Lemma~\ref{adp_planar_surfaces}, using the fact that the free groups 
$\pi_1(S \setminus Q_m)$ are RTFN (Example~\ref{ex:free_groups_are_RTFN}) 
and that their associated graded are free Lie rings by Magnus' theorem 
(Example~\ref{ex:free}), and thus have no torsion in their lower central series. 
Note that Magnus' theorem also holds for infinitely generated free groups, 
as one sees by either adapting its proof, or deducing it from the case of 
finitely generated groups (namely, if $F[X]$ is free on $X$, any relation 
in $\gr(F[X])$ involves a finite number of generators, corresponding to a 
finite subset $Y \subseteq X$. Then the Lie morphism $\gr(F[X]) \surj \gr(F[Y])$ 
induced by canonical retraction $F[X] \surj F[Y]$ sends this 
relation to itself, which means that is holds in the free Lie ring $\gr(F[Y])$, 
so it is trivial).
\end{proof}
			
\bibliographystyle{abbrv}
\bibliography{Ref_LCS_split}

\end{document}